\documentclass[11pt,a4paper]{article}
\usepackage{ifpdf}
\pdfoutput=1
\usepackage{amsmath}
\usepackage{amssymb}
\usepackage{mathrsfs}
\usepackage{color}
\usepackage{graphicx}
\usepackage[percent]{overpic}
\usepackage{url}
\usepackage{subfig}
\usepackage{enumerate}
\usepackage{overpic}
\usepackage{amsthm}
\usepackage{moreverb}
\usepackage[pdftex,colorlinks,bookmarksopen,bookmarksnumbered,citecolor=red,urlcolor=red]{hyperref}
\usepackage{authblk}

\def\b{{\bm b}}

\def\f{{\bm f}}

\def\n{{\bm n}}

\def\u{{\bm u}}
\def\v{{\bm v}}

\def\0{\boldsymbol{0}}

\def\vbar{\overline{\v}}

\def\dt{\partial_t}

\def\cl {\nonumber \\}
\def\el {\nonumber }

\newtheorem{rem}{Remark}[section]

\newcommand{\bm}[1]{\mbox{\boldmath{$#1$}}}
\def\div{\nabla\cdot}

\usepackage{verbatim}
\usepackage{graphicx}
\usepackage{epstopdf}

\begin{document}
\date{}
\title{A Hybrid Reduced Order Model for nonlinear LES filtering}


\author[1]{Michele Girfoglio\thanks{mgirfogl@sissa.it}}
\author[2]{Annalisa Quaini\thanks{quaini@math.uh.edu}}
\author[1]{Gianluigi Rozza\thanks{grozza@sissa.it}}
\affil[1]{SISSA, International School for Advanced Studies, Mathematics Area, mathLab, via Bonomea, Trieste 265 34136, Italy}
\affil[2]{Department of Mathematics, University of Houston, Houston TX 77204, USA}
\maketitle
\begin{abstract}
We develop a Reduced Order Model (ROM) for a Large Eddy Simulation (LES) approach
that combines a three-step algorithm called Evolve-Filter-Relax (EFR) with a computationally
efficient finite volume method. The main novelty of our ROM  lies
in the use within the EFR algorithm of a nonlinear, deconvolution-based indicator function that identifies 
the regions of the domain where the flow needs regularization.
The ROM we propose is a hybrid projection/data-driven strategy: a classical Proper Orthogonal Decomposition 
Galerkin projection approach for the reconstruction of the velocity and the pressure fields 
and a data-driven reduction method to approximate the indicator function used by the nonlinear differential filter.
This data-driven technique is based on interpolation with Radial Basis Functions. 
We test the performance of our ROM approach on two benchmark problems: 2D and 3D unsteady flow 
past a cylinder at Reynolds number $0 \leq Re \leq 100$. The accuracy of the ROM is assessed 
against results obtained with the full order model for velocity, pressure, indicator function
and time evolution of the aerodynamics coefficients.
\end{abstract}
\vspace*{0.5cm}

\section{Introduction}\label{sec:intro}

Although increasing computational power has become available recently, the resolution of systems of 
parametric Partial Differential Equations 
using classical discretization methods (e.g., finite element or finite volume methods, hereinafter called Full Order Models) is still unfeasible in several situations
where one needs to evaluate a certain output of interest for a large number of input configurations. This is the case of, e.g., real-time control problems, optimization problems, and uncertainty quantification. 
In this context, Reduced Order Models (ROMs) have been proposed as an efficient tool to significantly reduce the computational cost required by classical Full Order Models (FOMs).

The basic ROM framework consists of two steps. The first one is the so-called \emph{offline} phase, when a database
of several solutions is collected by solving the original FOM for selected parameter values. 
The second step is the \emph{online} phase, during which
the information obtained in the offline phase is used to quickly compute the solution for newly specified values of the parameters. 
The assumption underlying this framework is that the solution of a parametrized PDE (or a system of such equations)
lies on a low-dimensional manifold, which can be
approximated by a subspace spanned by the set of basis functions found in the offline phase. 
For a comprehensive review on ROMs, the reader is referred to, e.g., \cite{hesthaven2015certified, quarteroniRB2016, bennerParSys, Benner2015, Bader2016, ModelOrderReduction}.

In recent years, there has been a growing interest in constructing stable, accurate, and computationally efficient ROMs 
for the numerical simulation of higher Reynolds number flows. Traditional ROMs fail to maintain their promise
of reduced computational costs when the Kolmogorov $n$-width of
the solution manifold associated with the FOM is large, as is the case in convection-dominated flows. 
Indeed, if one choose to retain a large number of modes in order for the ROM to be accurate, then
the computational efficiency suffers. If the number of modes is otherwise kept low, a severe loss of information
hinders the accurate reconstruction of the solution. 
In fact, projection based ROMs of turbulent flows are affected by energy stability problems
related to the fact that proper orthogonal decomposition (POD) retains the modes biased toward large,
high-energy scales, while the turbulent kinetic energy is dissipated at level of
the small turbulent scales.  A possible way to tackle this challenging problem 
is to introduce dissipation via a closure model \cite{wang_turb, Aubry1988}. In \cite{couplet_sagaut_basdevant_2003},
it was shown theoretically and numerically that modes have a similar energy transfer mechanism to Fourier modes. 
Therefore, the use of Large Eddy Simulation (LES) could be beneficial. Following this idea, 
in this paper we develop a ROM for a LES filtering approach for flows at moderate
Reynolds numbers (a few hundreds or a few thousands). 

We focus on a nonlinear variant of the Evolve-Filter-Relax (EFR) algorithm \cite{Boyd1998283, Fischer2001265, Dunca2005,layton_CMAME,BQV}, which describes the effect of the small scales (i.e., the turbulent
scales that are not resolved by the discretization mesh) by a set of equations to be added to the discrete 
Navier-Stokes equations. This extra problem can be interpreted as a differential low-pass filter.
Modularity is an appealing feature of this approach since it can be implemented
without a major modification of a legacy Navier-Stokes solver.
A key role in the EFR algorithm is played by a nonlinear indicator function that identifies
the regions of the domain where the flow needs regularization (i.e., where and how much artificial viscosity is needed)
\cite{abigail_CMAME}. 
At the full order level, the EFR algorithm and its variant without relaxation (called EF)
have been extensively studied within a Finite Element framework. Instead, we choose to apply 
a computationally efficient Finite Volume (FV) method for the space discretization \cite{Girfoglio2019, Girfoglio2021}. 
The motivation for this choice is that many commercial codes are based on FV methods, Thus, a FV-based ROM 
would be appealing for a larger group of CFD practitioners. 
 
The large majority of the regularized ROMs adds the filtering step only at the reduced
order level, i.e.~the snapshots are obtained by Direct Numerical Simulation. See, e.g., to \cite{Xie2018_2,Wells2017,Gunzburger2019}. 
We proposed a different approach in \cite{Girfoglio_JCP, Girfoglio_ROM_Fluids, Strazzullo2021} 
for the EF/EFR algorithm: 
we apply the filter step at both the full and reduced order level, i.e.~we 
generate snapshots data with under-refined meshes. 
Such an approach provides a ROM that is fully consistent with the FOM since the same mathematical framework 
is used during both the \emph{offline} and \emph{online} stage. 
We use the POD basis related to the evolve velocity to approximate the filtered velocity 
and we compute the reduced pressure field with a Poisson Pressure Equation (PPE) method 
\cite{Stabile2018, akhtar2009stability}. The main difference between this work and  \cite{Girfoglio_JCP, Girfoglio_ROM_Fluids, Strazzullo2021}
lies in the indicator function. Indeed,
the EF or EFR algorithms in \cite{Girfoglio_JCP, Girfoglio_ROM_Fluids, Strazzullo2021} 
adopt a linear indicator function. While this was a necessary simplification that allowed us to focus on other challenges 
posed by a ROM differential filter, it is known that a linear indicator function is insufficiently selective 
as it applies the same amount of artificial viscosity everywhere in the domain. Thus, in this paper we
extend our approach to a nonlinear indicator function.
We propose a hybrid projection/data-driven ROM that draws inspiration from the work in \cite{Hijazi2020},
where a ROM framework for the Reynolds-averaged Navier Stokes (RANS) equations is developed. We exploit a traditional projection method
for the computation of the reduced velocity and pressure fields, while we use a
data-driven technique to compute the reduced coefficients of the indicator function field.
This data-driven technique consists in an interpolation process with Radial Basis Functions \cite{Lazzaro2002}. 

We choose the above-mentioned hybrid approach because it is accurate and partially non-intrusive 
The  same does not hold true for two well known alternatives available in the literature. The first, easy alternative would be to use the same 
set of reduced coefficients for velocity, pressure, and indicator function. This approach has been considered 
for RANS in \cite{Lorenzi2016} and it was shown in \cite{Hijazi2020} to provide less accurate results than a hybrid procedure.
It is reasonable to assume this lack of accuracy could be observed for LES too. 
The second option would be to use a EIM/DEIM technique \cite{barrault04:_empir_inter_method,Chaturantabut2010}
for the reconstruction of the indicator function field but its implementation is problem-dependent and intrusive.
Instead, our hybrid procedure provides a unique computational pipeline for the development of efficient ROMs 
for flows at higher Reynolds number, no matter if modeled by LES or RANS. 
We test our approach on two benchmarks: 
2D \cite{turek1996, John2004} and 3D \cite{turek1996} flow past a cylinder with time-dependent Reynolds number $0 \leq Re(t) \leq 100$. 
We limit our investigation to the ROM reconstruction of the time evolution of the
system. Parametric studies, in particular related to key model parameters such as 
filtering radius and relaxation parameter, will be addressed in a future work. 
For both tests, we compare the evolution of velocity, pressure, and indicator function
with the corresponding FOM quantities. Additionally, we show the accuracy of our approach in the time
reconstruction of the the lift and drag coefficients.

The work is organized as follows. Sec.~\ref{sec:FOM} describes 
the full order model and the numerical method we use for it. 
Sec.~\ref{sec:ROM} presents the ingredients of the reduced order model. 
The numerical examples are reported in Sec. \ref{sec:results}. 
Sec. \ref{sec:conclusions} provides conclusions and perspectives.

\section{The full order model}\label{sec:FOM}


We consider the motion of an incompressible viscous fluid in a fixed domain $\Omega \subset \mathbb{R}^D$ with $D = 2, 3$ over a time 
interval of interest ($t_0$, $T$) $\subset \mathbb{R}^+$. The flow is described by the incompressible Navier-Stokes equations (NSE):
\begin{align}
\rho\,\dt \u + \rho\,\div \left(\u \otimes \u\right) - 2\mu \Delta\u + \nabla p & = \f\quad \mbox{ in }\Omega \times (t_0,T),\label{eq:ns-mom}\\
\div \u & = 0\quad\, \mbox{ in }\Omega \times(t_0,T),\label{eq:ns-mass}
\end{align}
where $\rho$ is the fluid density, $\mu$ is the dynamic viscosity, $\u$ is velocity, and $p$ is the pressure. 
Problem \eqref{eq:ns-mom}-\eqref{eq:ns-mass} is endowed with 
suitable boundary conditions
\begin{align}
&\u = \u_D  &&\mbox{on } \partial\Omega_D \times(t_0,T), \label{eq:bc-d} \\ 
&(2\mu \nabla\u - p\mathbf{I})\n = \0  &&\mbox{on } \partial\Omega_N \times(t_0,T), \label{eq:bc-n} 
\end{align}
and the initial data $\u= \u_0$ in $\Omega \times\{t_0\}$. Here $\overline{\partial\Omega_D}\cup\overline{\partial\Omega_N}=\overline{\partial\Omega}$ and $\partial\Omega_D \cap\partial\Omega_N=\emptyset$. In addition, $\u_D$ and $\u_0$ are given.

In order to characterize the flow regime under consideration, we define the Reynolds number as
\begin{equation}\label{eq:re}
Re = \frac{U L}{\nu},
\end{equation}
where $\nu=\mu/\rho$ is the \emph{kinematic} viscosity of the fluid, and $U$ and $L$ are characteristic macroscopic velocity and length, respectively.
We are interested in simulations flows at moderately large Reynolds numbers, for which 
flow disturbances cannot be neglected and Reynolds-averaged Navier-Stokes (RANS) models \cite{pope} 
are inaccurate. For this reason, we choose to work with a Large Eddy Simulation (LES) approach that is described
next. 

\subsection{The Evolve-Filter-Relax algorithm}\label{sec:FOM_2}


Our LES approach is a modular algorithm, called Evolve-Filter-Relax (EFR) \cite{layton_CMAME}, 
that adds a differential filter to the Navier-Stokes equations \eqref{eq:ns-mom}-\eqref{eq:ns-mass}.
This algorithm comes from the decoupling of the time-discrete Leray model \cite{Leray1934}. For the implementation the EFR algorithm, we chose the C++ finite volume library OpenFOAM\textsuperscript{\textregistered} \cite{Weller1998}.

Let $\Delta t \in \mathbb{R}$, $t^n = t_0 + n \Delta t$, with $n = 0, ..., N_T$ and $T = t_0 + N_T \Delta t$. 
We denote by $y^n$ the approximation of a generic quantity $y$ at the time $t^n$. 
The EFR algorithm  reads as follows: given the velocities $\u^{n-1}$ and $\u^{n}$, at $t^{n+1}$:
\begin{enumerate}[i)]
\item \textit{Evolve}: find intermediate velocity and pressure $(\v^{n+1},p^{n+1})$ such that
\begin{align}
&\rho\, \dfrac{3}{2\Delta t}\v^{n+1} + \rho\, \div \left(\v^n \otimes \v^{n+1}\right) - 2\mu\Delta\v^{n+1} +\nabla p^{n+1} = \b^{n+1},\label{eq:evolve-1.1}\\
& \div \v^{n+1} = 0\label{eq:evolve-1.2},
\end{align}
with boundary conditions
\begin{align}
&\v^{n+1} = \u_D^{n+1}  \quad \mbox{on } \partial\Omega_D \times(t_0,T), \label{eq:bc-d_1} \\ 
&(2\mu \nabla\v^{n+1} - p^{n+1}\mathbf{I})\n = \0  \quad\mbox{on } \partial\Omega_N \times(t_0,T), \label{eq:bc-n_1} 
\end{align}
and initial condition $\v^0 = \u_0$ in $\Omega \times\{t_0\}$. In eq.~\eqref{eq:evolve-1.1}, we set $\b^{n+1} = \rho(4\u^n - \u^{n-1})/(2\Delta t)$.
Notice that this step corresponds to a time discretization with Backward Differentiation Formula of order 2 (BDF2)
of problem \eqref{eq:ns-mom}-\eqref{eq:ns-mass}.

\begin{rem}
We adopt a first order extrapolation for the convective velocity although a BDF2 scheme 
is used for the time discretization of problem \eqref{eq:ns-mom}-\eqref{eq:ns-mass}. 
This is what the NSE solvers in OpenFOAM do, so discretization 
\eqref{eq:evolve-1.1}-\eqref{eq:evolve-1.2} would make it a fair comparison 
between EFR and NSE algorithms. 
\end{rem}

\item \textit{Filter}: find filtered velocity $\vbar^{n+1}$ such that
\begin{align}
&-\alpha^2\div \left( a(\v^{n+1}) \nabla\vbar^{n+1}\right) +\vbar^{n+1} = \v^{n+1}, \label{eq:filter}
\end{align}
with boundary conditions
\begin{align}
& \vbar^{n+1} = \u_D^{n+1}  \quad \mbox{on } \partial\Omega_D \times(t_0,T), \label{eq:bc-filter-d_1} \\ 
& \nabla\vbar^{n+1}\n = \0  \quad\mbox{on } \partial\Omega_N \times(t_0,T).\label{eq:bc-filter-n_1}
\end{align}
In eq. \eqref{eq:filter}, $\vbar$ is the \emph{filtered velocity} and $\alpha$ can be 
interpreted as the \emph{filtering radius} (that is, the radius of the neighborhood where 
the filter extracts information from the unresolved scales). \emph{Indicator function} $a(\cdot)$ is such that:
\begin{align*}
a(\v)\simeq 0 & \mbox{ where the velocity $\v$ does not need regularization;}\\
a(\v)\simeq 1 & \mbox{ where the velocity $\v$ does need regularization.}
\end{align*}
Different choices of $a(\cdot)$ have been proposed and compared in \cite{Borggaard2009,layton_CMAME,O-hunt1988,Vreman2004,Bowers2012}. 
We choose indicator function:
\begin{align}
a(\v) = \left|  \v - F(\v) \right|, \label{eq:a_D0_a_D1}
\end{align}
where we take $F(\v)$ to be the linear Helmholtz filter operator, i.e.~ $F(\v) = \tilde{\u}$ with: 
\begin{equation}\label{eq:vtilde}
\tilde{\u}^{n+1} - \alpha^2  \Delta \tilde{\u}^{n+1} = \v^{n+1}. 
\end{equation}
For more details on this indicator function, we refer to \cite{BQV}.
\item \textit{Relax}: set 
\begin{align}
\u^{n+1}&=(1-\chi)\v^{n+1} + \chi\vbar^{n+1}, \label{eq:relax-1}
\end{align}
where $\chi\in(0,1]$ is a relaxation parameter.
\end{enumerate}
We consider $\u^{n+1}$ the approximation of the velocity $t^{n+1}$. It is possible to show that the above
EFR algorithm is equivalent to a generic viscosity model in LES \cite{Olshanskii2013}.

\begin{rem}
In this paper, we consider a simplified filter problem with respect to our previous 
work \cite{BQV,Girfoglio2019,Girfoglio_JCP,Girfoglio2021}, where we forced the filtered
velocity $\vbar^{n+1}$ to be solenoidal. We are releasing this constraint as it leads to 
a substantial simplification and computational time savings since there is one less variable
(i.e., the Lagrange multiplier to enforce the incompressibility constraint). 
As noted in \cite{Ervin2010}, the incompressibility is exactly  preserved by the simplified differential filter 
\eqref{eq:filter} only for periodic conditions. 
Thus, in our case the end-of-step velocity $\u^{n+1}$ does not strictly satisfy mass conservation. 
However, we will show in Sec. \ref{sec:2D} that at discrete level the mass conservation error is very low.
\end{rem}

\begin{rem}
The EFR method has an appealing advantage over other LES models: 
it is modular, i.e.~it adds a differential problem to the Navier-Stokes problem
instead of extra terms in the Navier-Stokes equations themselves (like, e.g., the popular
variational multiscale approach \cite{BAZILEVS2007173}). Thus, thanks to the EF method
anybody with a Navier-Stokes solver could simulate higher Reynolds number
flows without major modifications to the software core. 
\end{rem}

Finally, we note that while we consider homogeneous Neumann boundary conditions
non-homogeneous boundary condition can of course be handled. See \cite{BQV} to learn more about this.




\subsection{Space discrete problem: a Finite Volume approximation}\label{sec:FOM_3}
For the space discretization of problems \eqref{eq:evolve-1.1}-\eqref{eq:bc-n_1} and \eqref{eq:filter}-\eqref{eq:vtilde},
we adopt a Finite Volume (FV) method. We partition the computational domain $\Omega$ into cells or control volumes $\Omega_i$, with $i = 1, \dots, N_{c}$, where $N_{c}$ is the total number of cells in the mesh. 
Let  \textbf{A}$_j$ be the surface vector of each face of the control volume, 
with $j = 1, \dots, M$. 

The fully discretized form of problem \eqref{eq:evolve-1.1}-\eqref{eq:evolve-1.2} reads: Find
$(\v_i^{n+1},p^{n+1}_i)$ such that
\begin{align}
&\rho\, \frac{3}{2\Delta t}\, \v^{n+1}_i + \rho\, \sum_j^{} \varphi_j^n \v^{n+1}_{i,j} - 2\mu \sum_j^{} (\nabla\v^{n+1}_i)_j \cdot \textbf{A}_j + \sum_j^{} p^{n+1}_{i,j} \textbf{A}_j  = {\bm b}^{n+1}_i \label{eq:disc_evolve1} \\ 
&\sum_j^{} (\nabla p^{n+1})_j \cdot \textbf{A}_j = \sum_j^{} (\textbf{H}(\v_i^{n+1}))_j \cdot \textbf{A}_j, \label{eq:disc_evolve2}
\end{align}
where:
\begin{align}
\textbf{H}(\v^{n+1}_i) = -\rho \sum_j^{} \varphi_j^n \v^{n+1}_{i,j} + 2\mu \sum_j^{} (\nabla\v^{n+1}_i)_j \cdot \textbf{A}_j + {\bm b}^{n+1}_i \quad \text{with} \quad \varphi_j^n = \v^{n}_j \cdot \textbf{A}_j. \label{eq:H}
\end{align}
In \eqref{eq:disc_evolve1}-\eqref{eq:H}, $\v^{n+1}_i$ and ${\bm b}^{n+1}_i$ denote the average velocity and source term in control volume $\Omega_i$, respectively. Moreover, we denote with $\v^{n+1}_{i,j}$ and $p^{n+1}_{i,j}$ the velocity and pressure
associated to the centroid of face $j$ normalized by the volume of $\Omega_i$.
For the solution of the linear system associated with \eqref{eq:disc_evolve1}-\eqref{eq:disc_evolve2} we used the PISO algorithm \cite{PISO}. 
The advantage of this algorithm is the decoupling of the computation of the pressure from the computation of the velocity, which 
results in low computational costs.

Next, we discretize filter problem \eqref{eq:filter}. We obtain: 
\begin{align}
\vbar^{n+1}_i - \alpha^2 \sum_j^{} a_j^{n+1}(\nabla\vbar^{n+1}_i)_j \cdot \textbf{A}_j = \v^{n+1}_i, \label{eq:filter1_disc}
\end{align}
where $\vbar^{n+1}_i$ is the average value of $\vbar^{n+1}$ in control volume $\Omega_i$ and $a_j^{n+1} = a(\v_j^{n+1})$. 
To compute $a_j^{n+1}$, we need to solve the Helmholtz filter problem \eqref{eq:vtilde}. Once discretized, it reads: Find 
the average value of $\tilde{\v} ^{n+1}$ in $\Omega_i$, i.e.~$\tilde{\v} ^{n+1}_i$, such that:
\begin{align}
\tilde{\v} ^{n+1}_i - \alpha^2 \sum_j^{} (\nabla \tilde{\v} ^{n+1}_i)_j \cdot \textbf{A}_j = \v^{n+1}_i. \label{eq:filter2_disc}
\end{align}
Obviously, problems \eqref{eq:filter1_disc} and \eqref{eq:filter2_disc} are easier to solve than problem 
\eqref{eq:disc_evolve1}-\eqref{eq:disc_evolve2}, i.e.~the filter problem is computationally much less 
demanding than the Navier-Stokes problem. 




\section{The reduced order model}\label{sec:ROM}
The Reduced Order Model (ROM) we propose can be seen as an extension to 
a LES framework of the model introduced in \cite{Hijazi2020} for RANS. 
The key idea is the following: we use a Galerkin projection method related to compute 
the reduced velocity and pressure fields, while we use an interpolation procedure based on 
Radial Basis Functions (RBF) for the computation of the reduced coefficients of the indicator function. 
We call this hybrid approach \emph{data-driven POD-Galerkin ROM}.
In Sec~\ref{sec:ROM_1} we describe the details of our approach
and in Sec.~\ref{sec:ROM_2} we present the strategy we choose 
for pressure stabilization at reduced order level. 

The ROM computations have been carried out using ITHACA-FV \cite{RoSta17}, an in-house open source C++ library. 

\subsection{Our data-driven POD-Galerkin method}\label{sec:ROM_1}

We approximate velocity fields $\v$ and $\vbar$, pressure field $p$, and indicator function $a$
as linear combinations of the dominant modes (basis functions), assumed to depend 
on space variables only, multiplied by scalar coefficients that depend only on time:
\begin{align}
\v \approx \v_r = \sum_{i=1}^{N_{v_r}} \beta_i(t) \bm{\varphi}_i(\bm{x}), \quad 
p \approx p_r = \sum_{i=1}^{N_{q_r}} \gamma_i(t) \psi_i(\bm{x}), \label{eq:ROM_1}\\
\vbar \approx \vbar_r = \sum_{i=1}^{N_{{\overline{v}}_r}} \overline{\beta_i}(t) \bm{\varphi}_i(\bm{x}), \quad
a \approx a_r = \sum_{i=1}^{N_{a_r}} \delta_i(t) \eta_i(\bm{x}). \label{eq:ROM_2}
\end{align}
In \eqref{eq:ROM_1}-\eqref{eq:ROM_2}, $N_{\Phi_r}$ denotes the cardinality of a reduced basis for the space 
$\Phi$ belongs to. 

Using \eqref{eq:ROM_1} to approximate $\v^{n+1}$ and $q^{n+1}$ in \eqref{eq:evolve-1.1}-\eqref{eq:evolve-1.2},
we obtain
\begin{align}
\rho\, \frac{3}{2\Delta t}\, \v_r^{n+1} + \rho\, \div \left(\v_r^{n} \otimes \v_r^{n+1}\right) - 2\mu\Delta\v_r^{n+1} +\nabla p_r^{n+1} = \b_r^{n+1} , \label{eq:red-1.1} \\
 \div \v_r^{n+1} = 0. \label{eq:red-1.2}
\end{align}
Then, using \eqref{eq:ROM_2} to approximate $\vbar^{n+1}$ and $a(\v^{n+1})$ in \eqref{eq:filter} 
we get:
\begin{align}
\vbar_r^{n+1}  - \alpha^2 \div (a_r^{n+1} \nabla \vbar_r^{n+1}) = \v_r^{n+1}. \label{eq:red-2.1} 
\end{align}

\begin{rem}\label{rem:rem1}
As mentioned above, we use a data-driven interpolation for the approximation of indicator function $a$ 
defined in \eqref{eq:a_D0_a_D1}. 
For this reason, we do not need to compute a reduced order approximation of $\tilde{\v}$.
\end{rem}

\begin{rem}\label{rem:rem2}
We use the reduced basis $\bm{\varphi}_j$ associated to $\v$ also for the approximation of $\vbar$ in \eqref{eq:ROM_2}.
Thus, in the ROM velocity $\vbar$ is divergence free, although the same is not true in the FOM.
\end{rem}

In the literature, one can find several techniques to generate the reduced basis spaces, e.g.~Proper Orthogonal Decomposition (POD), the Proper Generalized Decomposition and the Reduced Basis with a greedy sampling strategy.
See, e.g., \cite{Rozza2008, ChinestaEnc2017, Kalashnikova_ROMcomprohtua, quarteroniRB2016, Chinesta2011, Dumon20111387, Huerta2020, ModelOrderReduction}. 
We choose to find the reduced basis by using the method of snapshots.
To this purpose, we solve the FOM described in Sec.~\ref{sec:FOM}  
for each time $t^k \in \{t^1, \dots, t^{N_s}\} \subset (t_0, T]$. 
The snapshots matrices are obtained from the full-order snapshots: 
\begin{align}\label{eq:space}
\bm{\mathcal{S}}_{{{\Phi}}} = [{{\Phi}}(t^1), \dots, {{\Phi}}(t^{N_s})] \in \mathbb{R}^{N_{\Phi_h} \times N_s} \quad
\text{for} \quad {{\Phi}} = \{\v, p, a\}, 
\end{align}
where the subscript $h$ denotes a solution computed with the FOM and $N_{\Phi_h}$
is the dimension of the space $\Phi$ belong to in the FOM. Note that ${\Phi}$ could be either a scalar or
a vector field. 
The POD problem consists in finding, for each value of the dimension of the POD space $N_{POD} = 1, \dots, N_s$, the scalar coefficients 
$c_1^1, \dots, c_1^{N_s}, \dots, c_{N_s}^1, \dots, c_{N_s}^{N_s}$ and functions ${\bm{\zeta}}_1, \dots, {\bm{\zeta}}_{N_s}$, that minimize the error between the snapshots and their projection onto the POD basis. In the $L^2$-norm, we have
\begin{align}
E_{N_{POD}} = \text{arg min} \sum_{i=1}^{N_s} ||{{\Phi}_i} - \sum_{k=1}^{N_{POD}} c_i^k {\bm{\zeta}}_k || \quad \forall N_{POD} = 1, \dots, N_s    \cl
\text{with} \quad ({\bm{\zeta}}_i, {\bm{\zeta}}_j)_{L_2(\Omega)} = \delta_{i,j} \quad \forall i,j = 1, \dots, N_s. \label{eq:min_prob}
\end{align}

It can be shown \cite{Kunisch2002492} that eq.~\eqref{eq:min_prob} is equivalent to the following eigenvalue problem
\begin{align}
\bm{\mathcal{C}}^{{\Phi}} \bm{Q}^{{\Phi}} &= \bm{Q}^{{\Phi}} \bm{\Lambda}^{{\Phi}}, \label{eq:eigen_prob} \\
\mathcal{C}_{ij}^\Phi &= ({\Phi}(t^i), {\Phi}(t^j))_{L_2(\Omega)} \quad \text{for} \quad i,j = 1, \dots, N_s,
\end{align}
where $\bm{\mathcal{C}}^{{\Phi}}$ is the correlation matrix computed from the snapshot matrix $\bm{\mathcal{S}}_{{{\Phi}}}$, $\bm{Q}^{{\Phi}}$ is the matrix of eigenvectors and $\bm{\Lambda}^{{\Phi}}$ is a diagonal matrix whose diagonal entries are the
eigenvalues of $\bm{\mathcal{C}}^{{\Phi}}$. 
Then, the basis functions are obtained as follows:
\begin{align}\label{eq:basis_func}
\bm{\zeta}_i = \dfrac{1}{N_s \Lambda_i^\Phi} \sum_{j=1}^{N_s} {\Phi}_j Q_{ij}^\Phi.
\end{align}
The POD modes resulting from the aforementioned methodology are:
\begin{align}\label{eq:spaces}
L_\Phi = [\bm{\zeta}_1, \dots, \bm{\zeta}_{N_{\Phi_r}}] \in \mathbb{R}^{N_{\Phi_h} \times N_{\Phi_r}},
\end{align}
where $N_{\Phi_r} < N_s$ are chosen according to the eigenvalue decay. 
The reduced order model can be obtained through a Galerkin projection of the governing equations onto the POD spaces. 

Let
\begin{align}
&M_{r_{ij}} = (\bm{\varphi}_i, \bm{\varphi}_j)_{L_2(\Omega)}, 
\quad A_{r_{ij}} = (\bm{\varphi}_i, \Delta \bm{\varphi}_j)_{L_2(\Omega)}, \quad B_{r_{ij}} = (\bm{\varphi}_i, \nabla \psi_j)_{L_2(\Omega)}, \label{eq:matrices_evolve1} \\
&P_{r_{ij}} = (\psi_i, \nabla \cdot \bm{\varphi}_j)_{L_2(\Omega)}, \quad G_{r_{ijk}} = (\bm{\varphi_i}, \nabla \cdot (\bm{\varphi_j} \otimes \bm{\varphi}_k))_{L_2(\Omega)}, 
\label{eq:matrices_evolve2}
\end{align}
where $\bm{\varphi}_i$ and $\psi_i$ are the basis functions in \eqref{eq:ROM_1}. The reduced algebraic system at time $t^{n+1}$
for problem \eqref{eq:red-1.1}-\eqref{eq:red-1.2} is: 
\begin{multline}
\rho\, \frac{3}{2\Delta t}\ \bm{M}_r \bm{\beta}^{n+1} + \rho (\bm{\beta}^{n})^T  \bm{G}_r \bm{\beta}^{n+1} - 2\mu \bm{A}_r \bm{\beta}^{n+1} + \bm{B}_r \bm{\gamma}^{n+1} = \\\dfrac{\rho}{2\Delta t} \bm{M}_r \left(\left(1-\chi\right)\left(4\bm{\beta}^{n} - \bm{\beta}^{n-1}\right) +  \chi\left(4\overline{\bm{\beta}}^{n} - \overline{\bm{\beta}}^{n-1}\right) \right), \label{eq:reduced_1} 
\end{multline}
\begin{align}
\bm{P}_r \bm{\beta}^{n+1} = 0, \label{eq:reduced_2}
\end{align}
where vectors $\bm{\beta}^{n+1}$ and $\bm{\gamma}^{n+1}$ contain the values of coefficients $\beta_i$ and $\gamma_i$ 
in \eqref{eq:ROM_1} at time $t^{n+1}$.

Next, let 
\begin{align}
&{A}_{r_{ijk}} = (\bm{{\varphi}}_i, \nabla \cdot \eta_j \nabla \bm{{\varphi}}_k)_{L_2(\Omega)}, \label{eq:matrices_filter1} 
\end{align}
where $\eta_i$ are the basis functions in \eqref{eq:ROM_2}. 
The reduced algebraic system at time $t^{n+1}$
for problem \eqref{eq:red-2.1} is
\begin{align}
&\bm{{{M}}}_r \bm{\overline{\beta}}^{n+1}  - \alpha^2 (\bm{\delta}^{n+1})^T\bm{{A}}_r \bm{\overline{\beta}}^{n+1}  = \bm{{{M}}}_r \bm{\beta}^{n+1} \label{eq:reduced2_1}, 
\end{align}
where vectors $\overline{\bm{\beta}}^{n+1}$ and $\bm{\delta}^{n+1}$ contain the values of coefficients $\overline{\beta}_i$ and $\delta_i$ in \eqref{eq:ROM_2} at time $t^{n+1}$.

The coefficients $\delta_i(t)$ in \eqref{eq:ROM_2} are computed with a data-driven approach that uses interpolation with Radial Basis Functions (RBF) \cite{Lazzaro2002}. The interpolation procedure is carried out for each mode separately. Let
$\eta_i$ be the indicator function mode under consideration. 
Function $G_i(t)$ that interpolates $\delta_i(t)$ using RBF functions can be written as: 
\begin{align}
G_i(t) = \sum_{j=1}^{N_s} w_{i,j} \zeta_{i,j}(|t - t^j), \quad \text{for} \quad i= 1, 2, \dots, N_{a_r}, ~t \in (t_0, T]\label{eq:RBF1}
\end{align}
where $w_{i,j}$ are suitable weights and $\zeta_{i,j}$ are the radial basis functions, 
which are chosen to be Gaussian functions. We observe that $\zeta_{i,j}$ is centered at time $t^j$. 
In order to compute the weights $w_{i,j}$, we use the fact that 
$G_i$ has to interpolate $\delta_i$ at time nodes 
$t^k$, i.e.:
\begin{align}
G_i(t^k) = \delta_i(t^k) \quad \text{for} \quad k = 1, 2, \dots, N_s. \el 
\end{align}
The coefficients $\delta_i(t^k)$ are obtained from projecting the $k$-th snapshot
onto the $i$-th mode: 
\begin{align}
\delta_i(t^k) = (a(t^k), \eta_i)_{L^2(\Omega)},\el 
\end{align}
$a(t^k)$ being the $k$-th column of the snapshot matrix $\bm{\mathcal{S}}_a$ \eqref{eq:space}.
Then, we have:
\begin{align}
G_i(t^k) = \sum_{j=1}^{N_s} w_{i,j} \zeta_{i,j} \left(|t^k - t^j|\right) = \delta_i(t^k). \el 
\end{align}
which can be written as a linear system:
\begin{align}
\bm{A}_i^\zeta \bm{w}_i= \bm{Y}_i, \quad \text{with } (\bm{A}_i^\zeta)_{kj} = \zeta_{i,j}(|t^k - t^j|). \label{eq:RBF5}
\end{align}
System \eqref{eq:RBF5} is solved \emph{offline} to get the weights $\bm{w}_i$. During the \emph{online} phase,
for every new time instant $t^*$ we compute ${\delta}_i(t^*)$ given by: 
\begin{align}
{\delta}_i(t^*) \approx G_i(t^*) = \sum_{j=1}^{N_s} w_{i,j} \zeta_{i,j} \left(|t^* - t^j|\right). \el 
\end{align}

 The initial conditions for the ROM algebraic system \eqref{eq:reduced_1}-\eqref{eq:reduced_2}, \eqref{eq:reduced2_1}
are obtained with a Galerkin projection of the initial full order conditions onto the POD basis spaces:
\begin{align}
{\beta^0}_i = (\v(\bm{x},t_0), \bm{\varphi}_i)_{L_2(\Omega)}, \cl
{\overline{\beta}^0}_i = (\vbar(\bm{x},t_0), \bm{\varphi}_i)_{L_2(\Omega)}. \el
\end{align}

Finally, we use the lifting function method \cite{Girfoglio_JCP} to account for 
non-homogeneous Dirichlet boundary conditions. 
The velocity snapshots are modified according to:
\begin{align}
\v'_h = \v_h - \sum_{j=1}^{N_{BC}} u_{{BC}_j}(t) \bm{\chi}_j(\bm{x}), \cl
\vbar'_h = \vbar_h - \sum_{j=1}^{N_{BC}} u_{{BC}_j}(t)\bm{\chi}_j(\bm{x}), \el
\end{align}
where $N_{BC}$ is the number of non-homogeneous Dirichlet boundary conditions, $\bm{\chi} (\bm{x})$ 
are the divergence free control functions that satisfy the boundary conditions, and $u_{{BC}_j}$ are suitable temporal coefficients. 
The POD is applied to the snapshots satisfying the  homogeneous boundary conditions 
and then the boundary value is added back:
\begin{align}
\v_r= \sum_{j=1}^{N_{BC}} u_{{BC}_j}(t)\bm{\chi}_j(\bm{x})  + \sum_{i=1}^{N_{v_r}} \beta_i(t) \bm{\varphi}_i(\bm{x}), \cl
\vbar_r= \sum_{j=1}^{N_{BC}} u_{{BC}_j}(t)\bm{\chi}_j(\bm{x})  + \sum_{i=1}^{N_{\overline{v}_r}} \overline{\beta_i}(t) \bm{\varphi}_i(\bm{x}). \el
\end{align}

\subsection{Pressure field reconstruction and stability}\label{sec:ROM_2}
In order to obtain a stable and accurate reconstruction of the pressure field at the reduced level, we choose to adopt the 
Poisson pressure equation (PPE) method used, e.g., in \cite{Stabile2018, Girfoglio_JCP, Girfoglio_ROM_Fluids}. 
To obtain the Poisson pressure equation, we take the divergence of eq.~\eqref{eq:evolve-1.1} and account for divergence free condition~\eqref{eq:evolve-1.2}:
\begin{align}
&\Delta p^{n+1} = -\rho\, \nabla \cdot \left(\div \left(\v^n \otimes \v^{n+1}\right)\right) + \div \bm{b}^{n+1}, \label{eq:system_def1_2} 
\end{align}
with boundary conditions \eqref{eq:bc-d_1} and:
\begin{align}
& \partial_{{n}} p^{n+1} = -2 \mu \bm{n} \cdot \left(\nabla \times \nabla \times \v^{n+1} \right) - \bm{n} \cdot \left(\rho\dfrac{3}{2\Delta t}\v^{n+1} - \bm{b}^{n+1}\right) \quad \mbox{on } \partial\Omega_N \times(t_0,T), \label{eq:system_def2_2}
\end{align}
where $\partial_{{n}}$ denotes the derivative with respect to the normal vector \bm{n}. In eq.~\eqref{eq:system_def1_2}, we retain the term $\div \bm{b}^{n+1}$ because at full order level the filtered velocity is not divergence free. For further details about the derivation of non-homogeneous Neumann conditions for the pressure field, we refer the reader to \cite{Orszag1986, JOHNSTON2004221}.

By using \eqref{eq:ROM_1} for the approximation of $\v^{n+1}$ and $p^{n+1}$ in \eqref{eq:system_def1_2},
we obtain
\begin{align}
&\Delta p^{n+1}_r = -\rho\, \nabla \cdot \left(\div \left(\v^n_r \otimes \v^{n+1}_r\right)\right). \label{eq:system_def1_2r}
\end{align}
The term $\div \bm{b}_r^{n+1}$ vanishes in eq.~\eqref{eq:system_def1_2r} 
because at reduced order level the filtered velocity is divergence free as explained in Remark \ref{rem:rem2}.

The matrix form of eq.~\eqref{eq:system_def1_2r} reads:
\begin{align}
&\bm{D}_r \bm{\gamma}^{n+1} + \rho \left({\bm{\beta}}^{n}\right)^T\bm{J}_r \bm{\beta}^{n+1} - 2\mu \bm{N}_r\bm{\beta}^{n+1} - 
\dfrac{3 \rho }{2\Delta t}{\boldsymbol{F}}_{r} \boldsymbol{\beta}^{n+1}\cl
& \quad \quad = \dfrac{\rho }{2\Delta t}{\boldsymbol{F}}_{r}\left( \left(1-\chi\right)\left(4{\boldsymbol{\beta}}^{n} - {\boldsymbol{\beta}}^{n-1}\right) + \chi\left(4\overline{\boldsymbol{\beta}}^{n} - \overline{\boldsymbol{\beta}}^{n-1}\right)\right), \label{eq:reduced_q}
\end{align}
where 
\begin{align}
&D_{r_{ij}} = (\nabla \psi_i, \nabla \psi_j)_{L_2(\Omega)}, \quad N_{r_{ij}} = (\bm{n} \times \nabla \psi_i, \nabla \times \bm{\varphi}_j) _{L_2(\partial \Omega)}, \label{eq:pressure_matrices1} \\
&F_{r_{ij}} = (
\psi _{i}, \boldsymbol{n} \cdot \boldsymbol{\varphi}_{j})_{L_{2}(
\partial \Omega )}, \quad {J}_{r_{ijk}} = (\nabla \psi_i, \nabla \cdot (\bm{\varphi}_j \otimes {\bm{\varphi}}_k))_{L_2(\Omega)}. \label{eq:pressure_matrices2} 
\end{align}

To conclude, the ROM algebraic system that has to be solved at every time step is \eqref{eq:reduced_1}, \eqref{eq:reduced2_1} and \eqref{eq:reduced_q}.

\section{Numerical results}\label{sec:results}

We test our approach on two well-known test cases \cite{John2004,turek1996}: 2D and 3D flow past a cylinder at $0 \leq Re(t) \leq 100$. 
Our goal is a thorough assessment of our ROM model in the reconstruction of the time evolution of the flow field. 

\subsection{2D flow past a cylinder}\label{sec:2D}
\begin{figure}[h]
\centering
\includegraphics[height=0.15\textwidth]{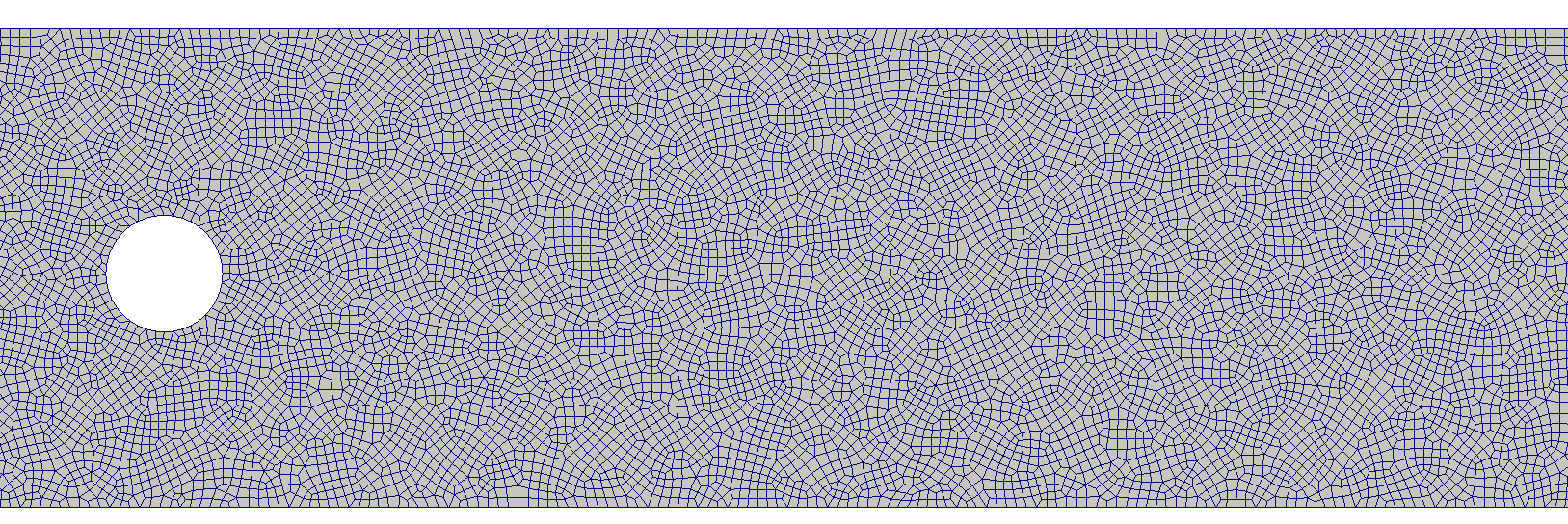}~
\includegraphics[height=0.145\textwidth]{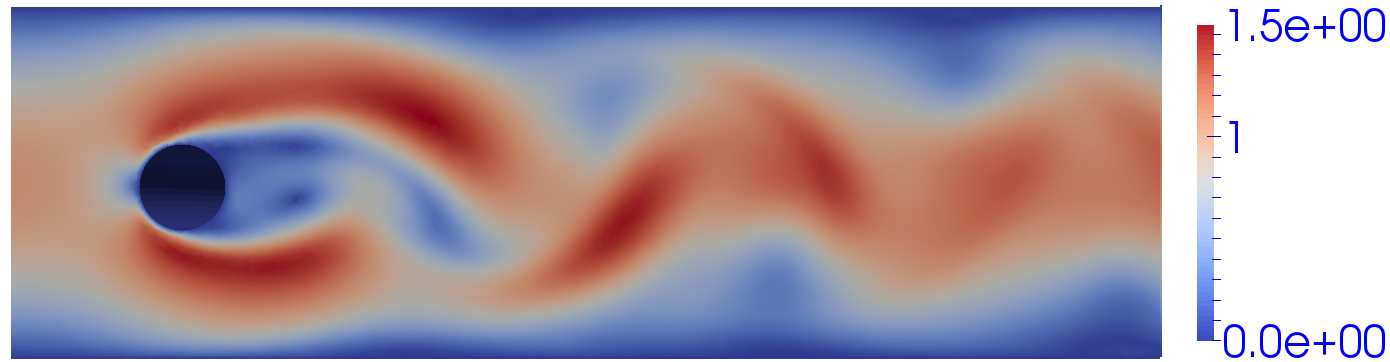}
\caption{2D flow past a cylinder: (left) part of the mesh under consideration and (right) illustrative representation of the velocity field 
in part of the domain at $t = 6$.
}
\label{fig:example_cyl}
\end{figure}

The computational domain is a 2.2 $\times$ 0.41 rectangular channel with a cylinder of radius 0.05 centered at (0.2, 0.2), 
when taking the bottom left corner of the channel as the origin of the axes. 
Fig.~\ref{fig:example_cyl} (left) shows part of the computational domain. 
The channel is filled with fluid with density $\rho = 1$ and viscosity $\mu = 10^{-3}$.
We impose a no slip boundary condition on the upper and lower wall and on the cylinder. At the inflow, we prescribe
the following velocity profile:
\begin{align}\label{eq:cyl_bc}
\v(0,y,t) = \left(\dfrac{6}{0.41^2} \sin\left(\pi t/8 \right) y \left(0.41 - y \right), 0\right), \quad y \in [0, 0.41], \quad t \in (0, 8],
\end{align}
and ${\partial p}/{\partial \n} = 0$. At the outflow we prescribe $\nabla \v \cdot  \n = 0$ and $p = 0$. 
We start the simulations from fluid at rest. 

The quantities of interest for this benchmark are the drag and lift coefficients \cite{John2004,turek1996}:
\begin{align}\label{eq:cd_cl}
C_D(t) = \dfrac{2}{\rho L_{r} {U}^2_{r}} \int_S \left(\left(2 \mu \nabla \u - q\boldsymbol{I}\right)
\cdot \boldsymbol{n}\right) \cdot \boldsymbol{t}~dS, \quad
C_L(t) = \dfrac{2}{\rho L_{r} {U}^2_{r}} \int_S \left(\left(2 \mu \nabla \u - q\boldsymbol{I}\right)
\cdot \boldsymbol{n}\right) \cdot \boldsymbol{n}~dS,
\end{align}
where $U_{r}= 1$ is the maximum velocity at the inlet/outlet, $L_r = 0.1$ is the cylinder diameter, 
$S$ is the cylinder surface, and $\boldsymbol{t}$ and $\boldsymbol{n}$ are the tangential and normal unit vectors
to the cylinder, respectively.

We have previously investigated this benchmark at the FOM level in \cite{Girfoglio2019} and the at ROM level 
in \cite{Girfoglio_JCP, Strazzullo2021}. However, while at the FOM level we considered the EFR algorithm with
a linear indicator function and with indicator function \eqref{eq:a_D0_a_D1} \cite{Girfoglio2019},  
the ROM study is limited to the EFR algorithm with a linear indicator function
\cite{Girfoglio_JCP, Strazzullo2021}. 
In \cite{Girfoglio2019}, we showed that at the FOM level the oscillatory pattern of the flow field (see Figure \ref{fig:example_cyl} (right)) 
can be recovered only with a nonlinear indicator function. 
Thus, 
we expect the results obtained with the ROM presented in this paper to be far more accurate than the results in \cite{Girfoglio_JCP}. 


We consider a hexaedral computational grid with $h_{min} =  4.2e-3$, $h_{avg} = 7.5e-3$ and $h_{max} = 1.1e-2$
for a total of $1.59e4$ cells. The quality of the mesh is high: it features very low values of maximum non-orthogonality (36$^\circ$), average non-orthogonality (4$^\circ$), skewnwss (0.7), and maximum aspect ratio (2). 
Fig.~\ref{fig:example_cyl} (left) shows a part of the mesh. We chose this mesh because it is the coarsest among all the meshes considered in \cite{Girfoglio2019} and thus the most challenging for our filtering approach.

\subsubsection{Validation of the FOM}

Before applying the ROM, we test the EFR algorithm at FOM level and compare its results with the ones produced by a NSE solver
in OpenFOAM. 
For the convective term, we use a second-order accurate central difference scheme that features low dissipation \cite{Lax1960}.
This is a difference with respect to \cite{Girfoglio2019} where we used a second-order accurate upwind scheme.  
We set $\Delta t = 1e-4$ which allows to obtain $CFL_{max} \approx 0.25$ at the time when the velocity reaches its maximum value. 
We set $\chi = \Delta t$ since this is a reasonable choice for academic problems such as the one we are considering \cite{layton_CMAME}. 
More realistic applications require a suitable formula to set $\chi$ \cite{BQV,Girfoglio2019}.
We set $\alpha = h_{avg}$. 

\begin{figure}[h]
\centering
\includegraphics[height=0.33\textwidth]{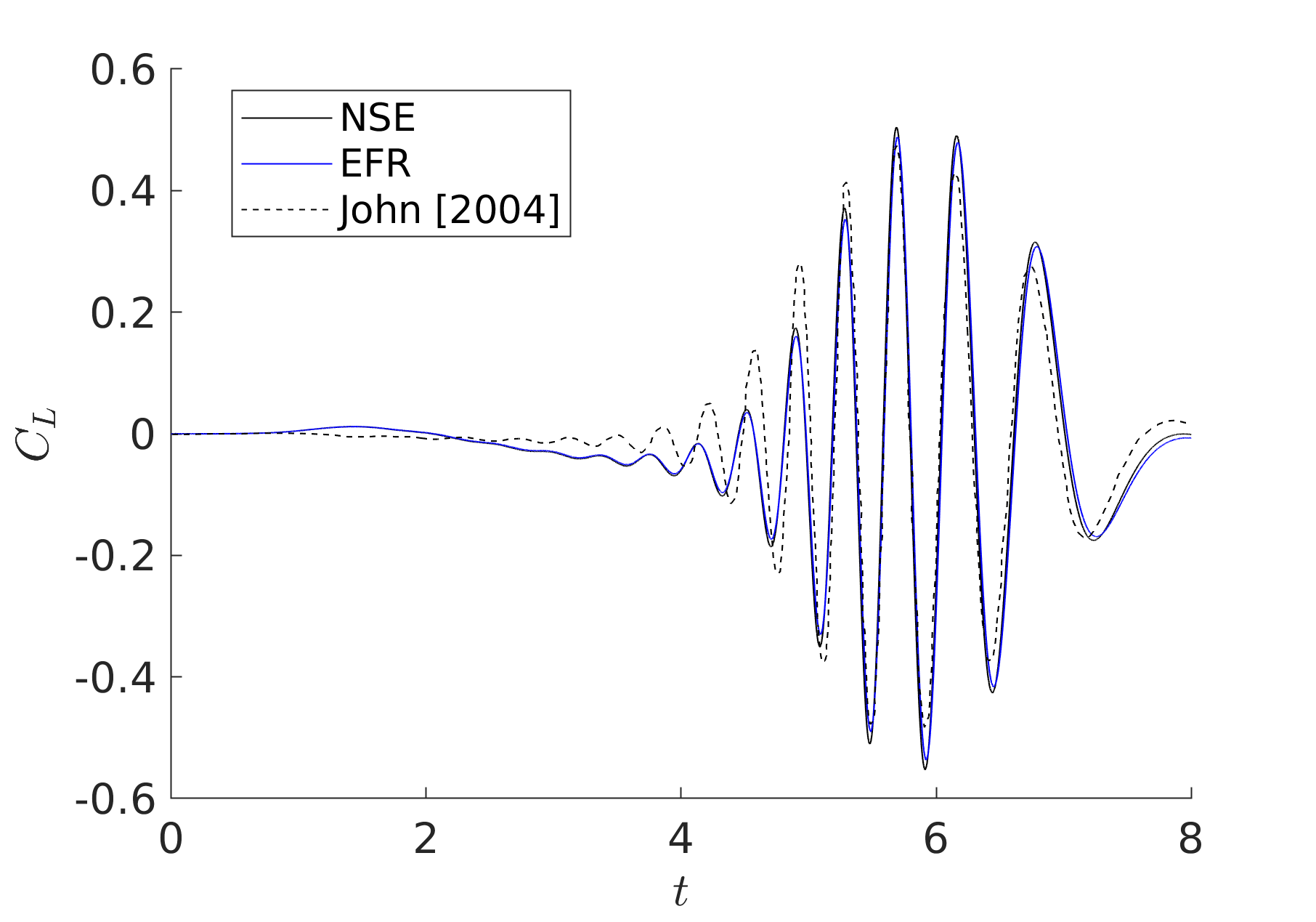}
\includegraphics[height=0.33\textwidth]{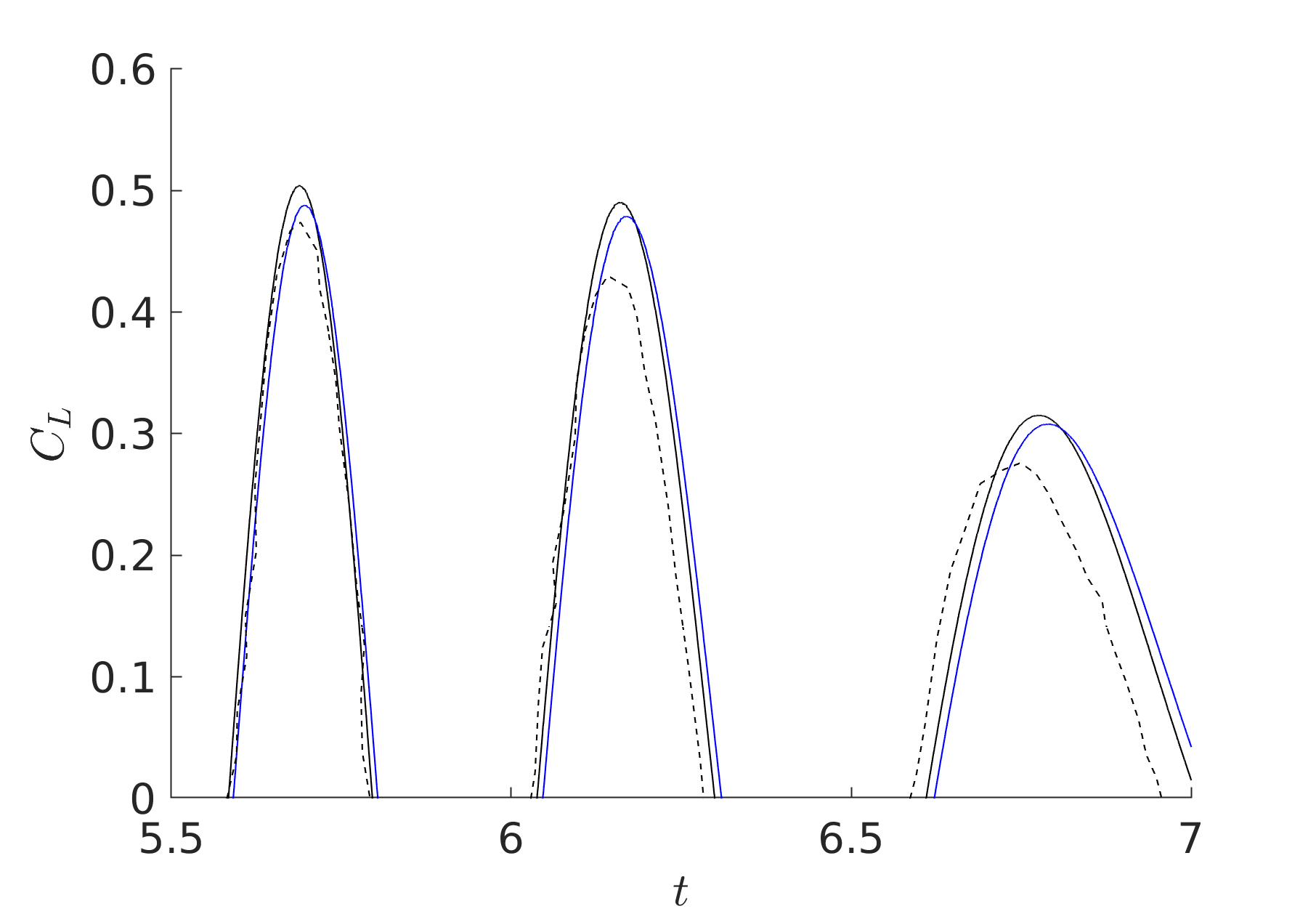}
\caption{2D flow past a cylinder: (left) time evolution of the lift coefficient given by the NSE and EFR for $\alpha = h_{avg}$ and $\chi = \Delta t$ compared against the results in \cite{John2004} and (right) close-up view.}
\label{fig:CL_FOM}
\end{figure}

Fig.~\ref{fig:CL_FOM} (left) shows the evolution of $C_L$ over time computed by EFR and NSE and a comparison
with the results from \cite{John2004}.
Fig.~\ref{fig:CL_FOM} (right) shows a close-up of the time window next to the time of maximum $C_L$. 
We observe that the lift coefficient computed with EFR is slightly closer the reference results from 
\cite{John2004}. To quantify this better agreement, we report in Table \ref{tab:CL_FOM} the computed values of the maximum 
lift coefficient and the corresponding time instant, together with the values from \cite{John2004}. 
Then, we can conclude that EFR is a little more accurate than NSE model when using a coarse mesh
even at the low Reynolds numbers we are considering.

\begin{table}
\centering
\begin{tabular}{lcc}
\multicolumn{2}{c}{} \\
\cline{1-3}
 & $t_{C_{L,max}}$ & $C_{L,max}$ \\
\hline
NSE & 5.688 & 0.504 \\
EFR & 5.699 & 0.488 \\
\cite{John2004} & 5.694 & 0.478 \\
\hline
\end{tabular}
\caption{2D flow past a cylinder: maximum $C_L$ and time instant at which it occurs given by NSE and EFR algorithms. 
Reference values from \cite{John2004} are also reported.}
\label{tab:CL_FOM}
\end{table}



\begin{figure}[h]
\centering
\includegraphics[height=0.35\textwidth]{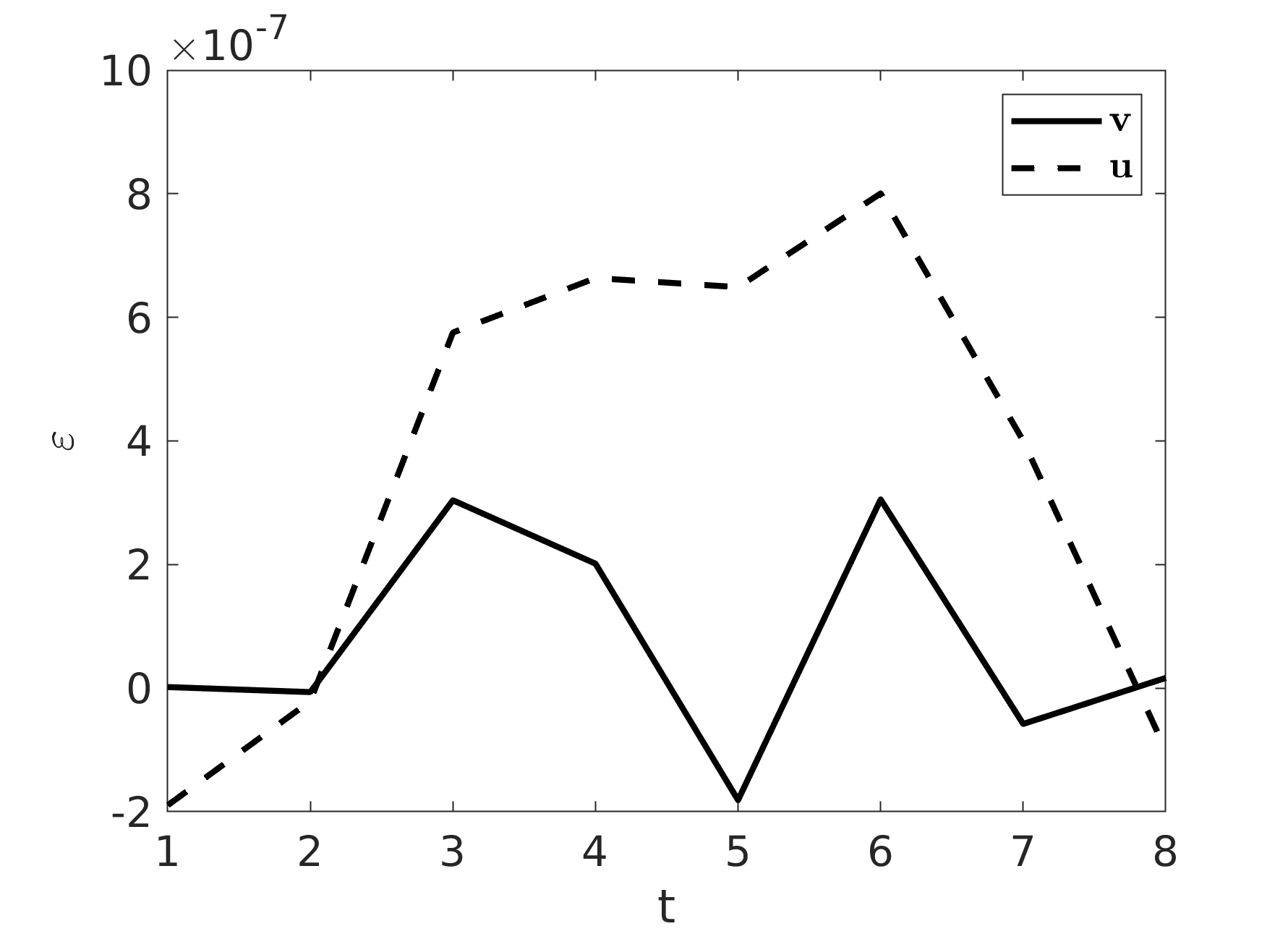}
\includegraphics[height=0.35\textwidth]{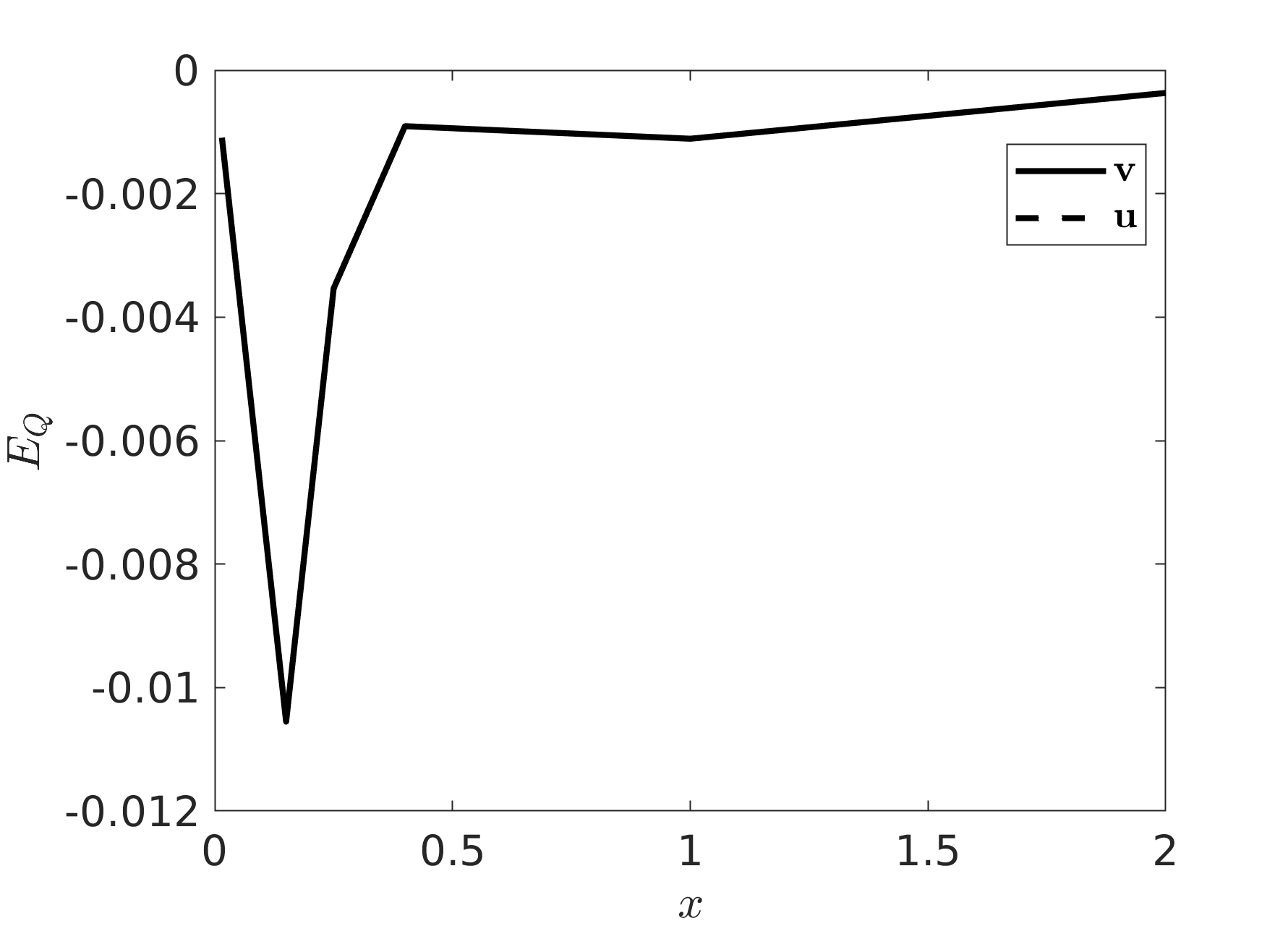}
\caption{2D flow past a cylinder: (left) mass conservation error \eqref{eq:err_mass} for the intermediate velocity $\v$ 
and the end-of-step velocity $\u$ at selected times in the time interval of interest and (right) mass conservation error \eqref{eq:err_mass_2} for the intermediate velocity $\v$ and the end-of-step velocity 
$\u$ for $t = 6$ at selected axial locations.}
\label{fig:err_div}
\end{figure}

Finally, we show in Fig. \ref{fig:err_div} (left) the mass conservation error over time defined as follows:
\begin{align}\label{eq:err_mass}
\varepsilon_{\bm{\Phi}} = \dfrac{1}{\Omega} \int_\Omega \div \bm{\Phi} d\Omega \quad \text{for} \quad \Phi = \{\v, \u\}.
\end{align}
We see that although the incompressibility constraint for the filtered velocity $\vbar$ is not enforced in the model, 
the mass conservation error for the end-of-step velocity $\u$ is comparable to the 
the mass conservation error for the intermediate velocity $\v$ (which, instead, is divergence free). We also evaluated 
the conservation of mass at the selected axial locations $x$ for $t = 6$ (i.e., when $\varepsilon_{\u}$ reaches its maximum value) 
using the error metric proposed in \cite{Passerini2013}:
\begin{align}\label{eq:err_mass_2}
E_{Q_{\bm{\Phi}}} = \dfrac{Q_{\bm{\Phi}} - \widetilde{Q}}{\widetilde{Q}} \quad \text{for} \quad \Phi = \{\v, \u\},
\end{align}
where $Q_{\bm{\Phi}}$ is the volumetric flow rate computed from the numerical axial velocity profiles and $\widetilde{Q}$ 
is the exact volumetric flow rate. 
Fig.~\ref{fig:err_div} (right) shows $E_{Q_{\bm{\v}}}$ and $E_{Q_{\bm{\u}}}$, which are 
overlapped over the whole axis. 

Next, we are going to validate our ROM approach. 
Since this benchmark is characterized by a flow field exhibiting a wide spectrum of frequencies, we adopt an idea 
proposed in \cite{Strazzullo2021}: we test the performances of our ROM approach 
(i) over the entire time window of interest $[0, 8]$ and (ii) over the second half of the time interval $[4, 8]$ where 
the high frequency modes are dominant.

\subsubsection{Validation of the ROM (i)}

\begin{figure}[h]
\centering
\includegraphics[height=0.35\textwidth]{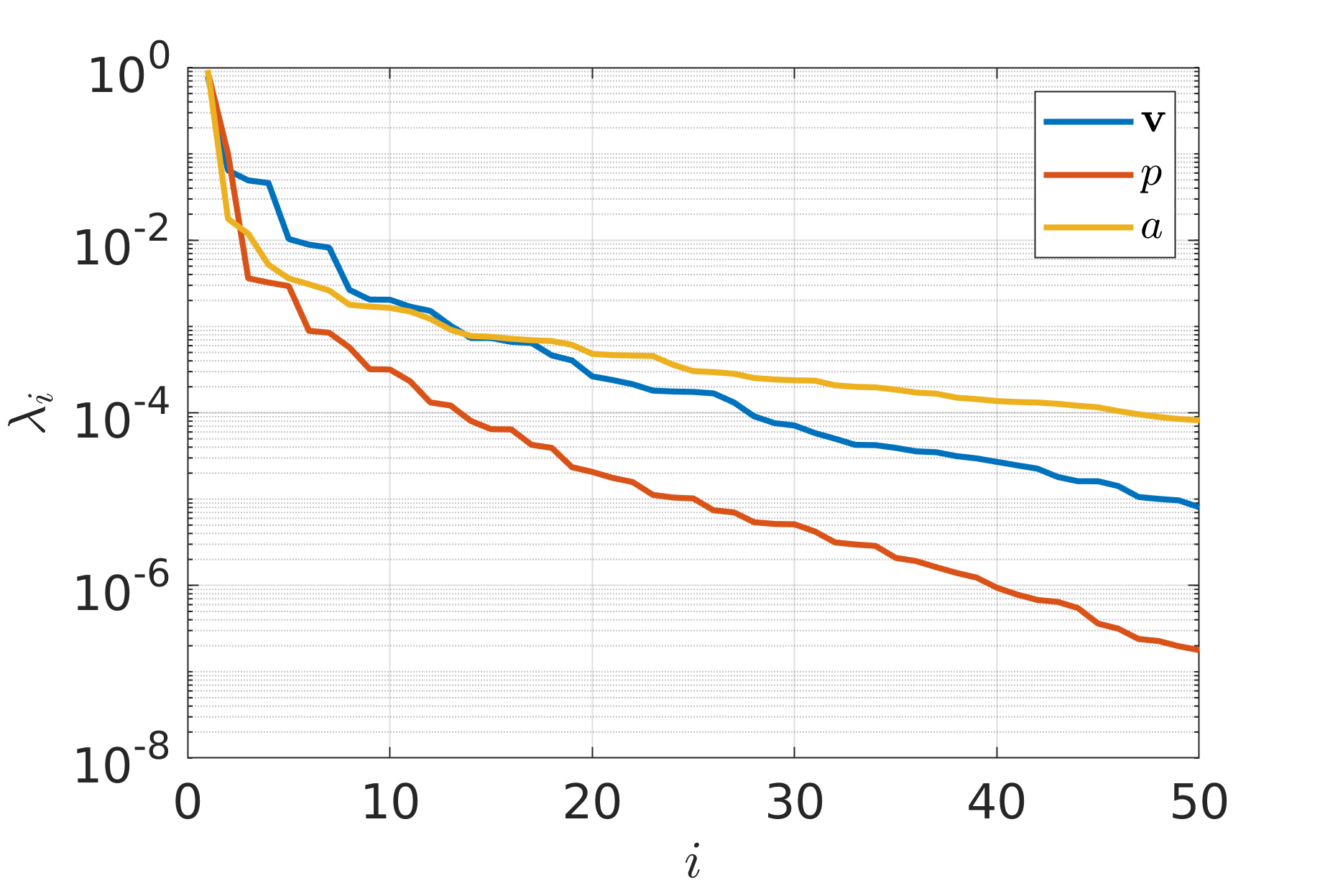}
\caption{2D flow past a cylinder - experiment (i): eigenvalue decay for the intermediate
velocity, pressure, and indicator function. 
}
\label{fig:eig_lift}
\end{figure}

We collect 400 FOM snapshots, one every 0.02 s, i.e.~we use an equispaced grid in time. 
Fig. \ref{fig:eig_lift} shows the eigenvalue decay for the intermediate velocity, pressure, and indicator function. 
For the ROM simulations, we collected the solutions every 0.01 s.
This means that the set of time samples includes the samples used in the offline stage and samples in between 
two consecutive offline samples. The reason for this choice is that we want to assess how accurate the reduced order
approximation is for time instants that were not in the training set. 

\begin{figure}
\centering
 \begin{overpic}[width=0.42\textwidth]{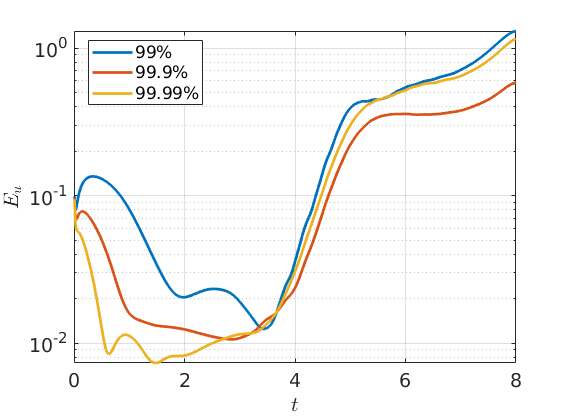}
      \end{overpic}
\begin{overpic}[width=0.42\textwidth]{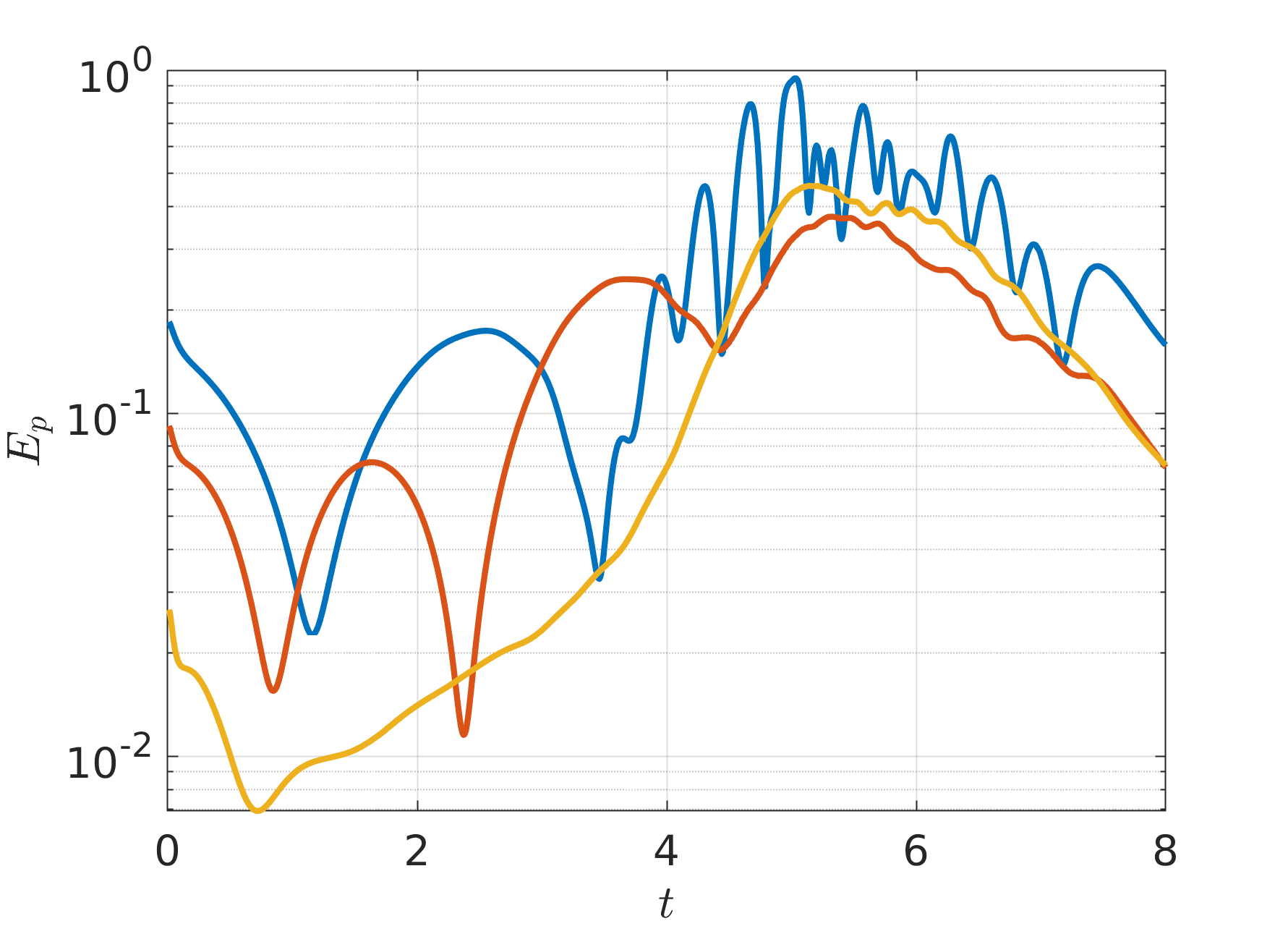}
      \end{overpic}\\ 
\begin{overpic}[width=0.42\textwidth]{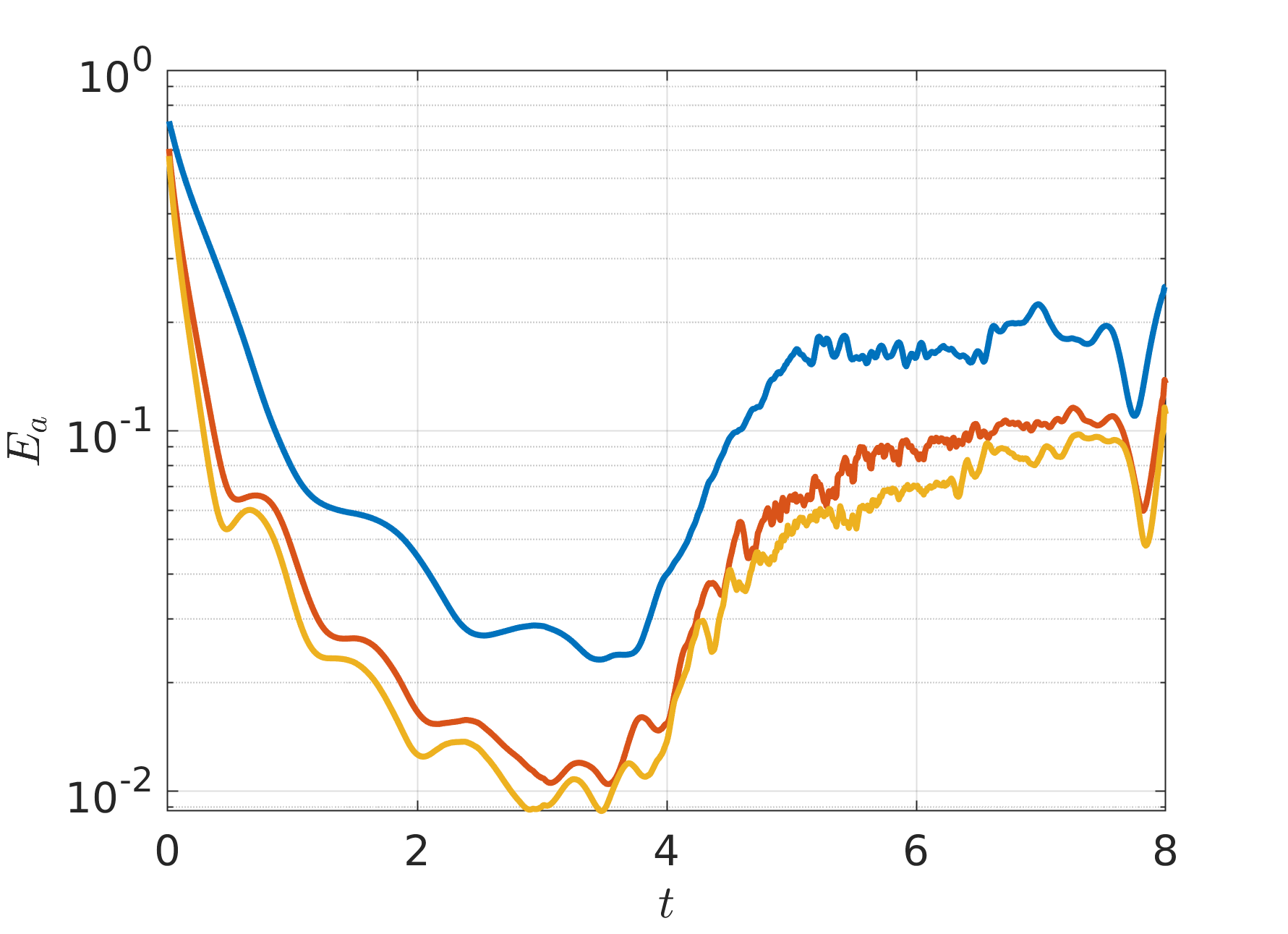}
      \end{overpic}\\ 
\caption{2D flow past a cylinder - experiment (i): time history of the relative $L^2$ error \eqref{eq:error1} 
for velocity $\u$ (top left), pressure $p$ (top right), and for the indicator function $a$ (bottom) for different 
thresholds of cumulative energy.}
\label{fig:err_t_global}
\end{figure}

We performed a convergence test as the number of modes increases. We considered three different energy 
thresholds using the first 50 most energetic POD modes: 99\% (11 modes for $\v$, 4 modes for $p$ and 15 modes for $a$), 
99.9\% (26 modes for $\v$, 10 modes for $p$ and 41 modes for $a$) and 
99.99\% (43 modes for $\v$, 22 modes for $p$ and 50 modes for $a$). We calculated the $L^2$ relative error: 
\begin{equation}\label{eq:error1}
E_{\Phi}(t) = \dfrac{||\Phi_h(t) - \Phi_r(t)||_{L^2(\Omega)}}{||{\Phi_h}(t)||_{L^2(\Omega)}},
\end{equation}
where $\Phi_h$ is a field computed with the FOM ($\u_h$, $p_h$ or $a_h$) 
and $\Phi_r$ is the corresponding field computed with the  ROM ($\u_r$, $p_r$ or $a_r$).
Figure ~\ref{fig:err_t_global} shows errors \eqref{eq:error1} and Table \ref{tab:errors1} reports
minimum, average, and maximum relative errors for 99.99\% of the cumulative energy. 
From Fig.~\ref{fig:err_t_global}, we observe that the majority of the relative errors is less than 1 
at all the time instants. The only exception is the velocity: 
we see that the error increases steeply past $t \approx 4$ and its value reaches 1 towards the end of the time interval.
The relative error for the velocity and pressure is significantly lower for $t<4$.
Both errors increase when the vortex shedding starts at around $t  = 4$.
The relative error associated to the indicator function seems to be less critical since its value remains 
below $10^{-1}$ for most of the time interval of interest for 99.9\% and 99.99\% of the cumulative energy.
Larger errors for the indicator function at the beginning of simulation might be due to the transient nature of the flow. 
This different behavior for the velocity and pressure errors on one side and the indicator
function on the other side could be explained by the fact that differente strategies
are used for the ROM reconstruction (a projection method for $\u$ and $p$ 
and an interpolation procedure for $a$).
Indeed, for what concerning the velocity, by moving from 99\%, to 99.9\% and 99.99\%, it becomes lower. 
We observe that going from 99\% to 99.9\% of the cumulative energy there is a general improvement of the errors, 
while such improvement lessens when  going from 99.9\% to 99.99\%.
From Fig.~\ref{fig:err_t_global} (top right), we see that the oscillations in the pressure error for $t>4$ 
are damped when a larger amount of energy snapshots is retained. 

\begin{table}
\centering
\begin{tabular}{lccc|ccc}
\multicolumn{2}{c}{} \\
\cline{1-7}
 & $\u$ & $p$ & $a$ & $\u$ & $p$ & $a$ \\
\hline
Maximum $E_\Phi$ & 0.09  & 0.07  & 0.58  & 1.15 & 0.46  & 0.12  \\
Average $E_\Phi$ & 0.01  & 0.02  & 0.04  & 0.5  & 0.26  & 0.06\\
Miminum $E_\Phi$ & 0.007 & 0.007  & 0.01  & 0.03  & 0.07  & 0.01\\
\hline
\end{tabular}
\caption{2D flow past a cylinder - experiment (i): maximum, average, and
minimum relative $L^2$ errors for the end-of-step velocity, pressure, and indicator function 
for 99.99\% of the cumulative energy. The first three columns refer to the first half of the time interval (i.e., [0, 4]), 
while the last three refer to the second half of the time interval (i.e., [4, 8]).}
\label{tab:errors1}
\end{table}

\begin{figure}
\centering
 \begin{overpic}[width=0.44\textwidth]{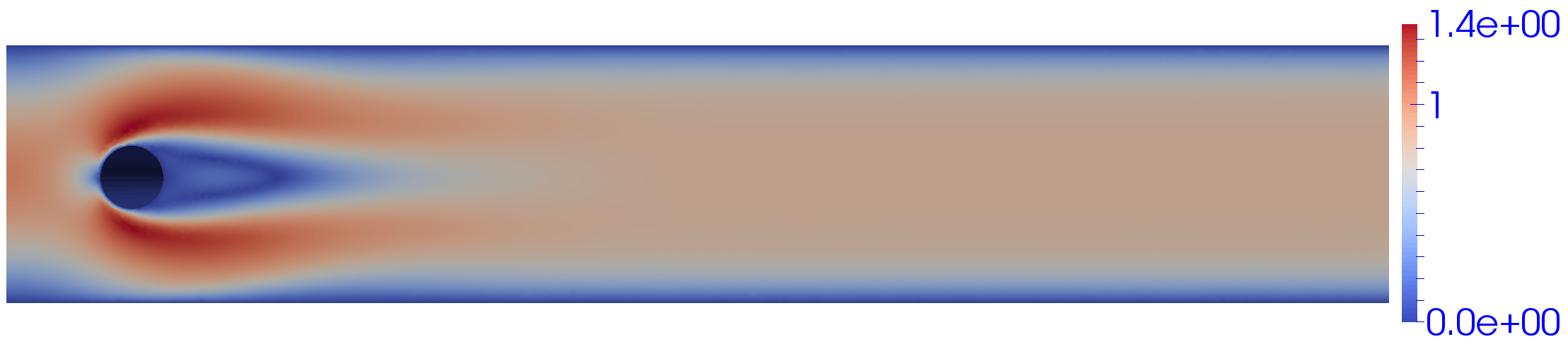}
        \put(35,20){FOM}
        \put(-10,10){$\u$}
      \end{overpic}
 \begin{overpic}[width=0.44\textwidth]{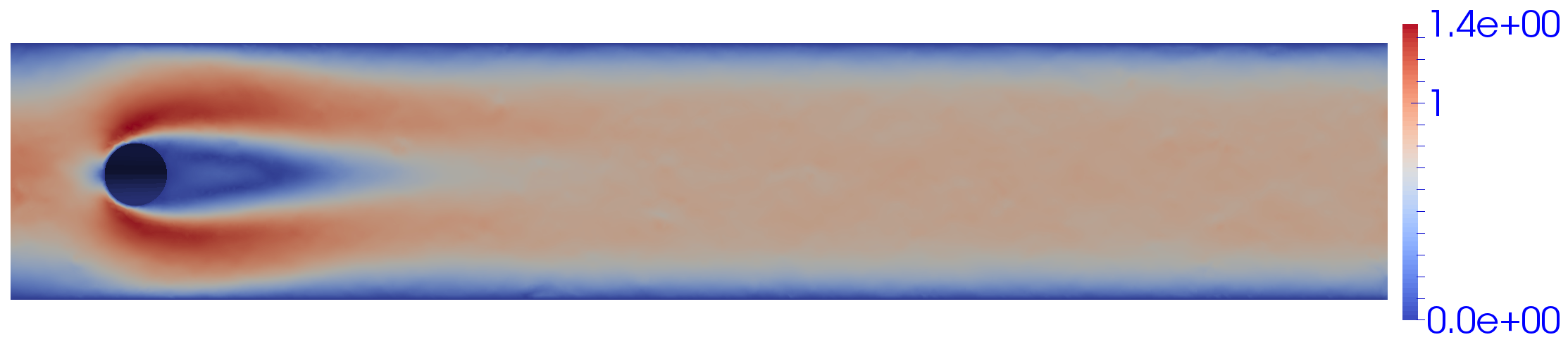}
         \put(35,20){ROM}
      \end{overpic}\\
 \begin{overpic}[width=0.44\textwidth]{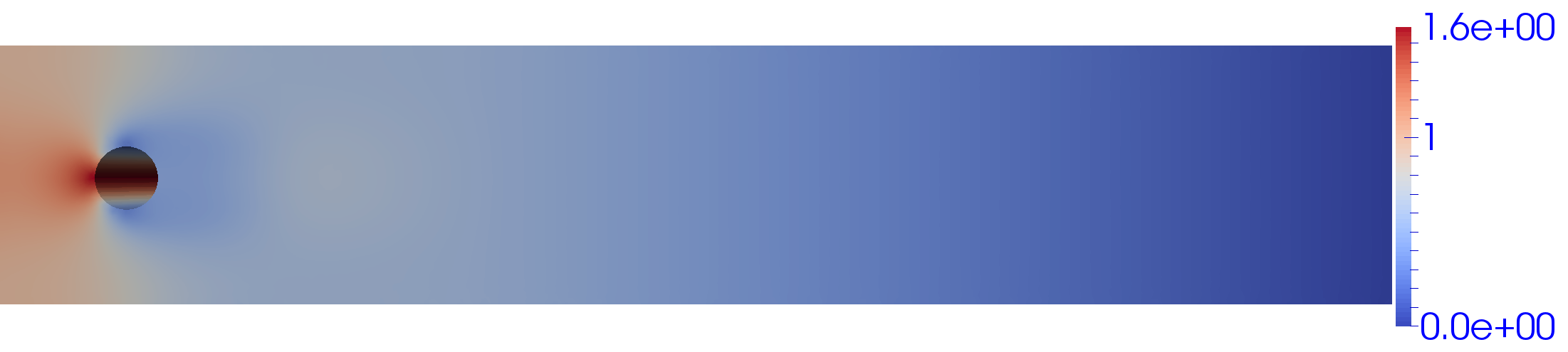}
        \put(-10,10){$p$}
      \end{overpic}
 \begin{overpic}[width=0.44\textwidth]{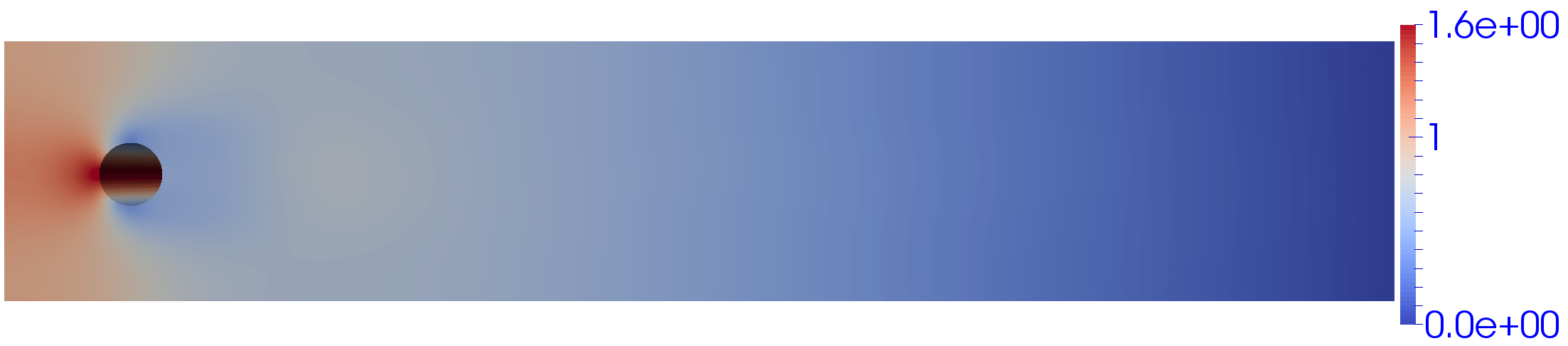}
      \end{overpic}\\
 \begin{overpic}[width=0.44\textwidth]{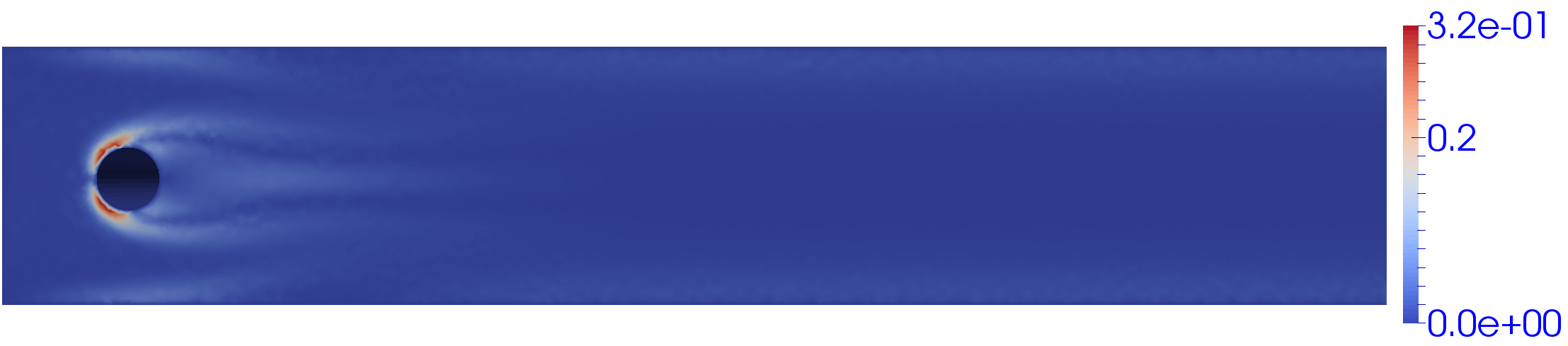}
        \put(-10,10){$a$}
      \end{overpic}
 \begin{overpic}[width=0.44\textwidth]{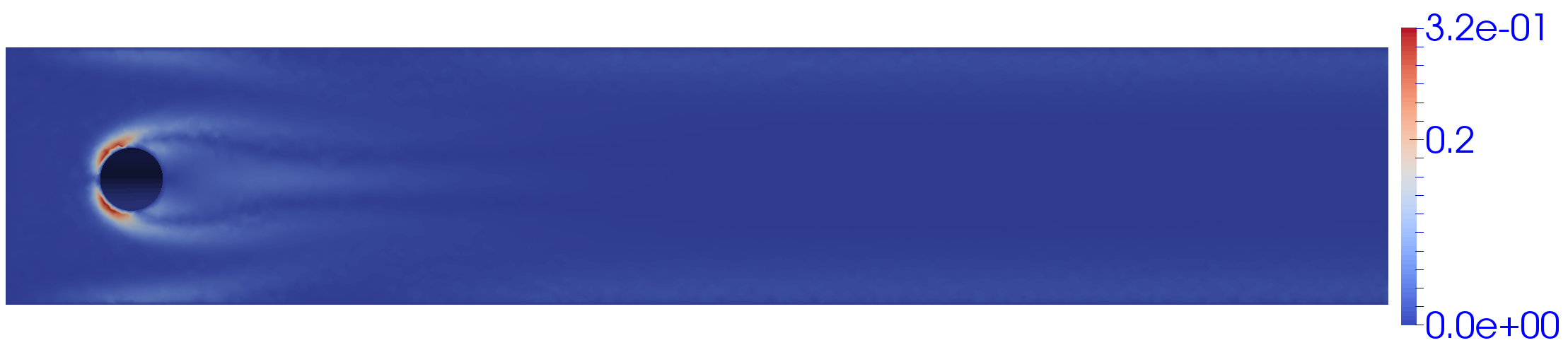}
      \end{overpic}\\
\caption{2D cylinder problem - experiment (i): comparison between 
velocity $\u$ (first raw), pressure (second raw), and indicator function (3nd raw)
computed by the FOM (left) and the ROM (right)  at $t = 1.9$ for 99.99\% of the cumulative energy.
}
\label{fig:comp_t_1_9}
\end{figure}

\begin{figure}
\centering
 \begin{overpic}[width=0.44\textwidth]{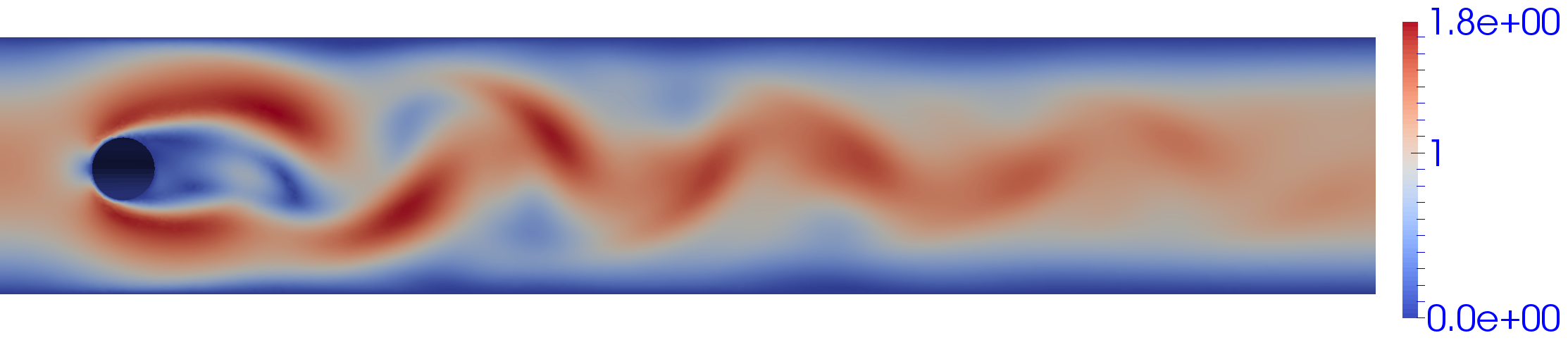}
        \put(35,20){FOM}
        \put(-10,10){$\u$}
      \end{overpic}
 \begin{overpic}[width=0.44\textwidth]{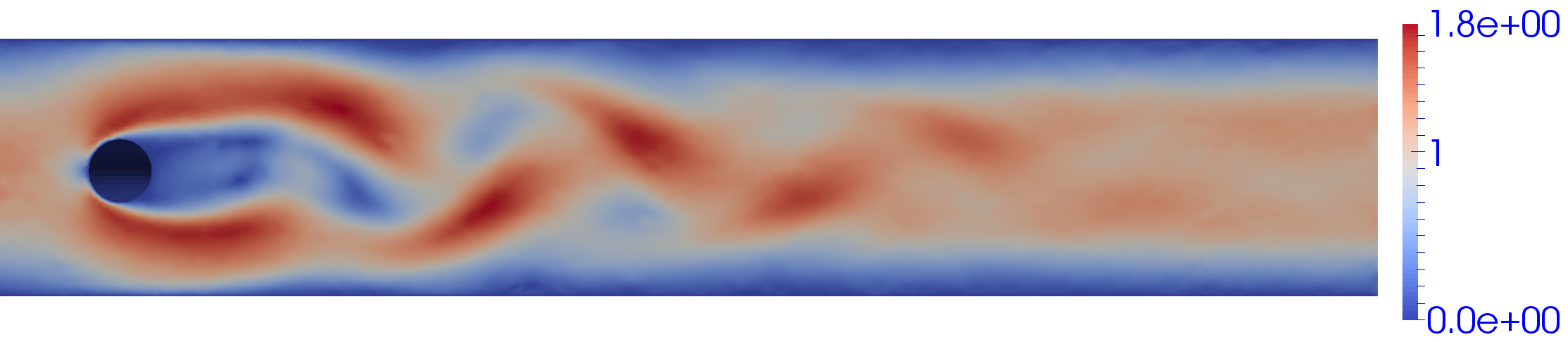}
         \put(35,20){ROM}
      \end{overpic}\\
 \begin{overpic}[width=0.44\textwidth]{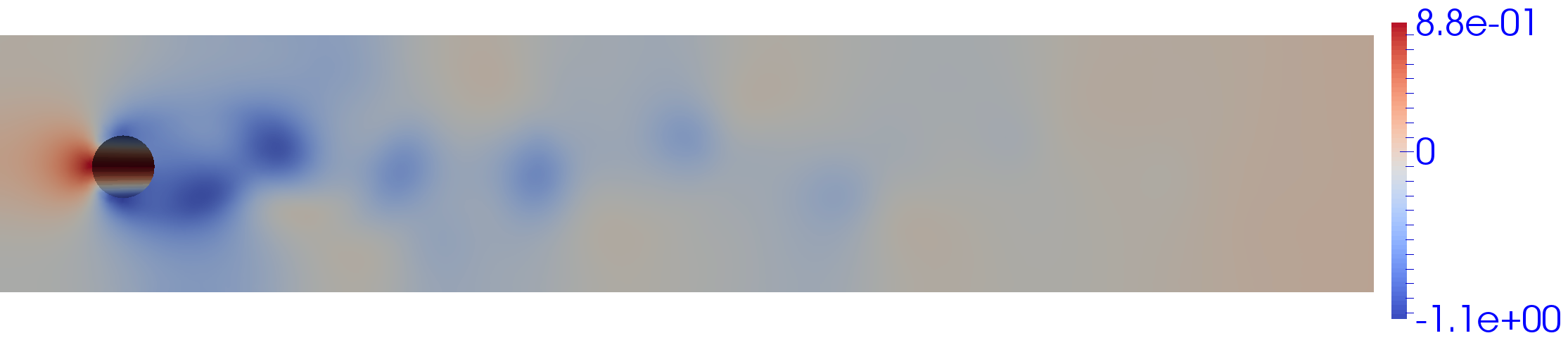}
        \put(-10,10){$p$}
      \end{overpic}
 \begin{overpic}[width=0.44\textwidth]{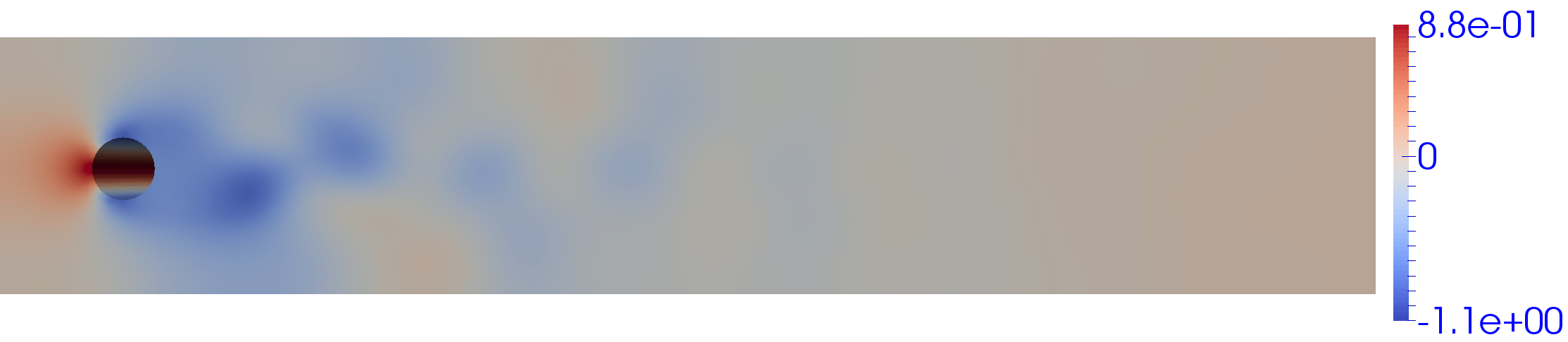}
      \end{overpic}\\
 \begin{overpic}[width=0.44\textwidth]{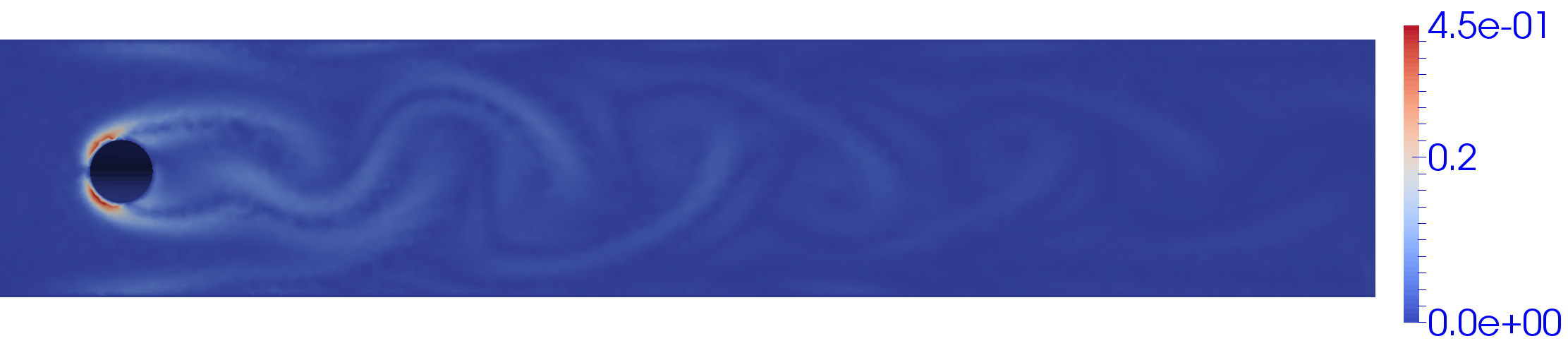}
        \put(-10,10){$a$}
      \end{overpic}
 \begin{overpic}[width=0.44\textwidth]{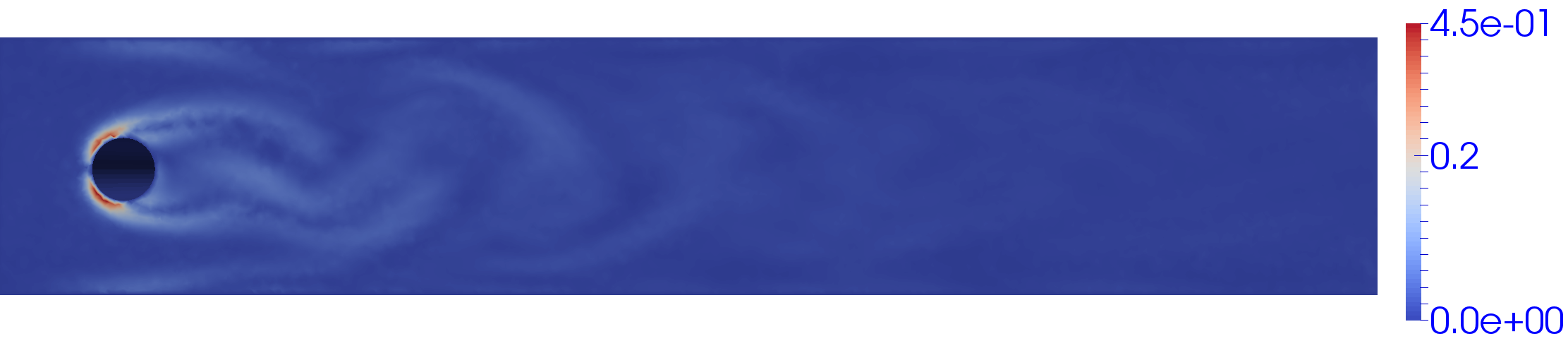}
      \end{overpic}\\
\caption{2D cylinder problem - experiment (i): comparison between 
velocity $\u$ (first raw), pressure (second raw), and indicator function (3nd raw)
computed by the FOM (left) and the ROM (right) at $t = 5.5$ for 99.99\% of the cumulative energy. 
}
\label{fig:comp_t_5_5}
\end{figure}

Figures~\ref{fig:comp_t_1_9} and~\ref{fig:comp_t_5_5} display a qualitative comparison between the computed FOM and ROM fields at two different times: $t = 1.9$ (first half of the time interval) and $t = 5.5$ (second half). As we can see from 
Fig.~\ref{fig:comp_t_1_9}, our ROM provides a good reconstruction of all the variables at $t = 1.9$. 
On the other hand, the ROM fails for provide an accurate approximation of 
velocity and pressure at $t = 5.5$, 
as shown in Fig.~\ref{fig:comp_t_5_5}.

\begin{figure}
\centering
 \begin{overpic}[width=0.45\textwidth]{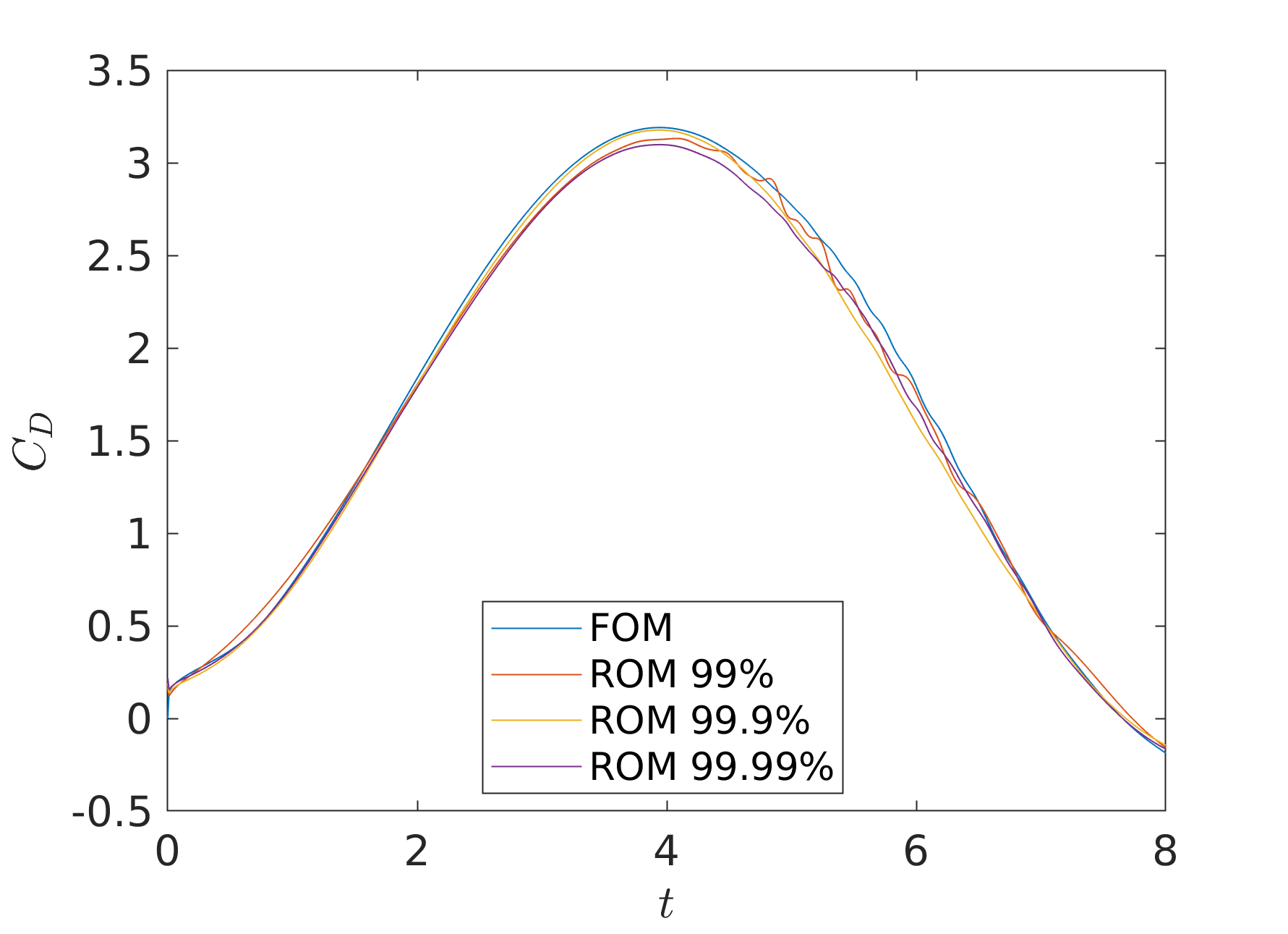}
      \end{overpic}
       \begin{overpic}[width=0.45\textwidth]{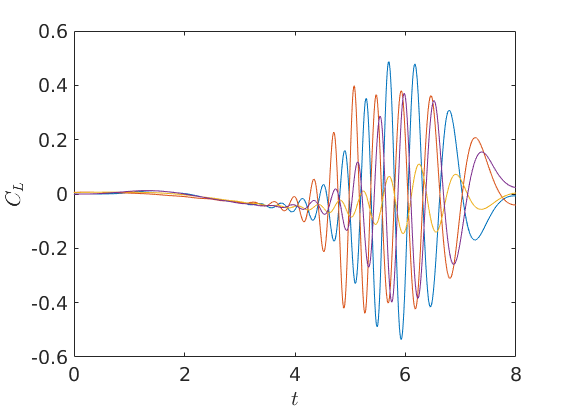}
      \end{overpic}\\
\caption{2D flow past a cylinder - experiment (i): aerodynamic coefficients $C_D$ (left) and $C_L$ (right) computed by FOM and ROM for different thresholds of cumulative energy.}
\label{fig:coeff_t_global}
\end{figure}

Figure~\ref{fig:coeff_t_global} reports a more quantitative comparison: the  aerodynamics coefficients \eqref{eq:cd_cl}
computed by FOM and ROM with different thresholds of cumulative energy.
We observe that the time evolution of the drag coefficient is correctly reconstructed by ROM, 
while the ROM reconstruction of the lift coefficient is accurate till about $t =4$. 
For $t>4$, the lift coefficient computed by the ROM is off in terms of both phase and amplitude
regardless of the percentage of retained energy.


\subsubsection{Validation of the ROM (ii)}

In this second experiment, we focus on the second half of the time interval of interest, $[4, 8]$ with the goal of
improving the ROM reconstruction of velocity and pressure in that time window.
We use the same database of FOM snapshots used in experiment (i), but we 
apply the POD only to the last 200 high-fidelity snapshots, i.e.~those related to $[4, 8]$.
Just like in experiment (i), for the ROM simulations 
we included both samples used in the training phase (i.e., the offline sampes) and samples in between the offline samples. 

\begin{figure}
\centering
 \begin{overpic}[width=0.45\textwidth]{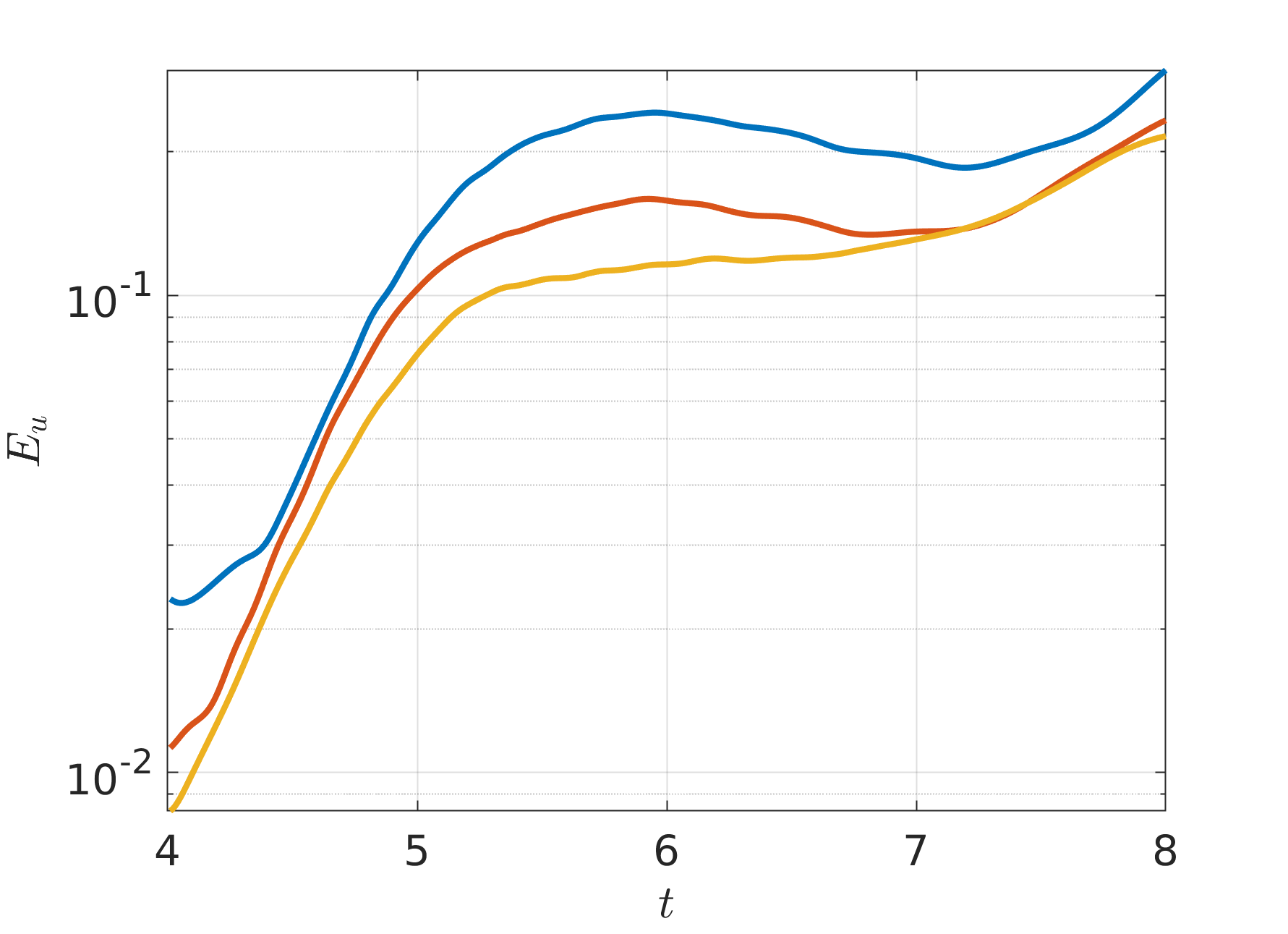}
      \end{overpic}
\begin{overpic}[width=0.45\textwidth]{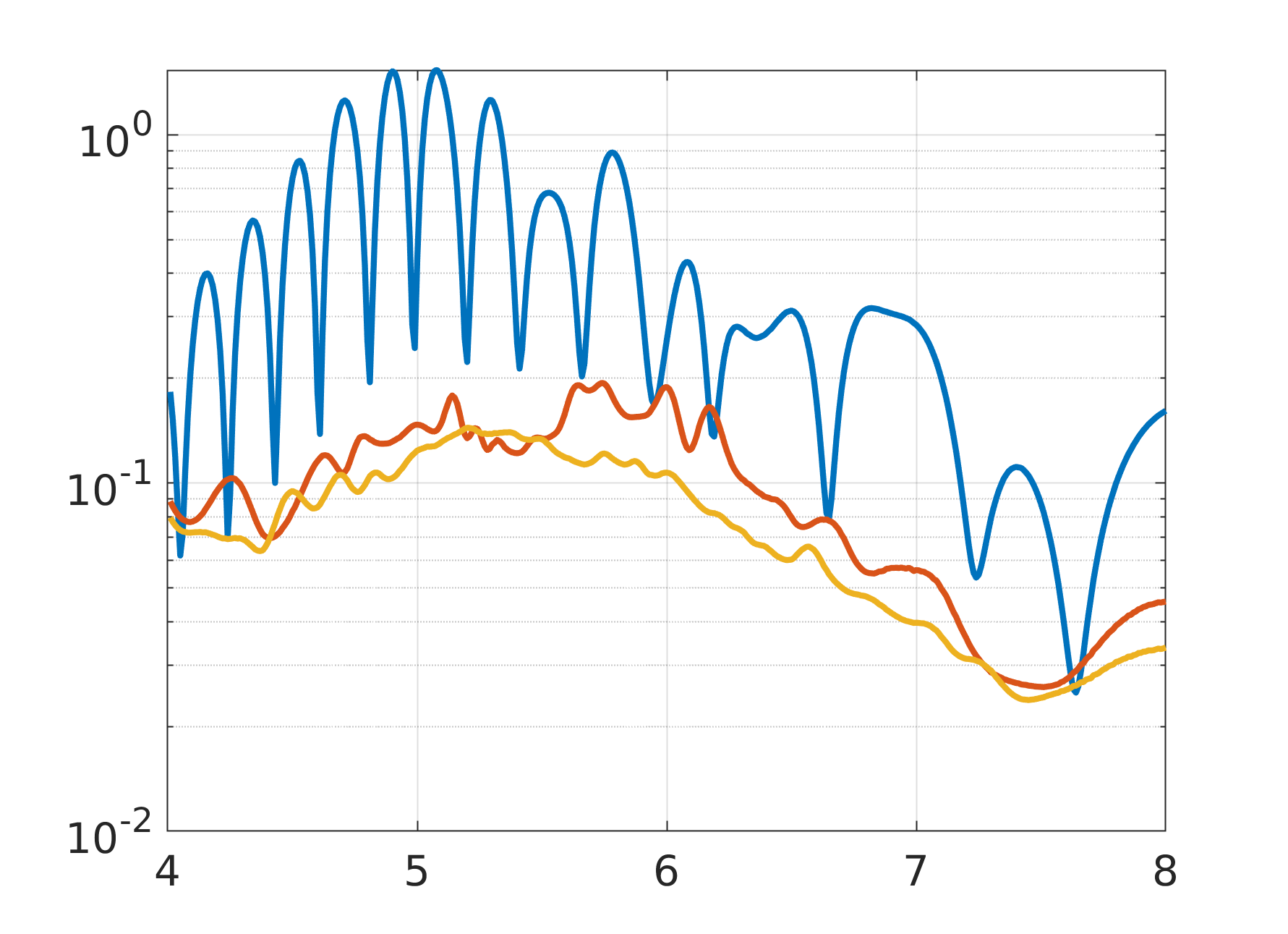}
      \end{overpic}\\ 
\caption{2D flow past a cylinder - experiment (ii): time history of the relative $L^2$ error \eqref{eq:error1} 
for velocity $\u$ (left) and pressure $p$ (right) for different thresholds of cumulative energy.}
\label{fig:err_t_nested}
\end{figure}

Also for this experiment, we considered the first 50 most energetic POD modes and 
performed a convergence test based on three different energy thresholds: 
99\% (11 modes for $\v$, 5 modes for $p$ and 20 modes for $a$), 
99.9\% (26 modes for $\v$, 12 modes for $p$ and 44 modes for $a$), 
and 99.99\% (42 modes for $\v$, 24 modes for $p$ and 50 modes for $a$). 
Fig.~\ref{fig:err_t_nested} shows errors \eqref{eq:error1} and Table \ref{tab:errors2} reports
minimum, average, and maximum relative errors for 99.99\% of the cumulative energy. 
We observe that the relative error for the velocity  reaches much lower values than in experiment (i): 
it drops below 0.22 during the entire time interval when the 99.99\% of the snapshots energy is retained. 
Moreover, there is a monotonic convergence as the number of the modes is increased. 
The improvement of the pressure reconstruction is also rather significative: 
compare Fig.~\ref{fig:err_t_global} (top right) with Fig.~\ref{fig:err_t_nested} (left). 
Once again, we note that there is not much difference in the relative errors
for velocity and pressure when going from  99.9\% to 99.99\% of the cumulative energy. 

\begin{table}
\centering
\begin{tabular}{lcc}
\multicolumn{2}{c}{} \\
\cline{1-3}
 & $\u$ & $p$ \\ 
\hline
Maximum $E_\Phi$ & 0.22  & 0.14  \\
Average $E_\Phi$ & 0.11  & 0.08  \\
Miminum $E_\Phi$ & 0.008  & 0.02  \\
\hline
\end{tabular}
\caption{2D flow past a cylinder - experiment (ii): 
maximum, average, and minimum relative $L^2$ errors for the end-of-step velocity and pressure
for 99.99\% of the cumulative energy.}
\label{tab:errors2}
\end{table}

\begin{figure}
\centering
 \begin{overpic}[width=0.47\textwidth]{img/uFOM_global_5_5s_cut.png}
        \put(35,20){FOM}
        \put(-10,10){$\u$}
      \end{overpic}
 \begin{overpic}[width=0.47\textwidth]{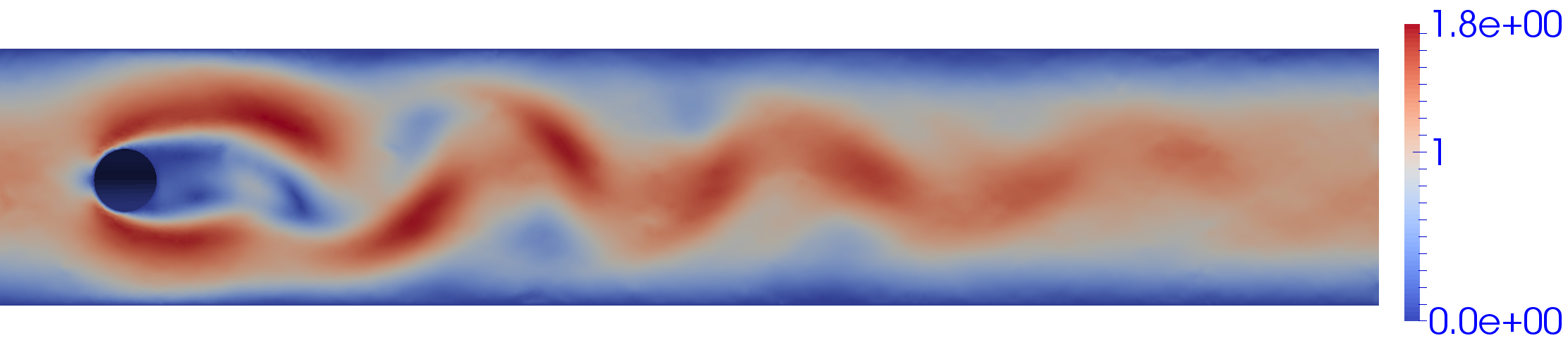}
         \put(35,20){ROM}
      \end{overpic}\\
 \begin{overpic}[width=0.47\textwidth]{img/pFOM_global_5_5s_cut.png}
        \put(-10,10){$p$}
      \end{overpic}
 \begin{overpic}[width=0.47\textwidth]{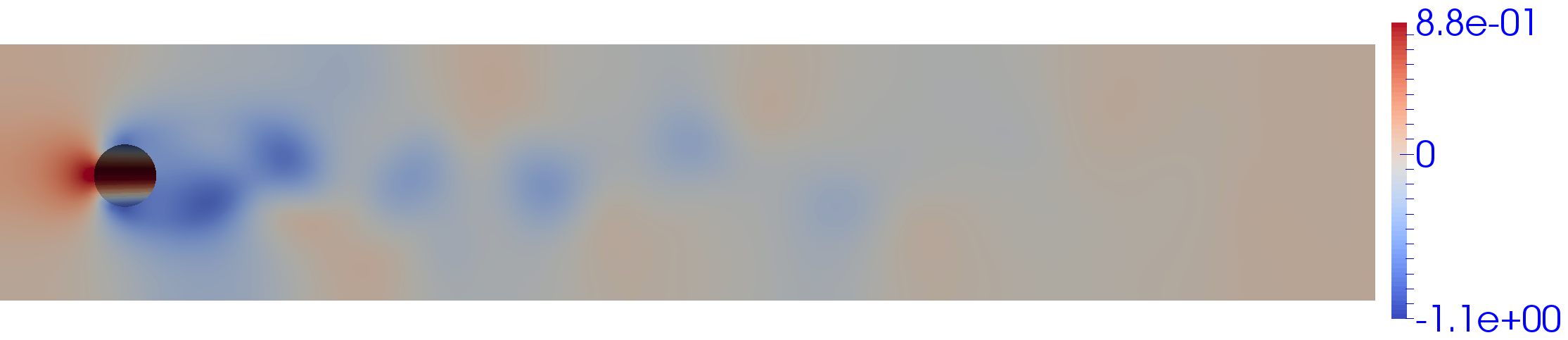}
      \end{overpic}\\
\caption{2D cylinder problem - experiment (ii): comparison between 
velocity $\u$ (first raw) and pressure (second raw)
computed by the FOM (left) and the ROM (right) at $t = 5.5$ for 99.99\% of the cumulative energy. }
\label{fig:comp_u_t_5_5_nested}
\end{figure}

\begin{figure}
\centering
       \begin{overpic}[width=0.6\textwidth]{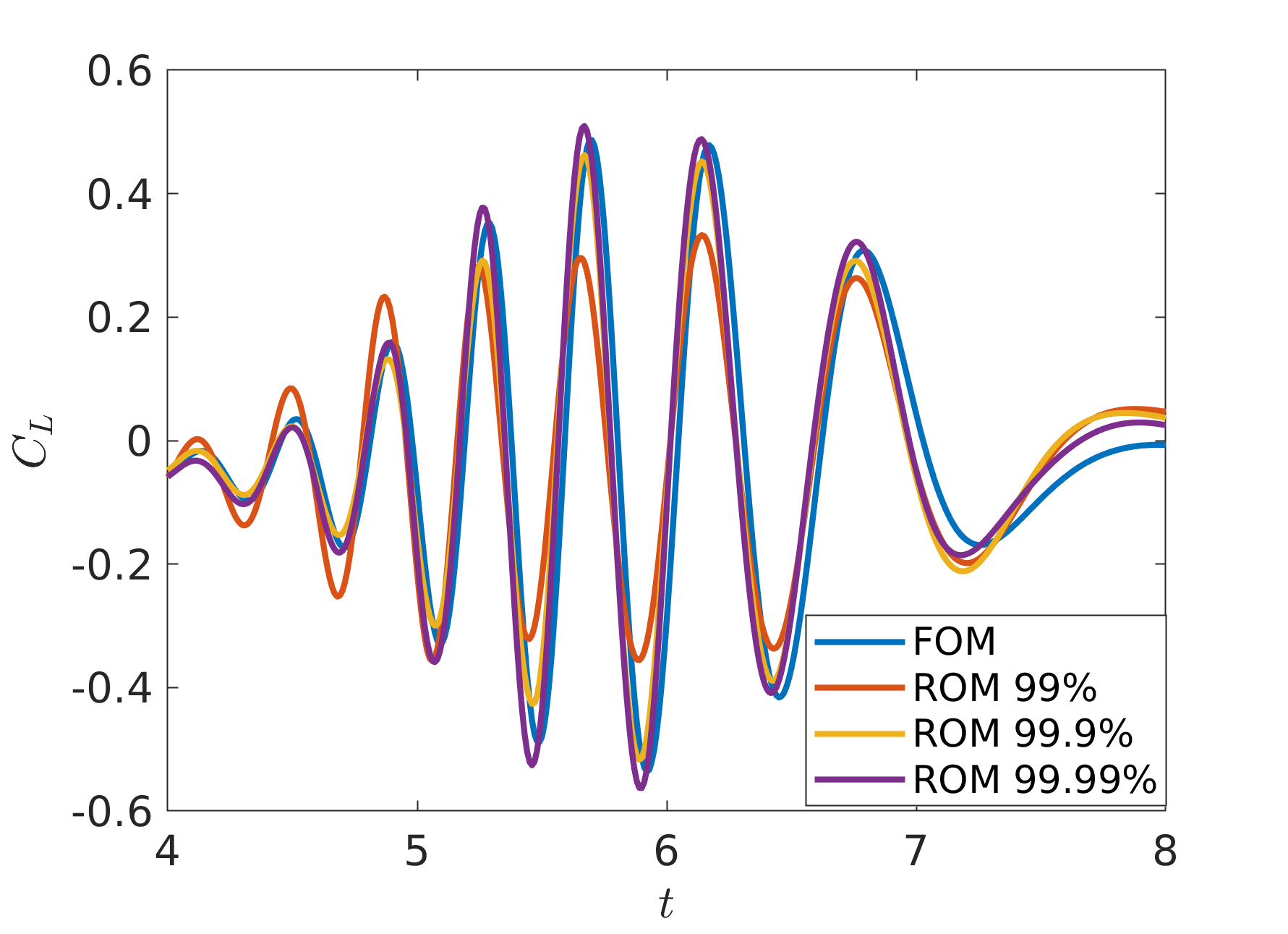}
      \end{overpic}\\
\caption{2D flow past a cylinder - experiment (ii): 
lift coefficient $C_L$ computed by FOM and ROM for different thresholds of cumulative energy. 
}
\label{fig:coeff_t_nested}
\end{figure}

Figure \ref{fig:comp_u_t_5_5_nested} shows the comparison between the computed FOM and ROM fields at $t = 5.5$. 
The ROM reconstructions of velocity and pressure are much more accurate
than in experiment (i). 
Figure~\ref{fig:coeff_t_nested} reports the quantitative FOM/ROM comparison for the lift coefficient $C_L$. 
The improvement with respect to experiment (i) is evident: the time evolution of $C_L$ computed
by the ROM is very accurate when the 99.9\% or 99.99\% of the energy is retained. 
For a further quantitative assessment, we computed the following error
\begin{align}\label{eq:error_coeff}
E_{C_L} = \dfrac{||C_L(t)^{FOM}- C_L(t)^{ROM}||_{L^2(4,8)}}{||C_L(t)^{FOM}||_{L^2(4,8)}}.
\end{align}
We obtain $E_{C_L} = 0.5, 0.39$, and 0.38, for 99\%, 99.9\% and 99.99\% of the cumulative energy, respectively.

We conclude by proving some information about the efficiency of our ROM approach. The total CPU time required by a FOM simulation is about $2300$ s, while the solution of the reduced algebraic systems for experiment (ii) with 99.99\% of the cumulative energy  
takes $570$ s. The resulting speed-up is about $4$, which is significantly lower than the speed-up observed in a 
ROM study for the EF algorithm with a linear indicator function \cite{Girfoglio_JCP, Girfoglio_ROM_Fluids}. 
The reason for this difference does not lie in the nonlinearity of the indicator function. Instead, such difference is due to the 
very large number of modes retained, which increases the size of the reduced dynamical system. 
The larger number of modes is related to the fact that the nonlinear filter combined with the EFR algorithm
preserves the vortex shedding observed with a DNS, while the linear filter combined with the EF algorithm introduces
too much artificial dissipation and the oscillatory flow is lost. Thus, the increased computational cost is the price to pay
for an accurate reconstruction of the flow. 

\subsection{3D flow past a cylinder}
\begin{figure}
\centering
 \begin{overpic}[width=0.3\textwidth]{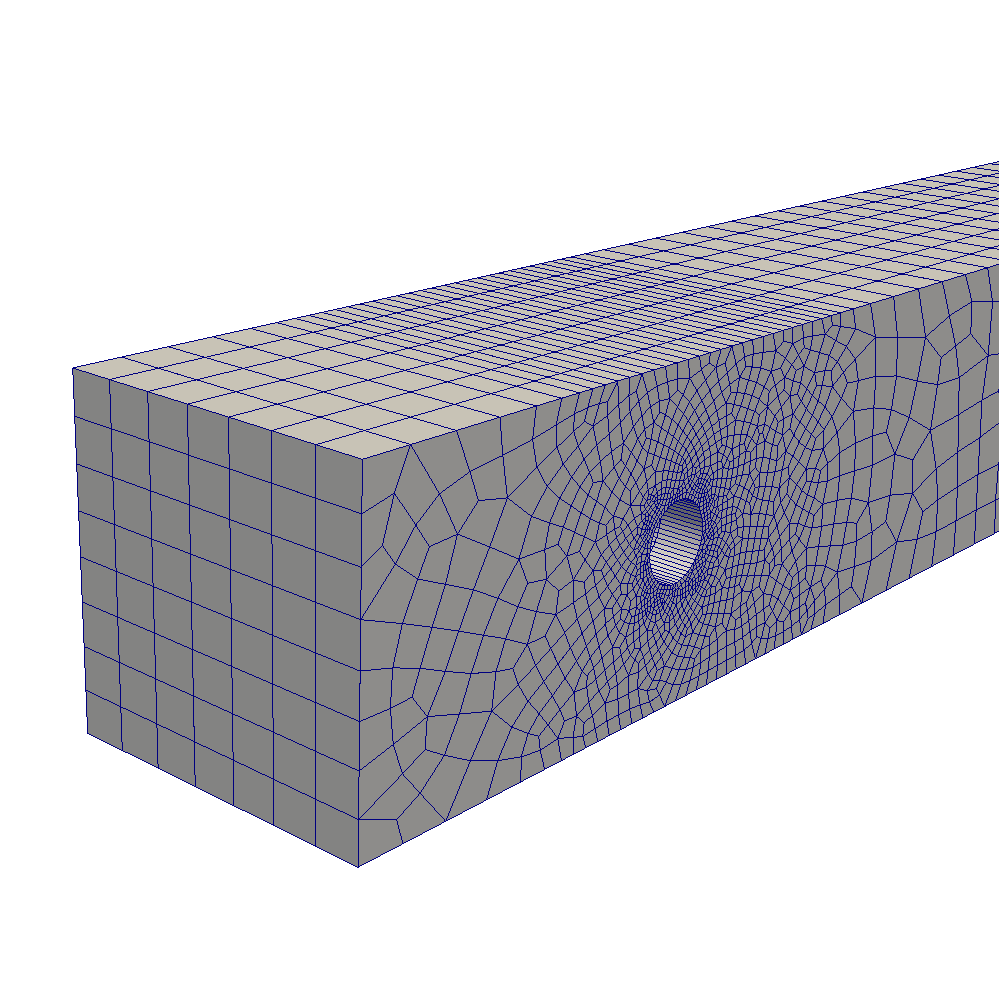}
      \end{overpic} \quad
\begin{overpic}[width=0.45\textwidth]{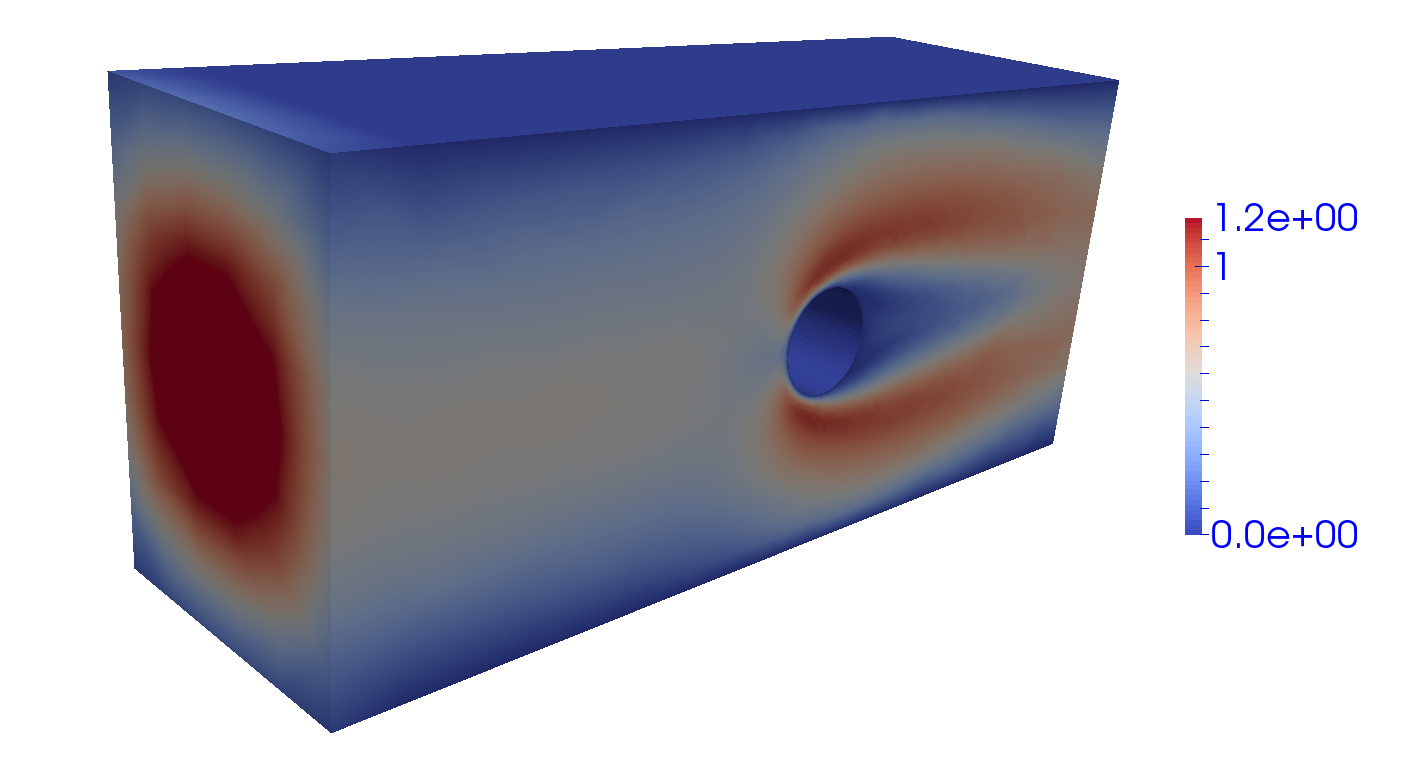}
      \end{overpic}\\
\caption{3D flow past a cylinder: (left) part of the mesh under consideration and (right) illustrative representation of the velocity field for $t = 6$. 
}
\label{fig:MESH3D}
\end{figure}

The 3D benchmark we consider has been studied for the first time in \cite{turek1996} 
and further investigated in \cite{Bayraktar2012, John2006}.
The computational domain is a 2.5 $\times$ 0.41 $\times$ 0.41 parallelepiped with a 
cylinder whose axis is parallel to the $z$-axis and center is located at (0.5, 0.2)
when taking the bottom left corner of the channel as the origin of the axes. Fig.~\ref{fig:MESH3D} (left) 
shows part of the computational domain. The channel is filled with fluid with density $\rho = 1$ 
and viscosity $\mu = 0.001$. We impose a no slip boundary condition on the channel
walls and on the cylinder. At the inflow, we prescribe the following velocity profile:
\begin{align}\label{eq:cyl_bc3d}
\u(0,y,z,t) = \left(\dfrac{36}{0.41^4} \sin\left(\pi t/8 \right) y z \left(0.41 - y \right) \left(0.41 - z \right), 0, 0\right), \quad y,z \in [0, 0.41], \quad t \in (0, 8].
\end{align}
In addition, on the channel walls, cylinder, and at the inlet we impose ${\partial p}/{\partial \n} = 0$ 
where $\n$ is the outward normal. At the outflow, we prescribe $\nabla \u \cdot  \n = 0$ and $p = 0$. 
Note that the Reynolds number is time dependent, with $0 \leq Re \leq 100$ \cite{turek1996, Bayraktar2012, John2006}. 
Like for the 2D benchmark, we start the simulations from fluid at rest.

We consider a hexahedral grid with $h_{min} = 9e-3$, $h_{avg} = 3.75e-2$, and $h_{max} = 6.6e-2$ and a total of 1.07e4 cells. 
The mesh features very low values of maximum non-orthogonality (34$^\circ$), average non-orthogonality (7$^\circ$), skewnwss (0.6), 
and maximum aspect ratio (25). In addition, the mesh is refined next to the cylinder, like the meshes used in \cite{Bayraktar2012, John2006}.
However, notice that this level of refinement is very far from the one required by a DNS \cite{Bayraktar2012, John2006}. Here, the main goal is to show that our ROM approach works well in 3D configurations too. We skip the 
validation of the EFR algorithm at FOM level, since it has already been carried out for the 2D case. 

Like in the 2D case, we use a second-order accurate Central Differencing (CD) 
scheme \cite{Lax1960} for the discretization of the convective term. We set $\Delta t = 5e-3$ \cite{Girfoglio_JCP}, $\alpha = h_{avg}$ (as 
for the 2D test case) and $\chi=\Delta t$ \cite{layton_CMAME}. The main difference with respect to the
2D test case is that the flow field does not exhibit an oscillatory pattern, as shown in Figure \ref{fig:MESH3D} (right). 
Thus, we perform only one numerical experiment and test the performances of our ROM approach 
over the entire time window of interest [0 8] only. 

\begin{figure}[h]
\centering
\includegraphics[height=0.35\textwidth]{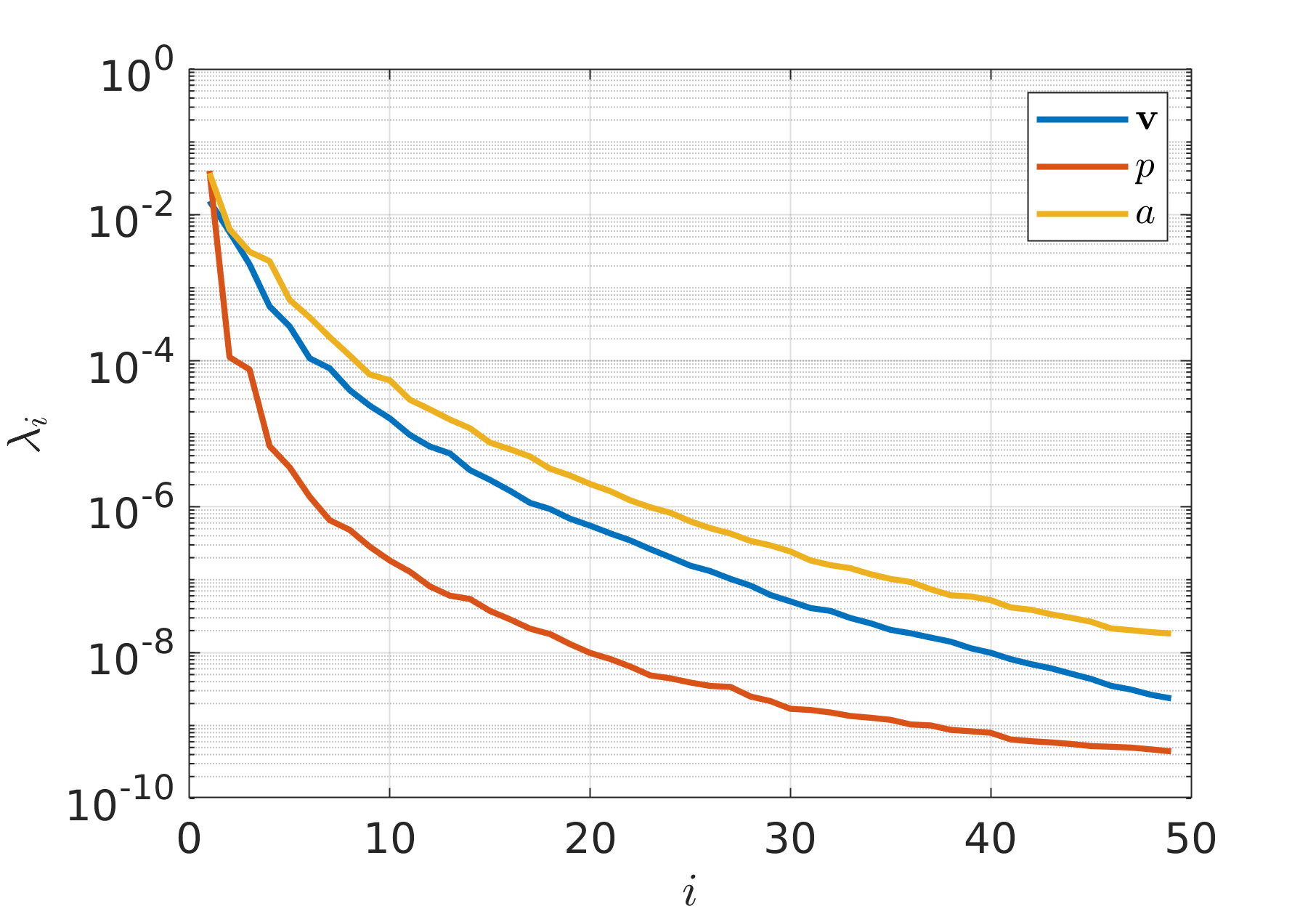}
\caption{3D flow past a cylinder: eigenvalue decay for the intermediate
velocity, pressure, and indicator function.}
\label{fig:eig_lift_3D}
\end{figure}

We collect 400 FOM snapshots for the training in the offline phase, one every 0.02 s (equispaced grid in time). 
Fig.~\ref{fig:eig_lift_3D} shows the eigenvalues decay for velocity $\v$, pressure and indicator function. By a comparison with Fig. \ref{fig:eig_lift}, we see that the decay is faster for all the variables in the 3D case. Therefore, a smaller number 
of basis functions needs to be considered. This is due to the fact that the 2D flow is more complex. 

\begin{figure}[h]
\centering
 \begin{overpic}[width=0.45\textwidth]{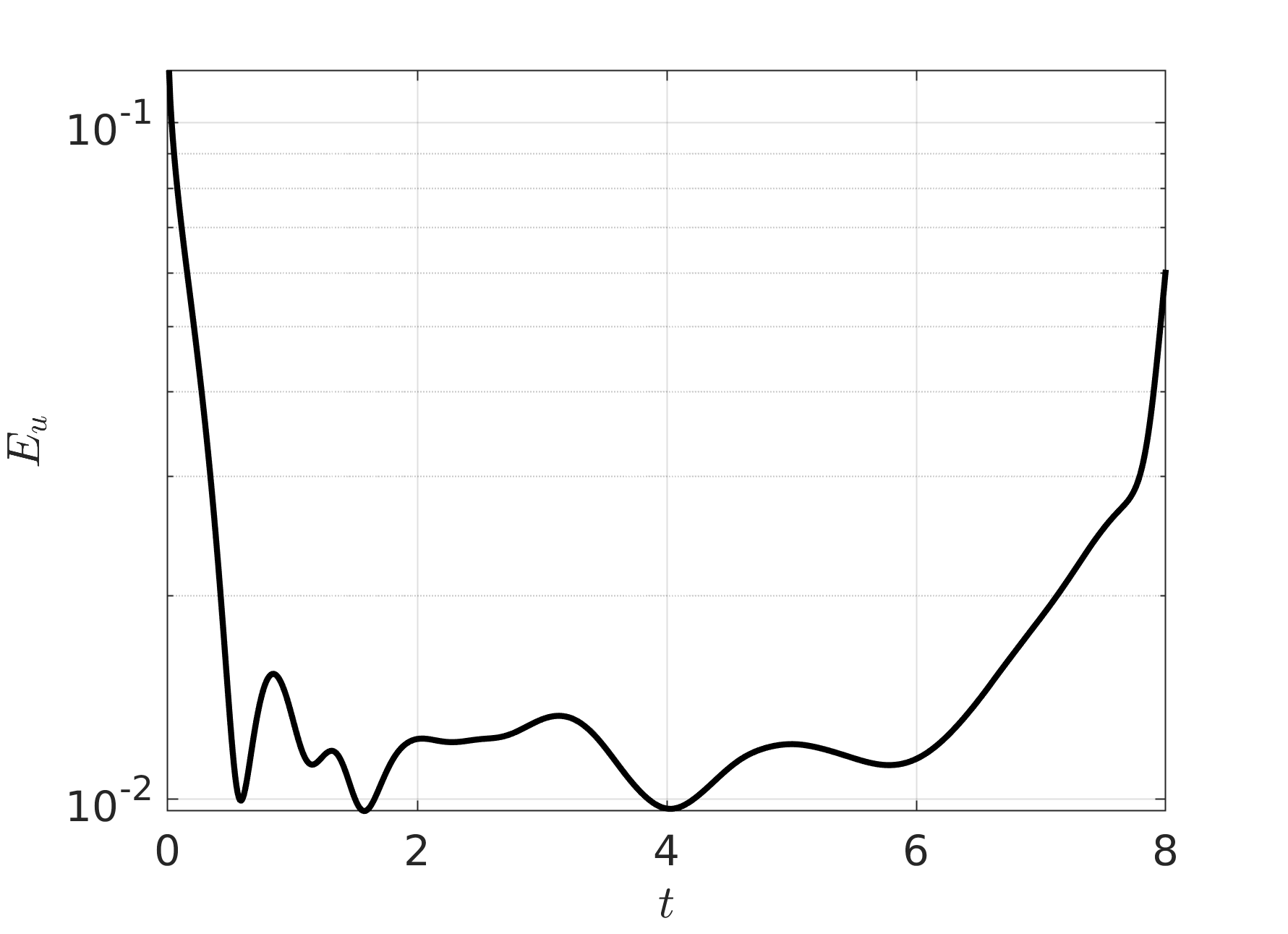}
      \end{overpic}
\begin{overpic}[width=0.45\textwidth]{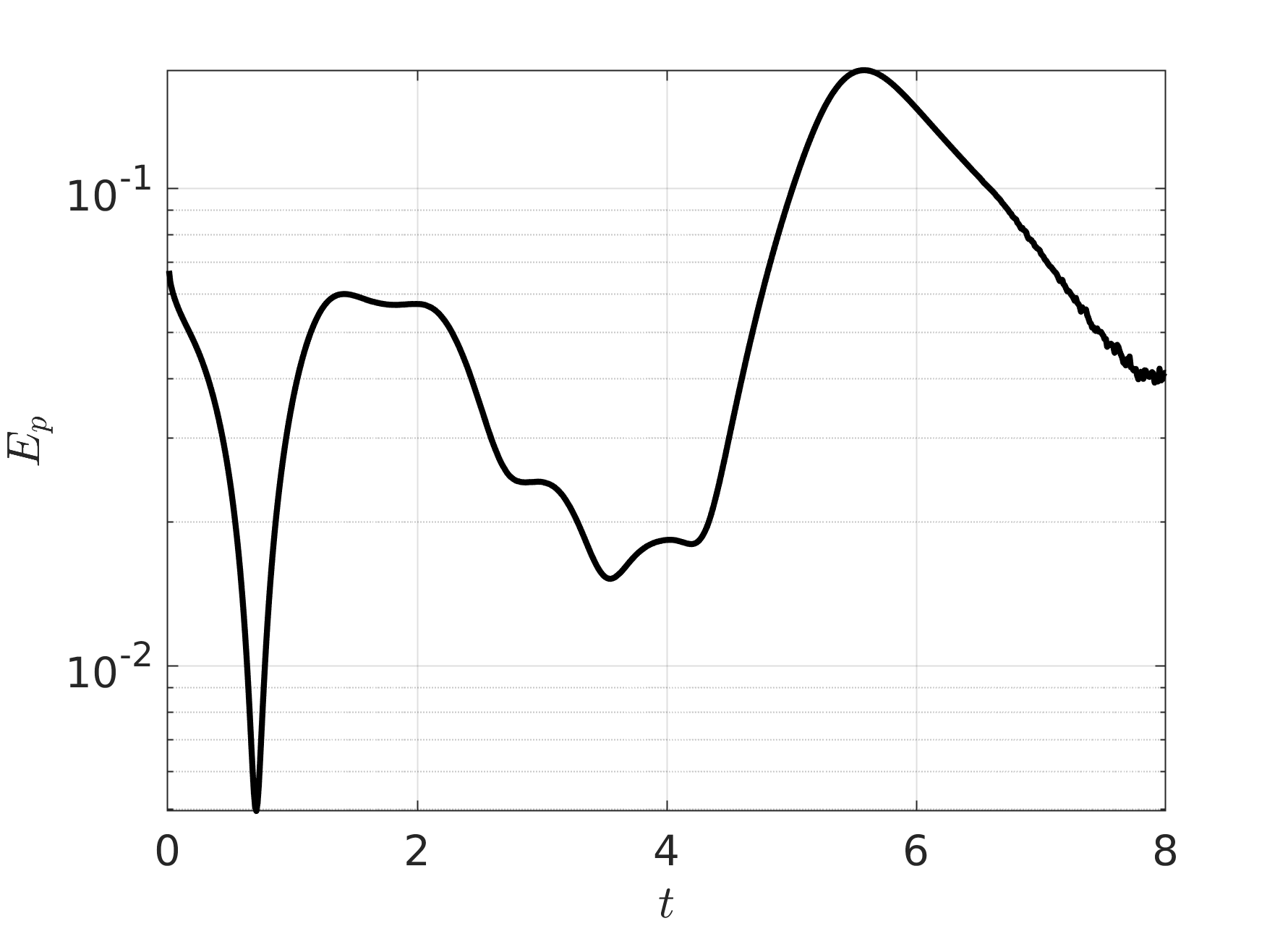}
      \end{overpic}\\ 
\begin{overpic}[width=0.45\textwidth]{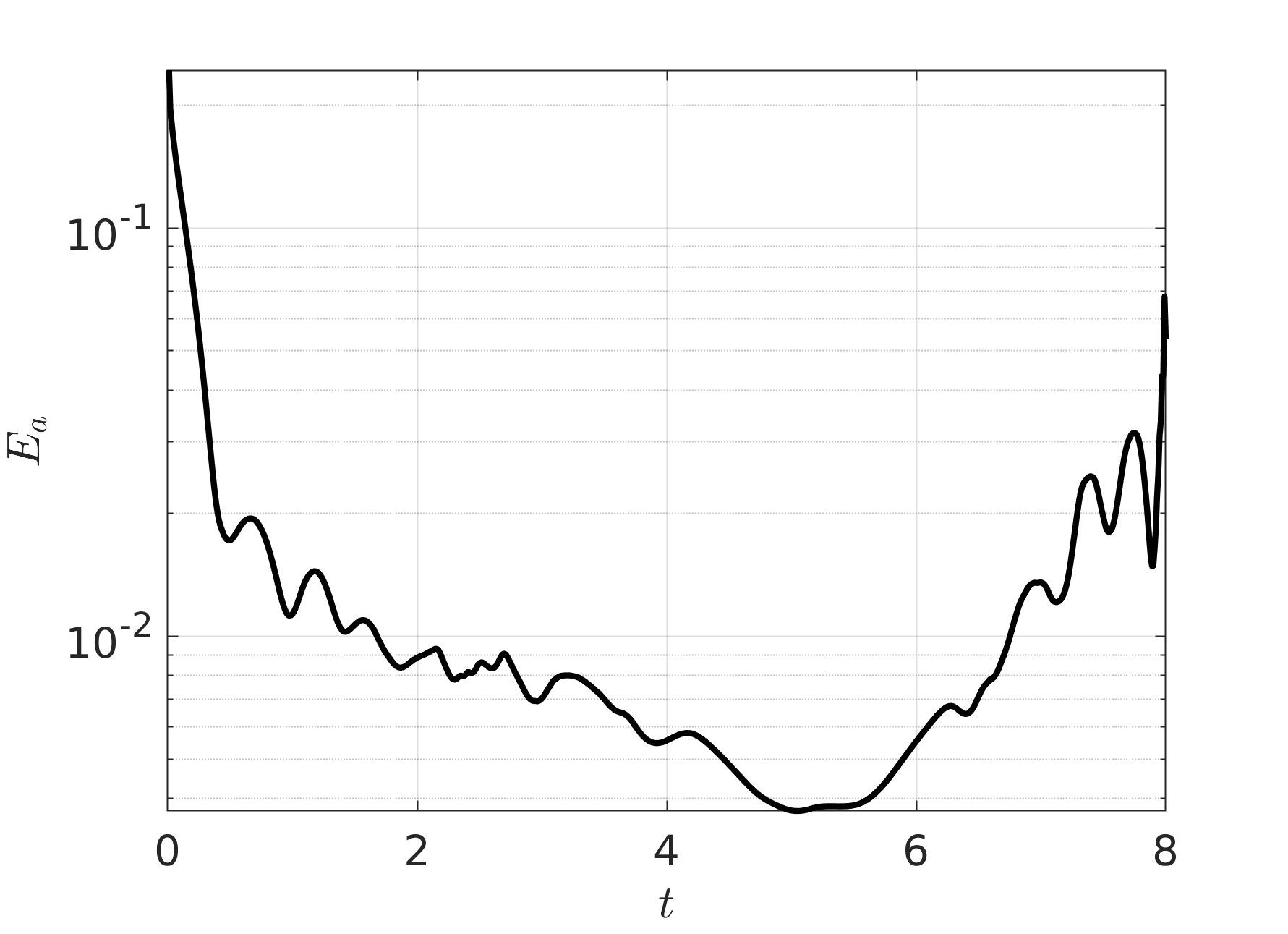}
     \end{overpic}\\ 
\caption{3D flow past a cylinder: time history of $L^2$ norm of the relative error \eqref{eq:error1} for velocity $\u$ (top left), 
pressure field (top right), and indicator function field (bottom).}
\label{fig:err_3D}
\end{figure}

For the ROM simulations, we collected a set of 800 equispaced temporal instants (i.e., one every 0.01 s), 
which includes both samples used in the offline stage and samples in between. 
In order to retain 99.99\% of the snapshots energy, we needs 9 modes for $\v$, 4 modes for $p$ and 12 modes for $a$. 
Fig.~\ref{fig:err_3D} shows errors \eqref{eq:error1} and Table \ref{tab:errors3D} reports
minimum, average, and maximum relative errors. 
The relative errors for the velocity $\u$ and the indicator functions stay below $10^{-1}$ for most of the time interval, 
expect at the beginning and the end of the simulation, i.e. when the flow pattern is affected by transient effects.
As for the pressure, we see a relative error larger than $10^{-1}$ at $t \approx 6$. 
These results are qualitatively similar to the ones showed by a Leray model in \cite{Girfoglio_JCP}. 
This is expected, since the EFR algorithm can be seen as a splitting scheme for the Leray model \cite{BQV}. 

\begin{table}[h]
\centering
\begin{tabular}{lccc}
\multicolumn{3}{c}{} \\
\cline{1-4}
 & $\u$ & $p$ & $a$\\ 
\hline
Maximum $E_\Phi$ & 0.12  & 0.18  & 0.24 \\
Average $E_\Phi$ & 0.02  & 0.06  & 0.01 \\
Miminum $E_\Phi$ & 0.01  & 0.005 & 0.004 \\
\hline
\end{tabular}
\caption{3D flow past a cylinder: 
maximum, average, and minimum relative $L^2$ errors \eqref{eq:error1} for the end-of-step velocity, pressure, and indicator function.}
\label{tab:errors3D}
\end{table}

\begin{figure}
\centering
 \begin{overpic}[width=0.44\textwidth]{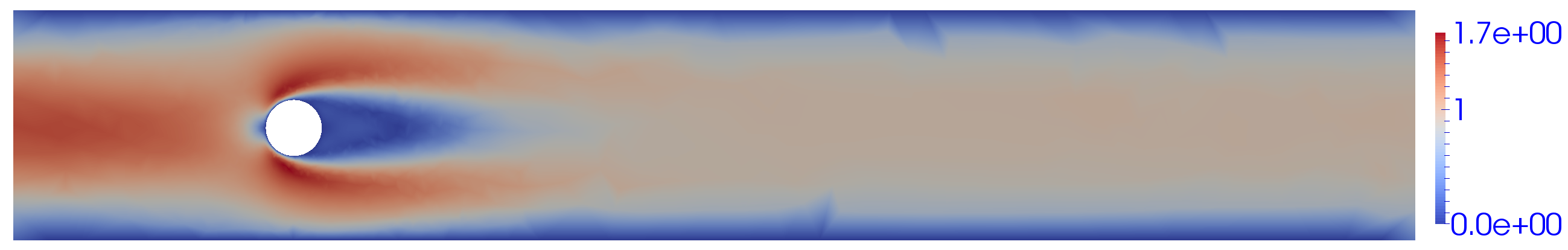}
        \put(35,18){FOM}
        \put(-8,7){$\u$}
      \end{overpic}
 \begin{overpic}[width=0.44\textwidth]{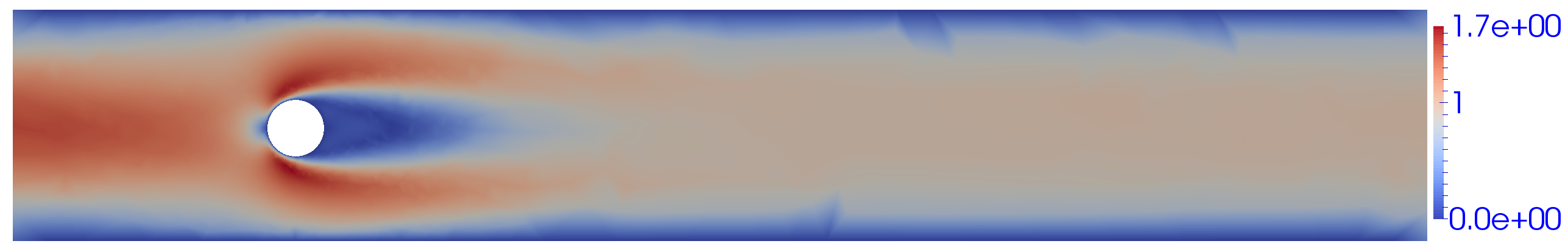}
        \put(35,18){ROM}
      \end{overpic}\\
      \vskip .2cm
 \begin{overpic}[width=0.44\textwidth]{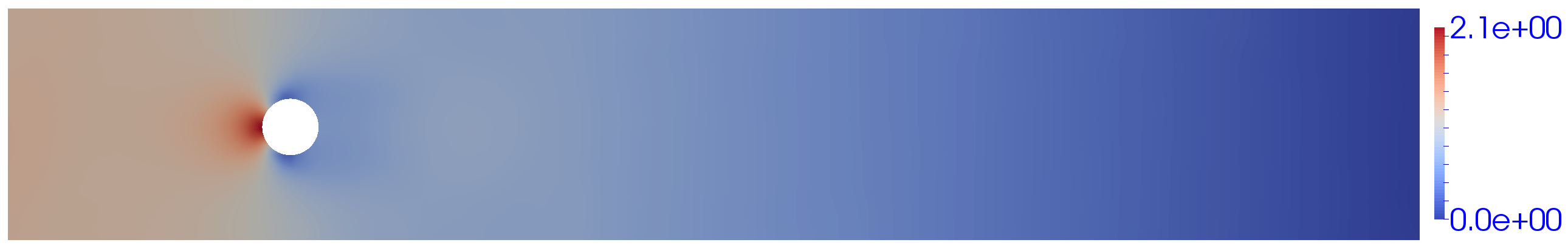}
        \put(-8,7){$p$}
      \end{overpic}
 \begin{overpic}[width=0.44\textwidth]{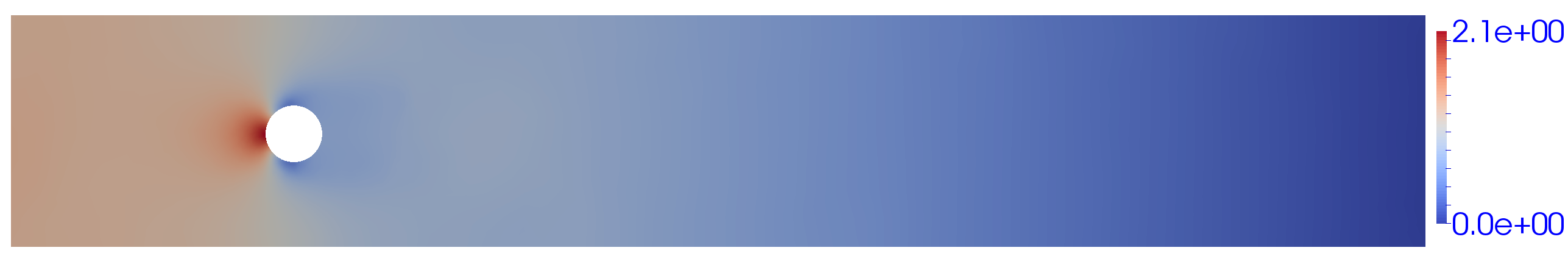}
      \end{overpic}\\
       \vskip .3cm
 \begin{overpic}[width=0.44\textwidth]{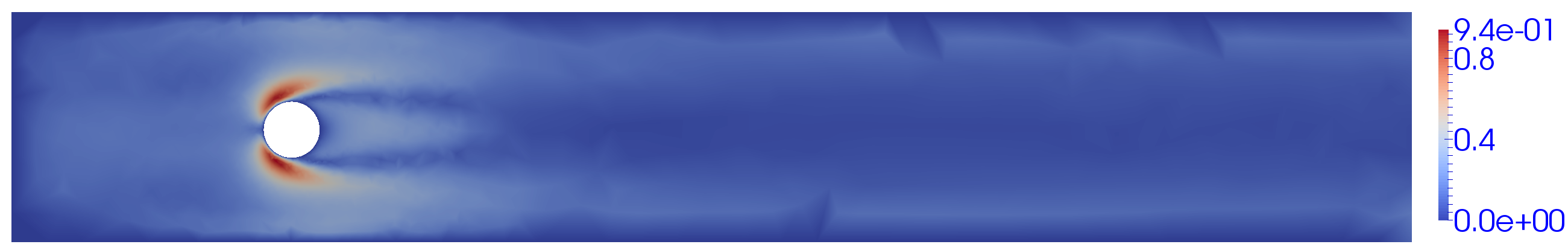}
        \put(-8,7){$a$}
      \end{overpic}
 \begin{overpic}[width=0.44\textwidth]{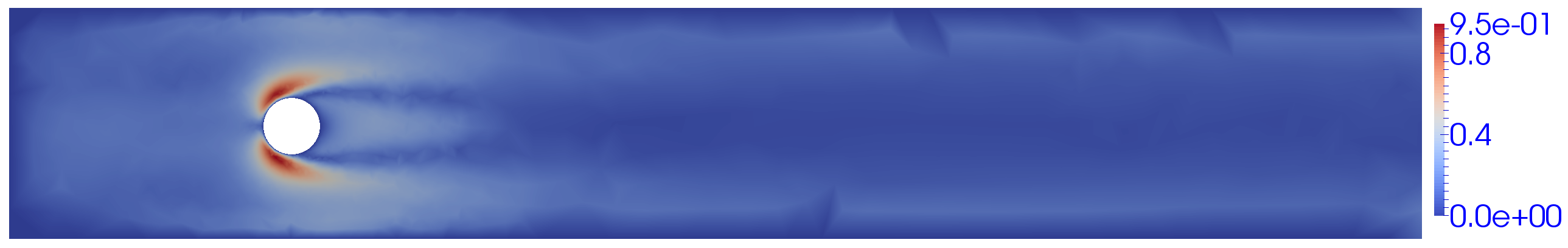}
      \end{overpic}
\caption{3D flow past a cylinder: 
comparison between 
velocity $\u$ (first raw), pressure (second raw), and indicator function (3nd raw)
computed by the FOM (left) and the ROM (right) on the midsection ($z = 0.205$) at $t = 1.9$.}
\label{fig:comp_t_1_9_3D}
\end{figure}

\begin{figure}
\centering
 \begin{overpic}[width=0.44\textwidth]{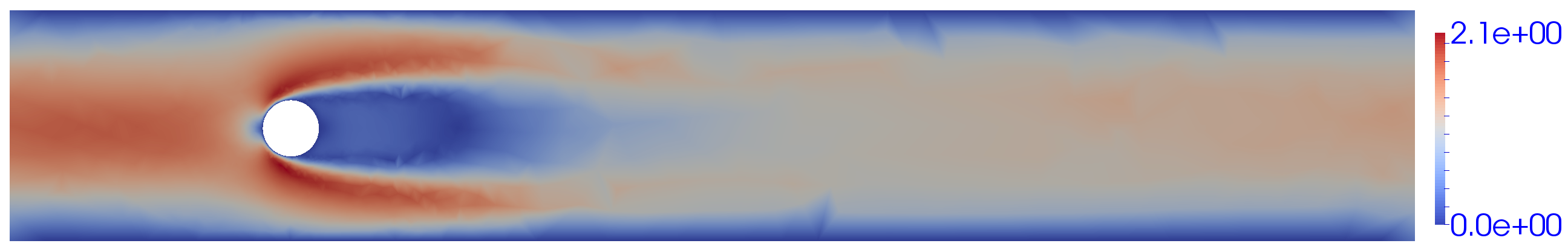}
         \put(35,18){FOM}
        \put(-8,7){$\u$}
      \end{overpic}
 \begin{overpic}[width=0.44\textwidth]{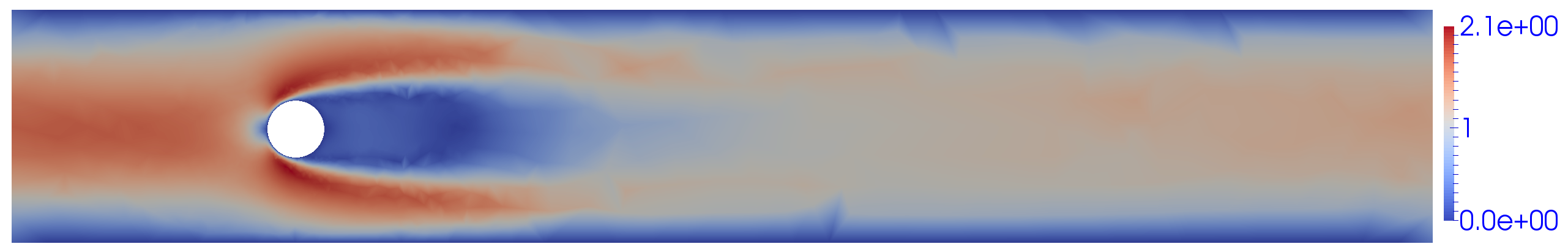}
        \put(35,18){ROM}
      \end{overpic}\\
        \vskip .3cm
 \begin{overpic}[width=0.44\textwidth]{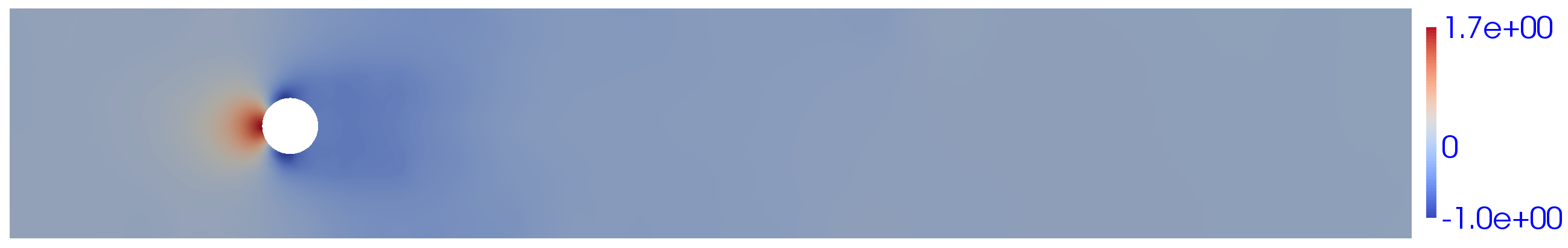}
        \put(-8,7){$p$}
      \end{overpic}
 \begin{overpic}[width=0.44\textwidth]{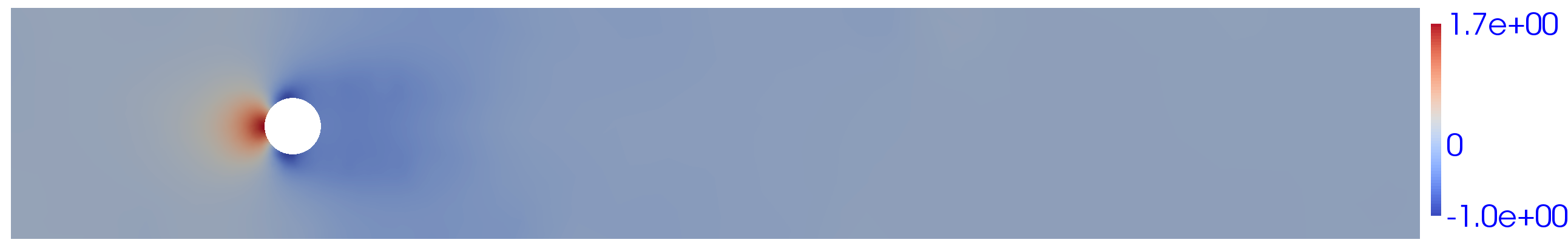}
      \end{overpic}\\
      \vskip .2cm
 \begin{overpic}[width=0.44\textwidth]{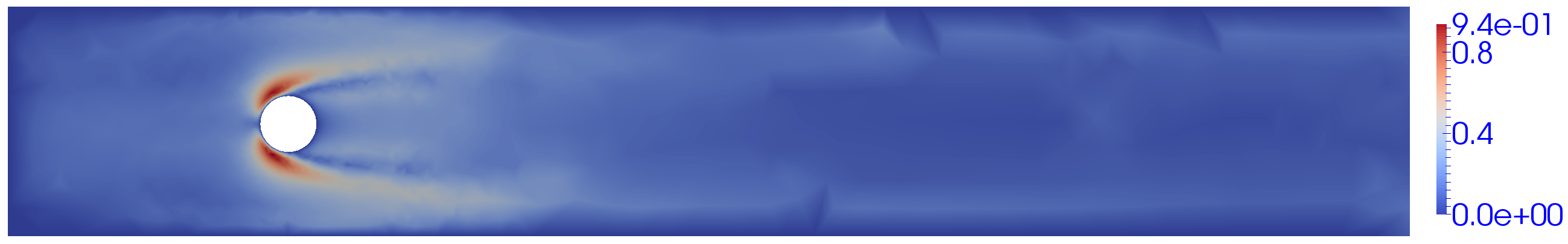}
        \put(-8,7){$a$}
      \end{overpic}
 \begin{overpic}[width=0.44\textwidth]{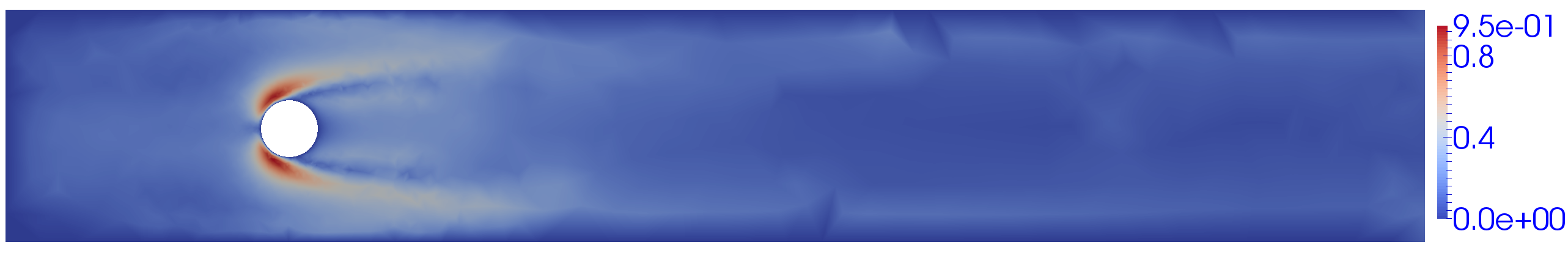}
      \end{overpic}
\caption{3D flow past a cylinder: comparison between 
velocity $\u$ (first raw), pressure (second raw), and indicator function (3nd raw)
computed by the FOM (left) and the ROM (right) on the midsection ($z = 0.205$) at $t = 5.5$.
}
\label{fig:comp_t_5_5_3D}
\end{figure}

\begin{figure}
\centering
 \begin{overpic}[width=0.44\textwidth]{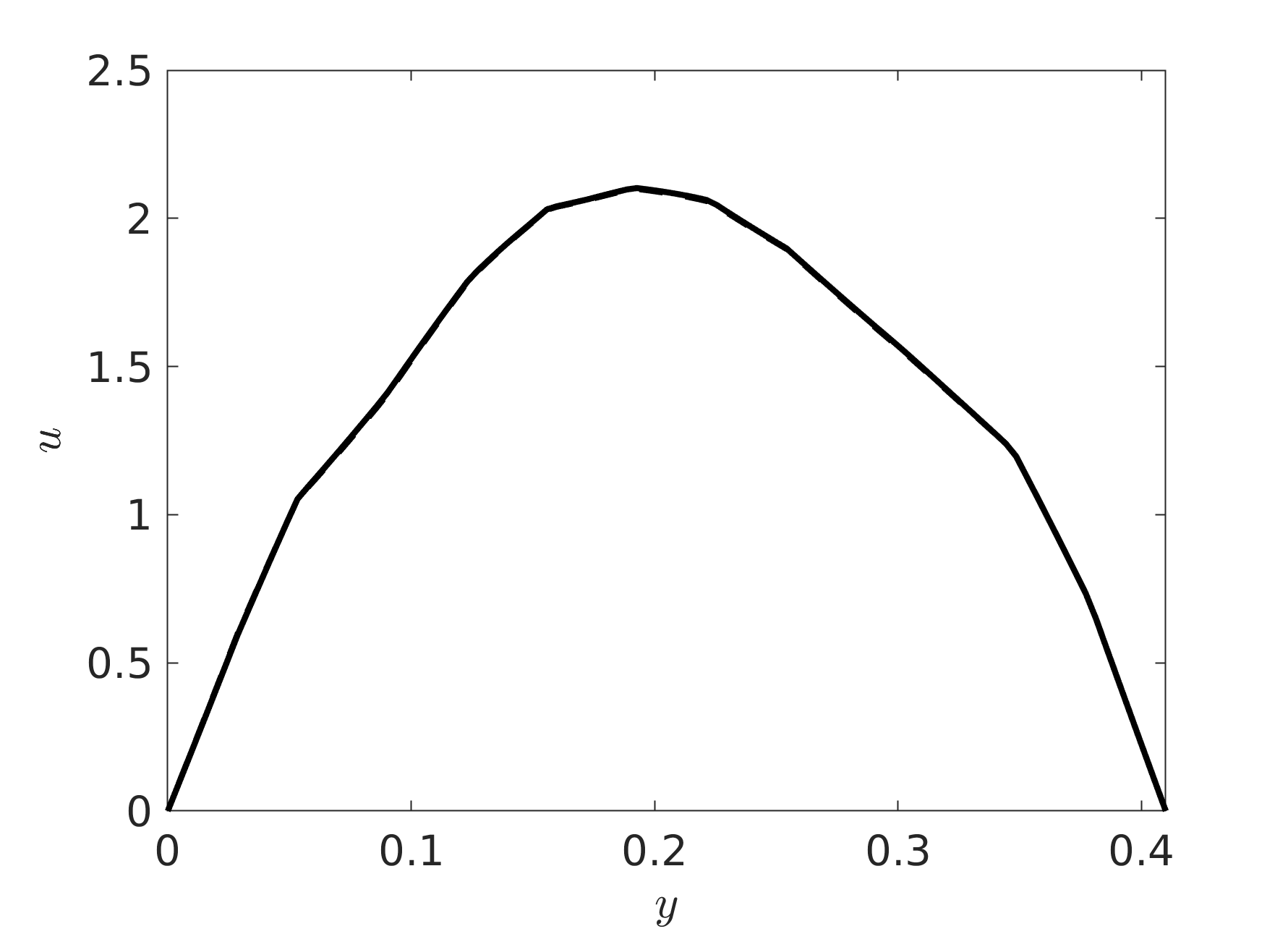}
      \end{overpic}
 \begin{overpic}[width=0.44\textwidth]{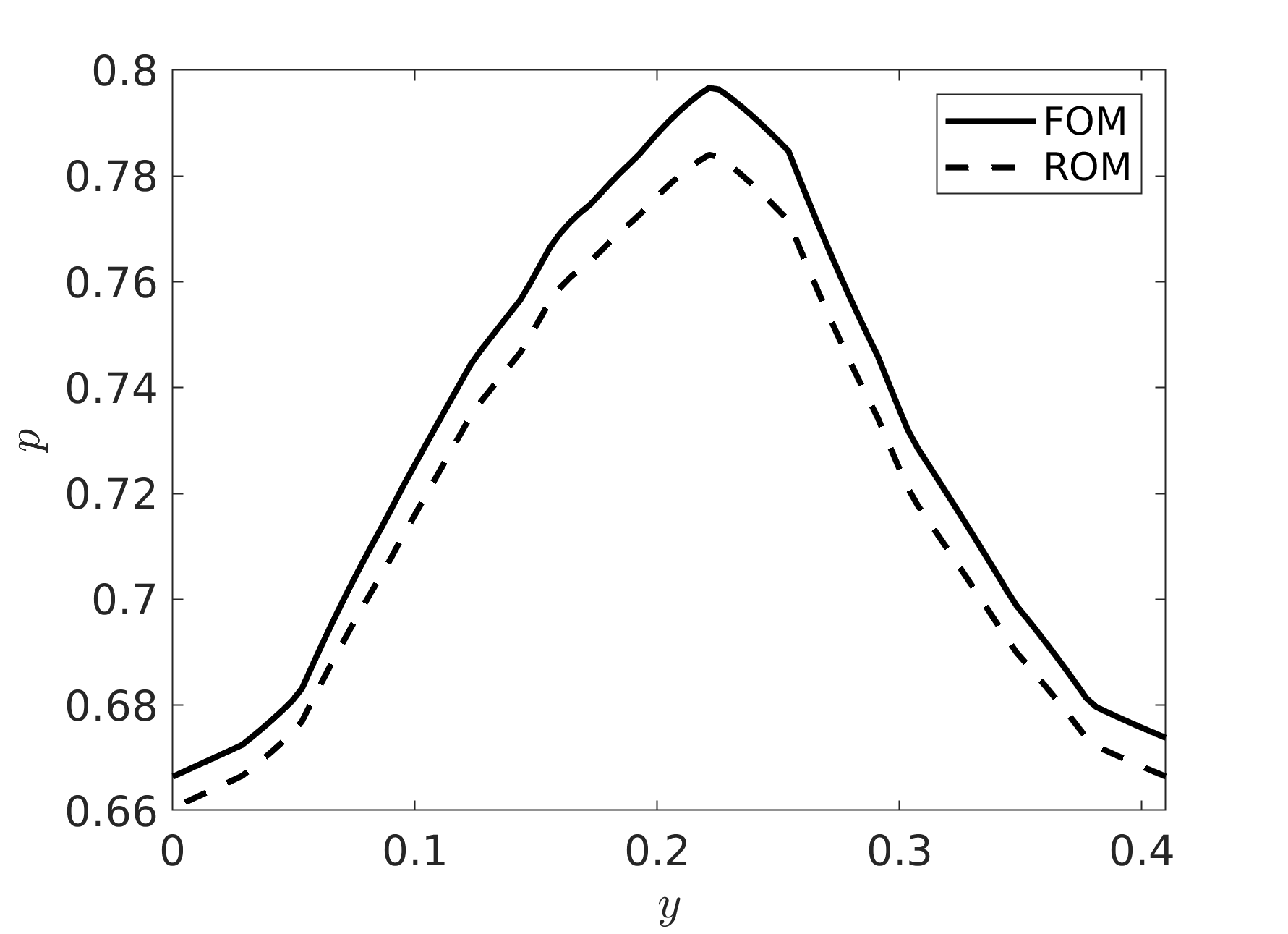}
      \end{overpic} \\
 \begin{overpic}[width=0.44\textwidth]{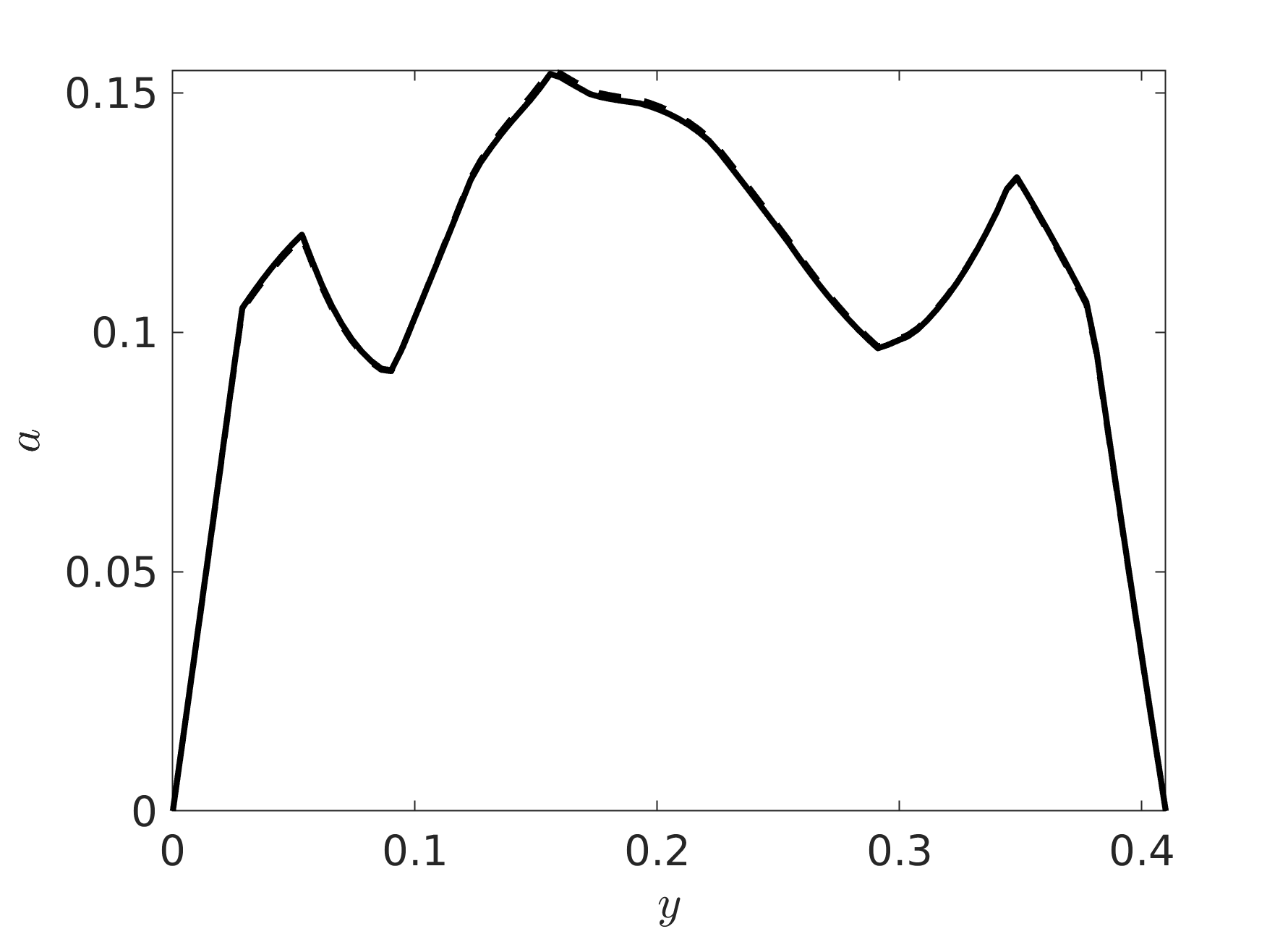}
      \end{overpic}\\
\caption{3D flow past a cylinder: velocity (top left), pressure (top rig), and indicator function (bottom) profiles computed by FOM and ROM
along line $x = 0.25$, $y \in$ [0 0.41], $z = 0.205$ at $t = 4$.}
\label{fig:comp_sec_3D_0_25}
\end{figure}

\begin{figure}
\centering
 \begin{overpic}[width=0.44\textwidth]{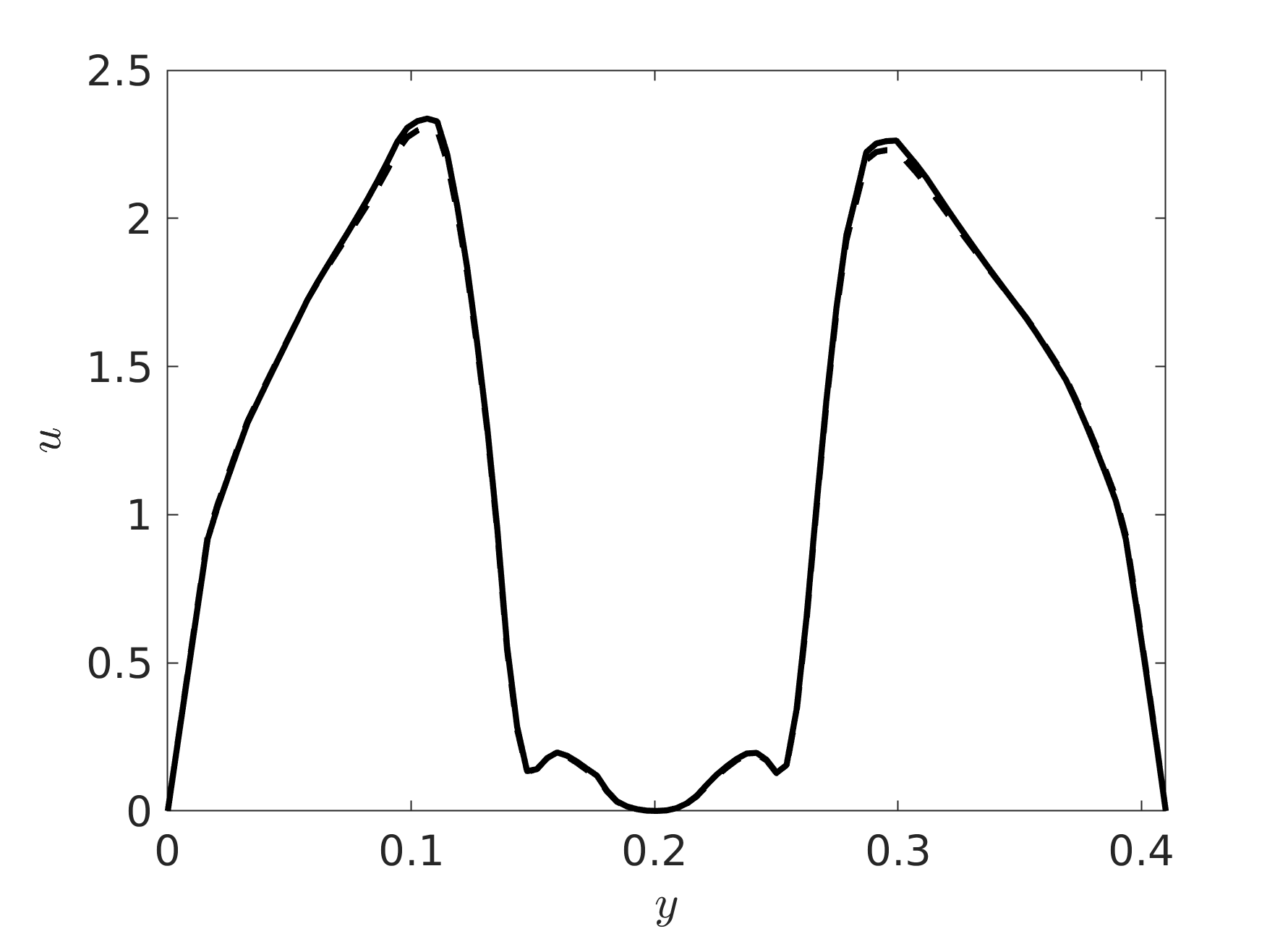}
      \end{overpic}
 \begin{overpic}[width=0.44\textwidth]{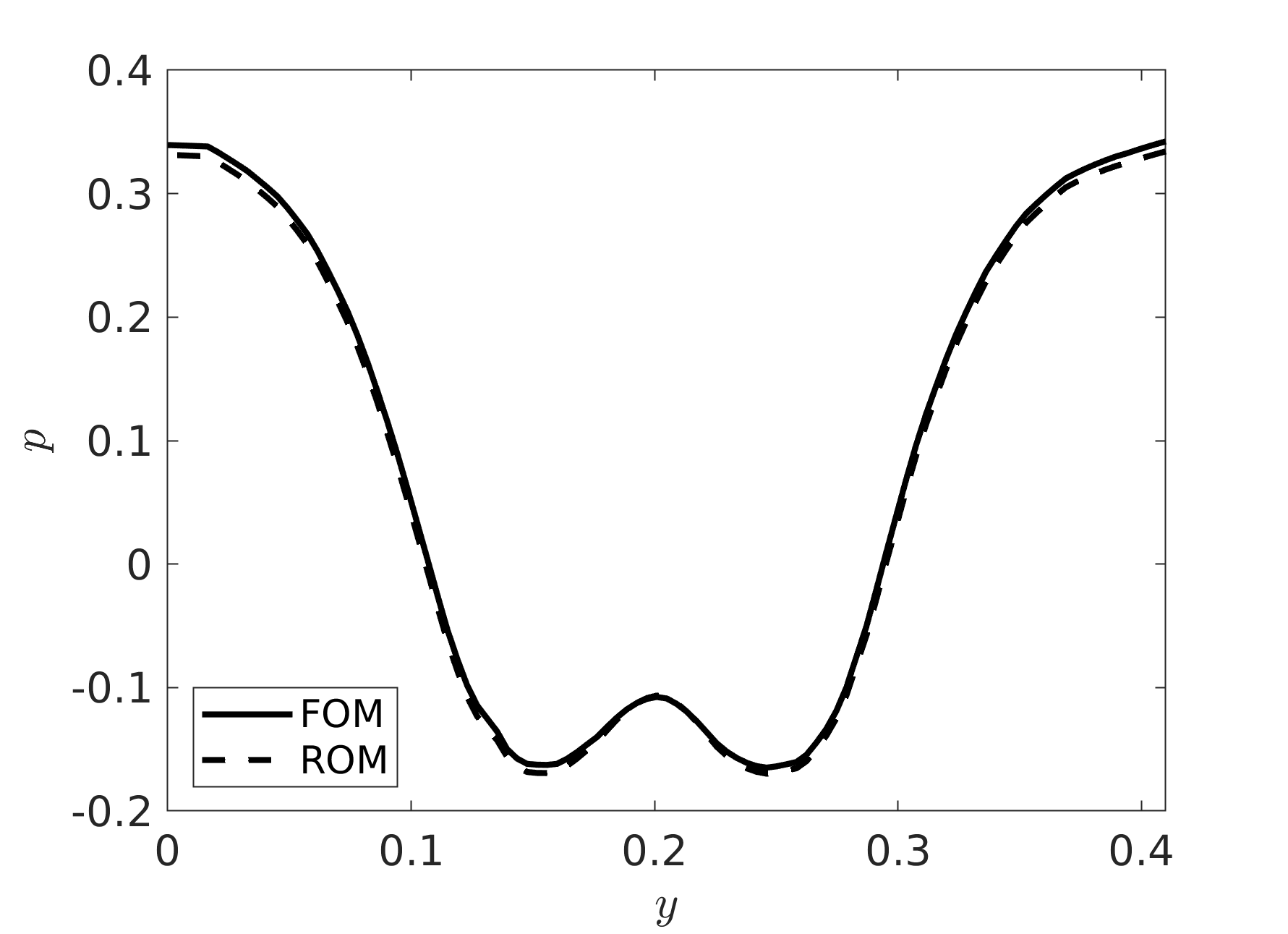}
      \end{overpic}\\
 \begin{overpic}[width=0.44\textwidth]{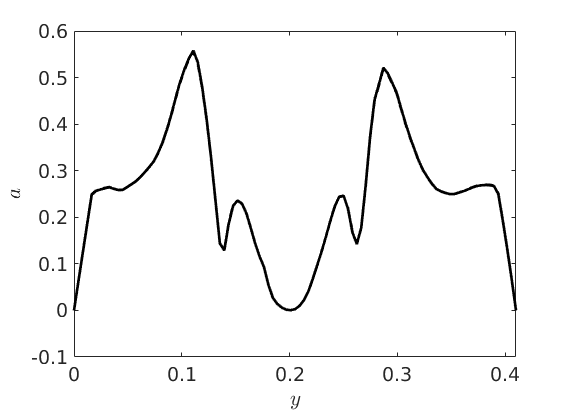}
      \end{overpic}\\
\caption{3D flow past a cylinder: velocity (top left), pressure (top rig), and indicator function (bottom) profiles computed by FOM and ROM
along line $x = 0.55$, $y \in$ [0 0.41], $z = 0.205$ at $t = 4$.}
\label{fig:comp_sec_3D_0_55}
\end{figure}

Figures \ref{fig:comp_t_1_9_3D} and \ref{fig:comp_t_5_5_3D} display a qualitative comparison 
for velocity $\u$, pressure, and indicator function computed by FOM and ROM 
on the midsection ($z = 0.205$) at times $t = 1.9$ and $t = 5.5$, respectively. 
Our ROM provides a good reconstruction of all the variables at both times.
For a further comparison, Figures \ref{fig:comp_sec_3D_0_25} and \ref{fig:comp_sec_3D_0_55} show the profiles of all the variables at 
$t = 4$ along a line upstream of the cylinder ($x = 0.25$, $y \in$ [0 0.41], $z = 0.205$)
and a line downstream of the cylinder ($x = 0.55$, $y \in$ [0 0.41], $z = 0.205$), respectively. 
We see great agreement between FOM and ROM profiles for all the variables with the exception 
of the pressure in Fig. \ref{fig:comp_sec_3D_0_25} (top right), for which we observe a small difference.

\begin{figure}
\centering
 \begin{overpic}[width=0.45\textwidth]{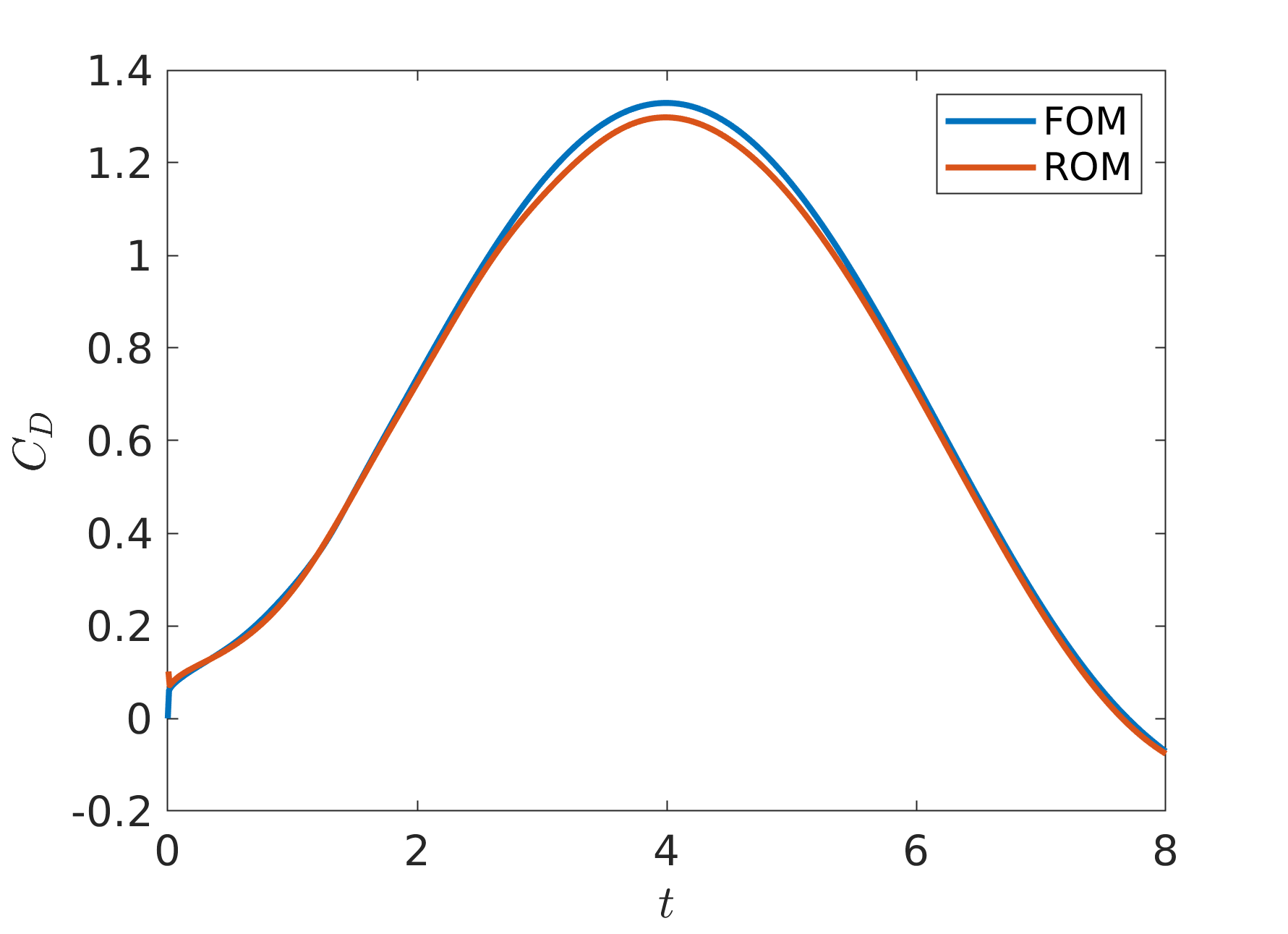}
      \end{overpic}
       \begin{overpic}[width=0.45\textwidth]{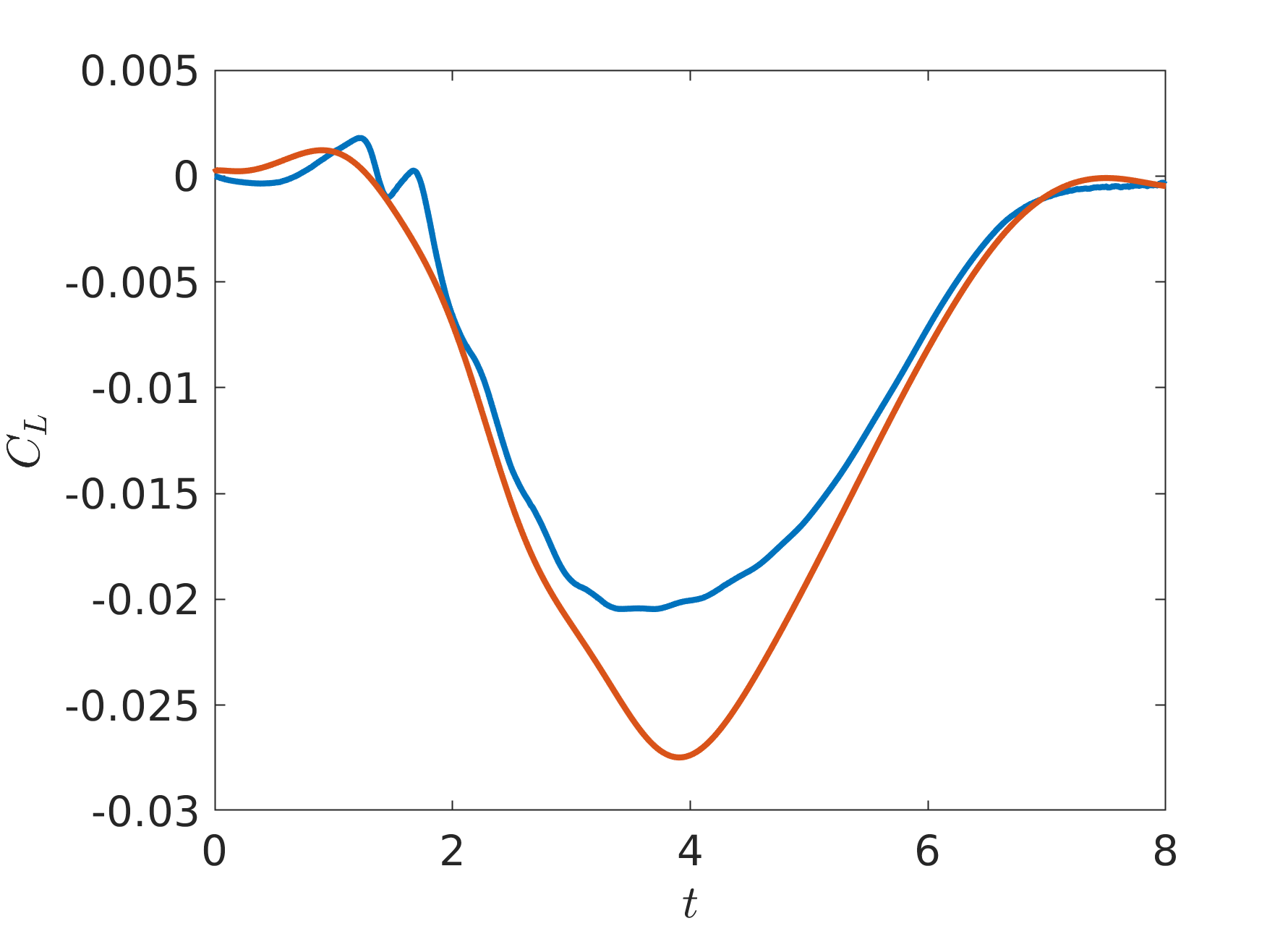}
      \end{overpic}\\
\caption{3D flow past a cylinder: aerodynamic coefficients $C_D$ (left) and $C_L$ (right) computed by FOM and ROM.}
\label{fig:coeff_3D}
\end{figure}

Fig.~\ref{fig:coeff_3D} reports the drag and lift coefficients computed by FOM and ROM. As for the 2D case, 
we see that $C_D$ is well reconstructed by the ROM while the reconstruction of the time evolution of $C_L$ is not as accurate. 
Error \eqref{eq:error_coeff} is $E_{C_L} = 0.25$. 

Finally, we comment on the computational costs. The total CPU time required by a FOM simulation is 460 s. 
Our ROM approach takes about 2 s for the solution of the reduced algebraic systems.
So we obtain a speed-up is about 230, 
which is much higher than the speed-up obtained for the 2D test because in the 3D test we retain a much lower number of modes.

\section{Conclusions and perspectives}\label{sec:conclusions}

This work presents an hybrid Reduced Order Method (ROM) for the Evolve-Filter-Relax (EFR) 
algorithm with a nonlinear indicator function: we exploit a data-driven reduction method based on interpolation
with Radial Basis Functions to approximate the indicator function and a classical POD-Galerkin projection approach 
for the reconstruction of the velocity and the pressure fields. 
This mixed strategy has been recently used for the RANS equations
and here we have extended it to  a LES framework. 
To the best our knowledge, it is the first time that a ROM nonlinear differential filter 
(i.e. a ROM spatial filter that uses an explicit lengthscale) is proposed.
We assessed our ROM approach through two classical benchmarks: 2D and 3D flow past a cylinder. 
We found that our ROM can capture the flow features and the evolution of the aerodynamics coeffiecients
with good accuracy when compared to the full order model. 
In addition, we showed the our approach is more computationally efficient in the 3D test. 

Moving forward, we will test our ROM approach with more realistic applications involving flow at 
larger Reynolds numbers (a few thousands) like, e.g., the FDA benchmark (see \cite{Girfoglio2019}). 
Moreover, we plan to run a parametric study for the viscosity and for key model parameters, such as $\delta$ and $\chi$.
This would help us understand how to set the model parameter to obtain the most accurate results
when compared to DNS.

\section{Acknowledgements}\label{sec:acknowledgements}
We acknowledge the support provided by the European Research Council Executive Agency by the Consolidator Grant project AROMA-CFD ``Advanced Reduced Order Methods with Applications in Computational Fluid Dynamics" - GA 681447, H2020-ERC CoG 2015 AROMA-CFD, PI G. Rozza, and INdAM-GNCS 2019-2020 projects.
This work was also partially supported by US National Science Foundation through grant DMS-1620384
and DMS-195353.

\bibliographystyle{plain}
\bibliography{ROM_NL_GQR}

\begin{thebibliography}{10}

\bibitem{akhtar2009stability}
I.~Akhtar, A.~H. Nayfeh, and C.~J. Ribbens.
\newblock On the stability and extension of reduced-order {Galerkin} models in
  incompressible flows.
\newblock {\em Theoretical and Computational Fluid Dynamics}, 23(3):213--237,
  2009.

\bibitem{Aubry1988}
N.~Aubry, P.~Holmes, J.~L. Lumley, and E.~Stone.
\newblock The dynamics of coherent structures in the wall region of a turbulent
  boundary layer.
\newblock {\em Journal of Fluid Mechanics}, 192:115--173, 1988.

\bibitem{Bader2016}
E.~Bader, M.K\"{a}rcher, M.~A. Grepl, and K.~Veroy.
\newblock {Certified Reduced Basis Methods for Parametrized Elliptic Optimal
  Control Problems with Distributed Controls}.
\newblock {\em {SIAM} Journal on Scientific Computing}, 75:276--307, 2018.

\bibitem{barrault04:_empir_inter_method}
M.~Barrault, N.~C. Nguyen, Y.~Maday, and A.~T. Patera.
\newblock An ``empirical interpolation'' method: Application to efficient
  reduced-basis discretization of partial differential equations.
\newblock {\em Comptes Rendus Mathematique}, 339:667--672, 2004.

\bibitem{Bayraktar2012}
E.~Bayraktar, O.~Mierka, and S.~Turek.
\newblock {Benchmark computations of 3D laminar flow around a cylinder with
  CFX, OpenFOAM and FeatFlow}.
\newblock {\em International Journal of Computational Science and Engineering},
  7:253 -- 266, 2012.

\bibitem{BAZILEVS2007173}
Y.~Bazilevs, V.M. Calo, J.A. Cottrell, T.J.R. Hughes, A.~Reali, and
  G.~Scovazzi.
\newblock Variational multiscale residual-based turbulence modeling for large
  eddy simulation of incompressible flows.
\newblock {\em Computer Methods in Applied Mechanics and Engineering},
  197(1):173--201, 2007.

\bibitem{Benner2015}
P.~Benner, S.~Gugercin, and K.~Willcox.
\newblock {A Survey of Projection-Based Model Reduction Methods for Parametric
  Dynamical Systems}.
\newblock {\em {SIAM} Review}, 57(4):483--531, 2015.

\bibitem{bennerParSys}
P.~Benner, M.~Ohlberger, A.~Patera, and K.~Rozza, G.and~Urban.
\newblock {\em {Model Reduction of Parametrized Systems}}, volume 1st ed. 2017
  of {\em MS\&A series}.
\newblock Springer, 2017.

\bibitem{ModelOrderReduction}
P.~Benner, W.~Schilders, S.~Grivet-Talocia, A.~Quarteroni, G.~Rozza, and L.~M.
  Silveira.
\newblock {\em Model Order Reduction}.
\newblock De Gruyter, Berlin, Boston, 2020.

\bibitem{BQV}
L.~Bertagna, A.~Quaini, and A.~Veneziani.
\newblock {Deconvolution-based nonlinear filtering for incompressible flows at
  moderately large Reynolds numbers}.
\newblock {\em International Journal for Numerical Methods in Fluids},
  81(8):463--488, 2016.

\bibitem{Borggaard2009}
J.~Borggaard, T.~Iliescu, and J.P. Roop.
\newblock A bounded artificial viscosity large eddy simulation model.
\newblock {\em SIAM Journal on Numerical Analysis}, 47:622--645, 2009.

\bibitem{Bowers2012}
A.~L. Bowers, L.~G. Rebholz, A.~Takhirov, and C.~Trenchea.
\newblock Improved accuracy in regularization models of incompressible flow via
  adaptive nonlinear filtering.
\newblock {\em International Journal for Numerical Methods in Fluids},
  70(7):805--828, 2012.

\bibitem{abigail_CMAME}
A.L. Bowers and L.G. Rebholz.
\newblock Numerical study of a regularization model for incompressible flow
  with deconvolution-based adaptive nonlinear filtering.
\newblock {\em Computer Methods in Applied Mechanics and Engineering},
  258:1--12, 2013.

\bibitem{Boyd1998283}
J.~P. Boyd.
\newblock Two comments on filtering (artificial viscosity) for {Chebyshev} and
  {Legendre} spectral and spectral element methods: Preserving boundary
  conditions and interpretation of the filter as a diffusion.
\newblock {\em Journal of Computational Physics}, 143(1):283--288, 1998.

\bibitem{Chaturantabut2010}
S.~Chaturantabut and D.C. Sorensen.
\newblock Nonlinear model reduction via discrete empirical interpolation.
\newblock {\em SIAM Journal on Scientific Computing}, 32(5):2737--2764, 2010.

\bibitem{ChinestaEnc2017}
F.~Chinesta, A.~Huerta, G.~Rozza, and K.~Willcox.
\newblock {Model Order Reduction}.
\newblock {\em Encyclopedia of Computational Mechanics, Elsevier Editor}, 2016.

\bibitem{Chinesta2011}
F.~Chinesta, P.~Ladeveze, and E.~Cueto.
\newblock {A Short Review on Model Order Reduction Based on Proper Generalized
  Decomposition}.
\newblock {\em Archives of Computational Methods in Engineering}, 18(4):395,
  2011.

\bibitem{couplet_sagaut_basdevant_2003}
M.~Couplet, P.~Sagaut, and C.~Basdevant.
\newblock {Intermodal energy transfers in a proper orthogonal decomposition
  {Galerkin} representation of a turbulent separated flow}.
\newblock {\em Journal of Fluid Mechanics}, 491:275–284, 2003.

\bibitem{Dumon20111387}
A.~Dumon, C.~Allery, and A.~Ammar.
\newblock {Proper General Decomposition (PGD) for the resolution of
  Navier-Stokes equations}.
\newblock {\em Journal of Computational Physics}, 230(4):1387--1407, 2011.

\bibitem{Dunca2005}
A.~Dunca and Y.~Epshteyn.
\newblock On the {Stolz-Adams} deconvolution model for the large-eddy
  simulation of turbulent flows.
\newblock {\em SIAM Journal on Mathematical Analysis}, 37(6):1890--1902, 2005.

\bibitem{Ervin2010}
V.~Ervin, W.~Layton, and M.~Neda.
\newblock Numerical analysis of filter based stabilization for evolution
  equations.
\newblock {\em SIAM Journal on Numerical Analysis}, 50:2307--2335, 2010.

\bibitem{Fischer2001265}
P.~Fischer and J.~Mullen.
\newblock Filter-based stabilization of spectral element methods.
\newblock {\em Comptes Rendus de l'Academie des Sciences - Series I -
  Mathematics}, 332(3):265--270, 2001.

\bibitem{Girfoglio2021}
M.~Girfoglio, A.~Quaini, and G.Rozza.
\newblock Fluid--structure interaction simulations with a {LES filtering
  approach in solids4Foam}.
\newblock https://arxiv.org/abs/2102.08011, 2021.

\bibitem{Girfoglio_ROM_Fluids}
M.~Girfoglio, A.~Quaini, and G.Rozza.
\newblock Pressure stabilization strategies for a {LES filtering Reduced Order
  Model}.
\newblock https://arxiv.org/abs/2106.15887, 2021.

\bibitem{Girfoglio2019}
M.~Girfoglio, A.~Quaini, and G.~Rozza.
\newblock {A Finite Volume approximation of the Navier-Stokes equations with
  nonlinear filtering stabilization}.
\newblock {\em Computers \& Fluids}, 187:27--45, 2019.

\bibitem{Girfoglio_JCP}
M.~Girfoglio, A.~Quaini, and G.~Rozza.
\newblock {A POD-Galerkin reduced order model for a LES filtering approach}.
\newblock {\em Journal of Computational Physics}, 436:110260, 2021.

\bibitem{Gunzburger2019}
M.~Gunzburger, T.~Iliescu, M.~Mohebujjaman, and M.~Schneier.
\newblock An evolve-filter-relax stabilized reduced order stochastic
  collocation method for the time-dependent {Navier--Stokes} equations.
\newblock {\em SIAM/ASA Journal on Uncertainty Quantification}, 7:1162--1184,
  2019.

\bibitem{hesthaven2015certified}
J.~S. Hesthaven, G.~Rozza, and B.~Stamm.
\newblock {\em {Certified Reduced Basis Methods for Parametrized Partial
  Differential Equations}}.
\newblock Springer International Publishing, 2016.

\bibitem{Hijazi2020}
S.~Hijazi, G.~Stabile, A.~Mola, and G.~Rozza.
\newblock Data-driven {POD-Galerkin} reduced order model for turbulent flows.
\newblock {\em Journal of Computational Physics}, 416:109513, 2020.

\bibitem{O-hunt1988}
J.C. Hunt, A.A. Wray, and P.~Moin.
\newblock Eddies stream and convergence zones in turbulent flows.
\newblock Technical Report CTR-S88, CTR report, 1988.

\bibitem{PISO}
R.~I. Issa.
\newblock Solution of the implicitly discretised fluid flow equations by
  operator-splitting.
\newblock {\em Journal of Computational Physics}, 62(1):40--65, 1986.

\bibitem{John2004}
V.~John.
\newblock {Reference values for drag and lift of a two dimensional
  time-dependent flow around a cylinder}.
\newblock {\em International Journal for Numerical Methods in Fluids},
  44:777--788, 2004.

\bibitem{John2006}
V.~John.
\newblock On the efficiency of linearization schemes and coupled multigrid
  methods in the simulation of a {3D} flow around a cylinder.
\newblock {\em International Journal for Numerical Methods in Fluids},
  50:845--862, 2006.

\bibitem{JOHNSTON2004221}
H.~Johnston and J.-G. Liu.
\newblock {Accurate, stable and efficient {Navier--Stokes} solvers based on
  explicit treatment of the pressure term}.
\newblock {\em Journal of Computational Physics}, 199(1):221--259, 2004.

\bibitem{Kalashnikova_ROMcomprohtua}
I.~Kalashnikova and M.~F. Barone.
\newblock {On the stability and convergence of a Galerkin reduced order model
  (ROM) of compressible flow with solid wall and far-field boundary treatment}.
\newblock {\em International Journal for Numerical Methods in Engineering},
  83(10):1345--1375, 2010.

\bibitem{Kunisch2002492}
K.~Kunisch and S.~Volkwein.
\newblock {Galerkin proper orthogonal decomposition methods for a general
  equation in fluid dynamics}.
\newblock {\em SIAM Journal on Numerical Analysis}, 40(2):492--515, 2002.

\bibitem{Lax1960}
P.D. Lax and B.~Wendroff.
\newblock System of conservation laws.
\newblock {\em Communications on Pure and Applied Mathematics}, 13:217--237,
  1960.

\bibitem{layton_CMAME}
W.~Layton, L.G. Rebholz, and C.~Trenchea.
\newblock Modular nonlinear filter stabilization of methods for higher
  {Reynolds} numbers flow.
\newblock {\em Journal of Mathematical Fluid Mechanics}, 14:325--354, 2012.

\bibitem{Lazzaro2002}
D.~Lazzaro and L.~Montefusco.
\newblock Radial basis functions for the multivariate interpolation of large
  scattered data sets.
\newblock {\em Journal of Computational and Applied Mathematics}, 140:521--536,
  2002.

\bibitem{Leray1934}
J.~Leray.
\newblock Essai sur le mouvement d'un fluide visqueux emplissant l'espace.
\newblock {\em Acta Mathematica}, 63:193--248, 1934.

\bibitem{Lorenzi2016}
S.~Lorenzi, A.~Cammi, L.~Luzzi, and G.~Rozza.
\newblock {POD-Galerkin method for finite volume approximation of Navier-Stokes
  and RANS equations}.
\newblock {\em Computer Methods in Applied Mechanics and Engineering},
  311:151--179, 2016.

\bibitem{Olshanskii2013}
M.A. Olshanskii and X.~Xiong.
\newblock A connection between filter stabilization and eddy viscosity models.
\newblock {\em Numerical Methods for Partial Differential Equations},
  29(6):2061--2080, 2013.

\bibitem{Orszag1986}
S.~A. Orszag, M.~Israeli, and M.O. Deville.
\newblock {Boundary conditions for incompressible flows}.
\newblock {\em Journal of Scientific Computing}, 1(1):75--111, 1986.

\bibitem{Passerini2013}
T.~Passerini, A.~Quaini, U.~Villa, A.~Veneziani, and S.~Canic.
\newblock Validation of an open source framework for the simulation of blood
  flow in rigid and deformable vessels.
\newblock {\em International Journal for Numerical Methods in Biomedical
  Engineering}, 29(11):1192--1213, 2013.

\bibitem{pope}
S.B Pope.
\newblock {\em Turbulent flows}.
\newblock Cambridge University Press, Cambridge, 2000.

\bibitem{quarteroniRB2016}
A.~Quarteroni, A.~Manzoni, and F.~Negri.
\newblock {\em {Reduced Basis Methods for Partial Differential Equations}}.
\newblock Springer International Publishing, 2016.

\bibitem{Rozza2008}
G.~Rozza, D.~B.~P. Huynh, and A.~T. Patera.
\newblock {Reduced Basis Approximation and a Posteriori Error Estimation for
  Affinely Parametrized Elliptic Coercive Partial Differential Equations}.
\newblock {\em Archives of Computational Methods in Engineering}, 15(3):229,
  2008.

\bibitem{RoSta17}
G.~Stabile and G.~Rozza.
\newblock {ITHACA-FV - In real Time Highly Advanced Computational Applications
  for Finite Volumes}.
\newblock Accessed: 2018-01-30.

\bibitem{Stabile2018}
G.~Stabile and G.~Rozza.
\newblock {Finite volume POD-Galerkin stabilised reduced order methods for the
  parametrised incompressible Navier–Stokes equations}.
\newblock {\em Computer \& Fluids}, 173:273--284, 2018.

\bibitem{Strazzullo2021}
M.~Strazzullo, F.~Ballarin, M.~Girfoglio, T.~Iliescu, and G.~Rozza.
\newblock Evolve-filter-relax based reduced order models for convection
  dominated flows in the finite element fashion.
\newblock {\em in preparation}.

\bibitem{Huerta2020}
V.~Tsiolakis, M.~Giacomini, R.~Sevilla, C.~Othmer, and A.~Huerta.
\newblock Parametric solutions of turbulent incompressible flows in openfoam
  via the proper generalised decomposition.
\newblock https://arxiv.org/abs/2006.07073, 2020.

\bibitem{turek1996}
S.~Turek and M.~Sch\"afer.
\newblock Benchmark computations of laminar flow around cylinder.
\newblock In E.H. Hirschel, editor, {\em Flow Simulation with High-Performance
  Computers II}, volume~52 of {\em Notes on Numerical Fluid Mechanics}. Vieweg,
  1996.

\bibitem{Vreman2004}
A.W. Vreman.
\newblock An eddy-viscosity subgrid-scale model for turbulent shear flow:
  {Algebraic} theory and applications.
\newblock {\em Physics of Fluids}, 16(10):3670--3681, 2004.

\bibitem{wang_turb}
Z.~Wang, I.~Akhtar, J.~Borggaard, and T.~Iliescu.
\newblock {Proper orthogonal decomposition closure models for turbulent flows:
  A numerical comparison}.
\newblock {\em Computer Methods in Applied Mechanics and Engineering},
  237–240:10--26, 2012.

\bibitem{Weller1998}
H.~G. Weller, G.~Tabor, H.~Jasak, and C.~Fureby.
\newblock A tensorial approach to computational continuum mechanics using
  object-oriented techniques.
\newblock {\em Computers in physics}, 12(6):620--631, 1998.

\bibitem{Wells2017}
D.~Wells, Z.~Wang, X.~Xie, and T.~Iliescu.
\newblock An evolve-then-filter regularized reduced order model for
  convection-dominated flows.
\newblock {\em International Journal for Numerical Methods in Fluids},
  84:598--615, 2017.

\bibitem{Xie2018_2}
X.~Xie, F.~Bao, and C.~Webster.
\newblock Evolve filter stabilization reduced-order model for stochastic
  burgers equation.
\newblock {\em Fluids}, 3:84, 2018.

\end{thebibliography}

\end{document}